\documentclass[12pt,a4paper,leqno]{article}
\usepackage{setspace} \usepackage{etoolbox} \usepackage{comment}
\usepackage{amsfonts}
\usepackage[margin=2cm]{geometry}
\usepackage{pdflscape}
\usepackage{times}
\usepackage{amsthm,amssymb}
\usepackage{amsmath,xypic}
\usepackage[cal=boondoxo,calscaled=.96,scr=rsfs]{mathalpha} \DeclareMathAlphabet{\mathscrbf}{OMS}{mdugm}{b}{n} \usepackage{multicol}
\usepackage{array}
\usepackage{hyphenat}
\usepackage[usenames,dvipsnames]{color}
\usepackage{rotating}
\usepackage{longtable}
\usepackage{caption}
\usepackage{float}
\usepackage{tcolorbox}

\usepackage{epigraph}
\usepackage{enumitem}
\setlist[enumerate]{listparindent=0.5in}
\newcommand{\be}{\begin{equation}}
\newcommand{\ee}{\end{equation}}
\newcommand{\bes}{\begin{equation*}}
\newcommand{\ees}{\end{equation*}}
\newcommand{\bea}{\begin{eqnarray}}
\newcommand{\eea}{\end{eqnarray}}
\newcommand{\beas}{\begin{eqnarray}}
\newcommand{\eeas}{\end{eqnarray}}
\newcommand{\ben}{\begin{note}}
\newcommand{\een}{\end{note}}
\newcommand{\bexl}{\vskip0.1em\noindent\hrulefill\vskip1em\begin{ExerciseList}}
\newcommand{\eexl}{\end{ExerciseList}\hrulefill}

\newcommand{\bthm}{\begin{theorem}}
\newcommand{\ethm}{\end{theorem}}
\newcommand{\bpro}{\begin{prop}}
\newcommand{\epro}{\end{prop}}
\newcommand{\bcor}{\begin{corollary}}
\newcommand{\ecor}{\end{corollary}}
\newcommand{\bcon}{\begin{conjecture}}
\newcommand{\econ}{\end{conjecture}}
\newcommand{\bp}{\begin{proof}}
\newcommand{\ep}{\end{proof}}
\newcommand{\blem}{\begin{lemma}}
\newcommand{\elem}{\end{lemma}}
\newcommand{\bn}{\begin{note}}
\newcommand{\en}{\end{note}}
\newcommand{\benum}{\begin{enumerate}}
\newcommand{\eenum}{\end{enumerate}}
\newcommand{\bed}{\begin{defn}}
\newcommand{\eed}{\end{defn}}
\newcommand{\brem}{\begin{remark}}
\newcommand{\erem}{\end{remark}}

\newcommand{\btik}{\begin{tikzpicture}\begin{axis}[scale=0.5,axis y line=center, axis x line=middle]}
\newcommand{\etik}{\end{axis}\end{tikzpicture}}

\let\into=\hookrightarrow
\let\mapsto=\longmapsto
\let\cong=\equiv

\newcommand{\upperRomannumeral}[1]{\uppercase\expandafter{\romannumeral#1}}

 \usepackage{tikz-cd}
\usepackage{graphicx}
\usepackage[pagewise]{lineno}
\usepackage{stackengine}
\usepackage{stmaryrd}
\usepackage{titlesec}
\newtoggle{arxiv}
\toggletrue{arxiv}
\iftoggle{arxiv}{
	\usepackage{url}
\usepackage{natbib}
\let\cite=\citep
	\typeout{In the arxiv mode}
	\usepackage{fancyhdr}
	\usepackage[colorlinks,citecolor=blue]{hyperref}
\usepackage{zref-clever}
	\let\Cref=\zcref
	\let\cref=\Cref
	\let\cite=\citep
}
{
}
\newtoggle{draft}
\togglefalse{draft}

\usepackage{xr}
\newcommand{\constranab}[2]{\cite[{#1\ref{0-#2}}]{joshi-anabelomorphy}}

\newcommand{\constrtwoh}[2]{\cite[{#1\ref{II5-#2}}]{joshi-teich-def}}

\newcommand{\constrthrB}[2]{\cite[{#1\ref{IIIB-#2}}]{joshi-teich-IIIB-rosetta}}

\newcommand{\constrfour}[2]{\cite[{#1\ref{IV-#2}}]{joshi-teich-abc}}

\newtheorem{theorem}[equation]{Theorem}      \newtheorem{lemma}[equation]{Lemma}          \newtheorem{corollary}[equation]{Corollary}  \newtheorem{proposition}[equation]{Proposition}

\theoremstyle{definition}
\newtheorem{conj}[equation]{Conjecture}

\newtheorem{question}[equation]{Question}

\theoremstyle{definition}
\newtheorem{defn}[equation]{Definition}
\theoremstyle{remark}

\theoremstyle{definition}
\newtheorem{remark}[equation]{Remark}

\numberwithin{equation}{section}

\zcsetup{cap} \AddToHook{env/theorem/begin}{\zcsetup{countertype={equation=theorem}}}
\AddToHook{env/corollary/begin}{\zcsetup{countertype={equation=corollary}}}
\AddToHook{env/lemma/begin}{\zcsetup{countertype={equation=lemma}}}
\AddToHook{env/proposition/begin}{\zcsetup{countertype={equation=proposition}}}
\AddToHook{env/prop/begin}{\zcsetup{countertype={equation=prop}}}
\AddToHook{env/theoremdef/begin}{\zcsetup{countertype={equation=theoremdef}}}
\AddToHook{env/defn/begin}{\zcsetup{countertype={equation=defn}}}
\AddToHook{env/definition/begin}{\zcsetup{countertype={equation=definition}}}
\AddToHook{env/remark/begin}{\zcsetup{countertype={equation=remark}}}
\AddToHook{env/rem/begin}{\zcsetup{countertype={equation=rem}}}
\AddToHook{env/conj/begin}{\zcsetup{countertype={equation=conj}}}
\AddToHook{env/conjecture/begin}{\zcsetup{countertype={equation=conjecture}}}
\AddToHook{env/question/begin}{\zcsetup{countertype={equation=question}}}
\AddToHook{env/example/begin}{\zcsetup{countertype={equation=example}}}

\zcRefTypeSetup{section}{
	Name-sg = \S ,
	name-sg = \S ,
	Name-pl = \S\S ,
	name-pl = \S\S ,
}

\zcRefTypeSetup{subsection}{
	Name-sg = \S ,
	name-sg = \S ,
	Name-pl = \S\S ,
	name-pl = \S\S ,
}

\zcRefTypeSetup{subsubsection}{
	Name-sg = \S ,
	name-sg = \S ,
	Name-pl = \S\S ,
	name-pl = \S\S ,
}

\zcRefTypeSetup{equation}{
	Name-sg = eq. ,
	name-sg = eq. ,
	Name-pl = eqns. ,
	name-pl = eqns. ,
}

\zcRefTypeSetup{definition}{
	Name-sg = Definition ,
	name-sg = definition ,
	Name-pl = Definitions ,
	name-pl = definitions ,
}

\zcRefTypeSetup{defn}{
	Name-sg = Definition ,
	name-sg = definition ,
	Name-pl = Definitions ,
	name-pl = definitions ,
}

\zcRefTypeSetup{theorem}{
	Name-sg = Theorem ,
	name-sg = theorem ,
	Name-pl = Theorems ,
	name-pl = theorems ,
}

\zcRefTypeSetup{theoremdef}{
	Name-sg = Theorem-Definition ,
	name-sg = theorem-definition ,
	Name-pl = Theorem-Definitions ,
	name-pl = theorem-definitions ,
}

\zcRefTypeSetup{proposition}{
	Name-sg = Proposition ,
	name-sg = proposition ,
	Name-pl = Propositions ,
	name-pl = propositions ,
}

\zcRefTypeSetup{lemma}{
	Name-sg = Lemma ,
	name-sg = lemma ,
	Name-pl = Lemmas ,
	name-pl = lemmas ,
}

\zcRefTypeSetup{remark}{
	Name-sg = Remark ,
	name-sg = remark ,
	Name-pl = Remarks ,
	name-pl = remarks ,
}

\zcRefTypeSetup{question}{
	Name-sg = Question ,
	name-sg = question ,
	Name-pl = Questions ,
	name-pl = questions ,
}

\zcRefTypeSetup{example}{
	Name-sg = Example ,
	name-sg = example ,
	Name-pl = Examples ,
	name-pl = examples ,
}

\zcRefTypeSetup{conj}{
	Name-sg = Conjecture,
	name-sg = conjecture ,
	Name-pl = Conjectures ,
	name-pl = conjectures ,
}

\zcRefTypeSetup{conjecture}{
	Name-sg = Conjecture,
	name-sg = conjecture ,
	Name-pl = Conjectures ,
	name-pl = conjectures ,
}

\newcommand{\para}{\subsection{}}
\titleformat{\subsection}[runin]{\normalfont\bfseries}{\S\ \thesubsection}{.5em}{}[{\ \ }]
\titlespacing{\subsection}{0pt}{1.5ex plus .1ex minus .2ex}{0pt}
\titleformat{\subsubsection}[runin]{\normalfont\bfseries}{\S\ \thesubsubsection}{.5em}{}[{\ \ }]
\titlespacing{\subsubsection}{0pt}{0.75ex plus .1ex minus .2ex}{0pt}

\newcommand{\nws}{\numberwithin{equation}{section}}
\newcommand{\nwss}{\numberwithin{equation}{subsection}}

\let\into=\hookrightarrow
\let\isom=\simeq

\let\tensor=\otimes
\newcommand{\A}{\mathscr{A}}
\newcommand{\abs}[1]{\left\vert#1\right\vert}

\newcommand{\bF}{{\bar{F}}}
\newcommand{\bQ}{{\bar{\Q}}}

\newcommand{\C}{{\mathbb C}}

\newcommand{\End}{\rm{End}}

\newcommand{\F}{{\mathbb F}}

\newcommand{\gal}{{\rm Gal}}

\newcommand{\hGm}{\widehat{\mathbb{G}}_m}

\newcommand{\mydot}{{\small{\bullet}}}
\newcommand{\N}{\mathscr{N}}

\newcommand{\Q}{{\mathbb Q}}
\newcommand{\R}{{\mathbb R}}

\newcommand{\spec}{{\rm Spec}}
\newcommand{\Spec}{{\rm Spec}}

\newcommand{\Z}{{\mathbb Z}}

\renewcommand{\O}{{\mathscr O}}
\renewcommand{\P}{{\mathbb P}}
\renewcommand{\wp}{{\mathfrak p}}

\newcommand{\fm}{{\mathfrak{M}}}

\newcommand{\invlim}{\varprojlim}
\let\fm=\fa
\renewcommand{\fm}{\mathfrak{m}}
\newcommand{\mapright}[1]{{\xymatrix{{}\ar[r]^{#1}&{}}}}

\renewcommand{\bpro}{\begin{proposition}}
	\renewcommand{\epro}{\end{proposition}}
\renewcommand{\bcon}{\begin{conj}}
	\renewcommand{\econ}{\end{conj}}
\newcommand{\ilim}{\varprojlim}
\title{Construction of Arithmetic Teichm\"uller Spaces I
}
\author{Kirti Joshi}

\newcommand{\Address}{\bigskip\noindent{\footnotesize\textsc{{Math. department, University of Arizona, 617 N Santa Rita, Tucson
		85721-0089, USA.}}\par\nopagebreak 
\noindent\textit{Email:}	\texttt{kirti@arizona.edu}}}

\begin{document}
\maketitle

\lhead{}

\iftoggle{draft}{\pagewiselinenumbers}{\relax}
\newcommand{\act}{\curvearrowright}
\newcommand{\lmp}{{\Pi\act\Ot}}
\newcommand{\lmpi}{{\lmp}_{\int}}
\newcommand{\lmpf}{\lmp_F}
\newcommand{\Om}{\O^{\times\mu}}
\newcommand{\Omf}{\O^{\times\mu}_{\bF}}
\renewcommand{\N}{\mathbb{N}}
\newcommand{\yoga}{Yoga}
\newcommand{\gl}[1]{{\rm GL}(#1)}
\newcommand{\bK}{\overline{K}}
\newcommand{\reptrip}{\rho:G_K\to\gl V}
\newcommand{\reptripp}[1]{\rho\circ\alpha:G_{\ifstrempty{#1}{K}{{#1}}}\to\gl V}
\newcommand{\benumlab}{\begin{enumerate}[label={{\bf(\arabic{*})}}]}
\newcommand{\ord}{\mathop{\rm ord}\nolimits}	
\newcommand{\kcs}{K^\circledast}
\newcommand{\lcs}{L^\circledast}
\renewcommand{\A}{\mathbb{A}}
\newcommand{\bfq}{\bar{\mathbb{F}}_q}
\newcommand{\tripod}{\P^1-\{0,1728,\infty\}}

\newcommand{\vseq}[2]{{#1}_1,\ldots,{#1}_{#2}}
\newcommand{\anab}[4]{\left({#1},\{#3 \}\right)\anabelmap\left({#2},\{#4 \}\right)}

\newcommand{\gln}{{\rm GL}_n}
\newcommand{\glo}[1]{{\rm GL}_1(#1)}
\newcommand{\glt}[1]{{\rm GL_2}(#1)}

\newcommand{\iut}{\cite{mochizuki-iut1, mochizuki-iut2, mochizuki-iut3,mochizuki-iut4}}
\renewcommand{\iut}{\cite{mochizuki-iut1-4}}
\newcommand{\iutabc}{\cite{mochizuki-iut1, mochizuki-iut2, mochizuki-iut3}}
\newcommand{\topics}{\cite{mochizuki-topics1,mochizuki-topics2,mochizuki-topics3}}
\newcommand{\present}{the present series of papers (\cite{joshi-anabelomorphy, joshi-teich, joshi-untilts,joshi-teich-estimates,joshi-teich-def})}
\newcommand{\Present}{The present series of papers (\cite{joshi-anabelomorphy, joshi-teich, joshi-untilts,joshi-teich-estimates,joshi-teich-def})}
\newcommand{\moccor}{\iutf{IUT3, Corollary 3.12}}

\newcommand{\linv}{\mathfrak{L}}
\newcommand{\bedef}{\begin{defn}}
\newcommand{\eedef}{\end{defn}}
\renewcommand{\act}[1][]{\overset{#1}{\curvearrowright}}
\newcommand{\bfx}{\overline{F(X)}}
\newcommand{\anabelmap}{\leftrightsquigarrow}
\newcommand{\ban}[1][G]{\mathscr{B}({#1})}
\newcommand{\pit}{\Pi^{temp}}
 
 \newcommand{\bL}{\overline{L}}
 \newcommand{\bkm}{\bK_M}
 \newcommand{\vbk}{v_{\bK}}
 \newcommand{\vbkm}{v_{\bkm}}
\newcommand{\ocs}{\O^\circledast}
\newcommand{\ot}{\O^\triangleright}
\newcommand{\ocsk}{\ocs_K}
\newcommand{\otk}{\ot_K}
\newcommand{\ok}{\O_K}
\newcommand{\oko}{\O_K^1}
\newcommand{\oks}{\ok^*}
\newcommand{\Qpb}{\overline{\Q}_p}
\newcommand{\Qpbh}{\widehat{\overline{\Q}}_p}
\newcommand{\tr}{\triangleright}
\newcommand{\ocpt}{\O_{\C_p}^\tr}
\newcommand{\ocpf}{\O_{\C_p}^\flat}
\newcommand{\sG}{\mathscr{G}}
\newcommand{\sX}{\mathscr{X}}
\newcommand{\sY}{\mathscr{Y}}
\newcommand{\sxfe}{\mathscr{X}_{F,E}}
\newcommand{\sxfep}{\mathscr{X}_{F,E'}}
\newcommand{\loglt}{\log_{\sG}}
\newcommand{\fc}{\mathfrak{t}}
\newcommand{\ku}{K_u}
\newcommand{\kup}{\ku'}
\newcommand{\kt}{\tilde{K}}
\newcommand{\sGpf}{\sG(\O_K)^{pf}}
\newcommand{\hgm}{\widehat{\mathbb{G}}_m}
\newcommand{\bE}{\overline{E}}
\newcommand{\syfe}{\mathscr{Y}_{F,E}}
\newcommand{\syfqp}[1]{\mathscr{Y}_{\cptl{#1},\Q_p}}
\newcommand{\syfqpe}[1]{\mathscr{Y}_{{#1},E}}
\newcommand{\fJ}{\mathfrak{J}}
\newcommand{\fM}{\mathfrak{M}}
\newcommand{\locvar}{local arithmetic-geometric anabelian variation of fundamental group of $X/E$ at $\wp$}
\newcommand{\fjxep}{\fJ(X,E,\wp)}
\newcommand{\fjxe}{\fJ(X,E)}
\newcommand{\ufjxe}{\underline{\fJ}(X,E)}
\newcommand{\fpc}[1]{\widehat{{\overline{\F_p(({#1}))}}}}
\newcommand{\cpt}{\C_p^\flat}
\newcommand{\cptl}[1]{\C_{p,{#1}}^\flat}
\newcommand{\fja}[1]{\fJ^{\rm arith}({#1})}
\newcommand{\ainfe}{A_{\inf,E}(\O_F)}
\renewcommand{\ainfe}{W_{\O_E}(\O_F)}
\newcommand{\gmh}{\widehat{\mathbb{G}}_m}
\newcommand{\sE}{\mathscr{E}}
\newcommand{\bpi}{B^{\varphi=\pi}}
\newcommand{\onto}{\twoheadrightarrow}

\newcommand{\cpmax}{{\C_p^{\rm max}}}
\newcommand{\xan}{X^{an}}
\newcommand{\yan}{Y^{an}}
\newcommand{\bPi}{\overline{\Pi}}
\newcommand{\bPit}{\bPi^{\rm{\scriptscriptstyle temp}}}
\newcommand{\Pit}{\Pi^{\rm{\scriptscriptstyle temp}}}
\renewcommand{\pit}[1]{\Pi^{\scriptscriptstyle temp}_{#1}}
\newcommand{\pitb}[1]{\overline{\Pi}^{\scriptscriptstyle temp}_{#1}}
\newcommand{\pitk}[2]{\Pi^{\scriptscriptstyle temp}_{#1;#2}}
\newcommand{\pio}[1]{\pi_1({#1})}
\newcommand{\fTeich}{\widetilde{\fJ(X/L)}}
\newcommand{\ssep}{\S\,} 

\setcounter{tocdepth}{1}

\tableofcontents
\togglefalse{draft}
\newcommand{\FF}{\cite{fargues-fontaine}}
\iftoggle{draft}{\pagewiselinenumbers}{\relax}

\newcommand{\attportion}{Sections~\ref{se:number-field-case}, \ref{se:construct-att}, \ref{se:relation-to-iut}, \ref{se:self-similarity} and \ref{se:applications-elliptic}}

\newcommand{\fgm}{\gmh}
\newcommand{\fgmcp}{\fgm(\O_{\C_p})}
\newcommand{\sgtok}{\widetilde{\sG(\O_K)}}
\newcommand{\sgtocp}{\widetilde{\sG(\O_{\C_p})}}
\newcommand{\ttxl}{\tilde{\Theta}_{X,\ell}}
\newcommand{\ttxlt}[1]{\tilde{\Theta}_{X,\ell,{#1}}}
\newcommand{\bpip}{B^{\varphi=p}}
\newcommand{\syQp}{\mathscr{Y}_{\cpt,\Q_p}}
\newcommand{\sxQp}{\mathscr{X}_{\cpt,\Q_p}}
\newcommand{\sgt}{\widetilde{\sG}}
\newcommand{\otmu}{\O_{\bE}^{\times\mu}}

\newcommand{\V}{\mathbb{V}}
\newcommand{\vl}{\V_L}
\newcommand{\vnl}{\V_L^{non}}
\newcommand{\vlnon}{\vnl}
\newcommand{\vlarch}{\mathbb{V}_L^{arc}}
\newcommand{\sL}{\mathscr{L}}
\newcommand{\sM}{\mathscr{M}}
\newcommand{\fjxlv}{\fJ(X,L_v)}
\newcommand{\tfjxl}{\widetilde{\fJ(X,L)}}

\newcommand{\iutthr}{\cite{mochizuki-iut1,mochizuki-iut2,mochizuki-iut3}}
\newcommand{\iutf}[1]{\cite[#1]{mochizuki-iut1-4}}
\renewcommand{\iutthr}{\iutf{IUT1-IUT3}}
\renewcommand{\iut}{\cite{mochizuki-iut1-4}}
\newcommand{\four}{Sections~\ref{se:grothendieck-conj}, \ref{se:untilts-of-Pi}, and \ref{se:riemann-surfaces}}
\newpage
\epigraphwidth0.55\textwidth
\epigraph{
	Then let winged Fancy wander\\
	Through the thought still spread beyond her:\\
	Open wide the mind’s cage-door,\\
	She’ll dart forth, and cloudward soar.\\
	O sweet Fancy! let her loose;}{\iftoggle{arxiv}{John Keats \cite{keats}}{\citeauthor{keats}}}
	
\newcommand{\PSL}{{\mathrm{PSL}}}
\newcommand{\fH}{\mathfrak{H}}
\newcommand{\sT}{\mathcal{T}}
 \interfootnotelinepenalty=10000 \addtocontents{toc}{\protect\setcounter{tocdepth}{3}}
\section{Introduction}

\subsection{Extending the analogy between Number Fields and Riemann surfaces}

As is well-known, from \cite{dedekind-weber}, one has an analogy between Riemann surfaces and number fields in which the field,  $\C(\Sigma)$,  of meromorphic functions on a connected, compact Riemann surface $\Sigma$ corresponds to a number field $L$:
$$
\C(\Sigma)\leftrightsquigarrow L.
$$ 
This analogy has played a central role in the evolution of Number Theory and Algebraic Geometry over the past hundred and fifty or so years since its inauguration. 

Classical Teichm\"uller Theory\footnote{Standard texts are \cite{ahlfors-quasiconf-book}, \cite{lehto-book}, \cite{nag-book}, \cite{imayoshi-book}.} extends Riemann's moduli Theory and provides a refined view of Riemann surfaces.  Therefore, it is not unreasonable to ask if a  Teichm\"uller Theory  exists for  number fields.

That such a Teichm\"uller theory exists was first asserted, but with far lesser quantitative precision, by Shinichi Mochizuki in \iutthr. Earliest indication, independent of Mochizuki's work, that such a theory does indeed exist came from \cite{joshi-untilts-2020} which approached the problem from an  analytic function theory \textit{optik}. 
This paper and \cite{joshi-teich-def} lays the foundations of the Theory of Arithmetic Teichm\"uller Spaces of Algebraic Varieties over Number Fields. This paper primarily deals with varieties over $p$-adic fields and number fields, while \cite{joshi-teich-def} deals exclusively with the Arithmetic Teichm\"uller Theory of Number Fields. Other papers: \cite{joshi-teich-estimates,joshi-teich-IIIB-rosetta,joshi-teich-IIIC-moccor,joshi-teich-abc-conj} deal with the principal application of this theory envisaged by Mochizuki in \iut.

My approach to establishing the existence of the arithmetic theory is based on the strong analogy with the classical theory of Riemann surfaces. For example, one knows that there exists Riemann surfaces of fixed topological type $(g,n)$ with isomorphic topological fundamental groups.  \Cref{thm:main2} and \Cref{thm:main} establish the $p$-adic analog of this well-known classical result by replacing complex holomorphic structures with Berkovich analytic spaces in such a way, that Berkovich analytic structures are deformed while keeping both: the (algebraic and the $p$-adic analytic) datum over finite extensions of $\Q_p$ fixed, and the \'etale fundamental group and the geometric \'etale fundamental group datum fixed (i.e remain isomorphic as topological groups). This approach provides a firm foundation for defining Arithmetic Holomorphic Structures in \Cref{def:arith-hol-strs}, and provides full compatibility with Anabelian Geometry (as the relevant fundamental groups remain fixed).

\clearpage
\subsection{Tables of Analogies between Classical and Arithmetic Teichm\"uller Theories}
Readers may find the following tables of analogies between the geometric and the arithmetic cases  useful in comparing the theorems and relevant concepts.
\vskip0.25in
\begin{tcolorbox}[label=tab:analogies1,grow to right by = 1cm, colback=white,coltitle=black,colframe=white!25!olive,fonttitle=\bfseries,title={Table of Analogies: Riemann Surfaces and Number Fields}]
\centering
		\begin{tabular}{>{\raggedright}p{1.5in}|>{\raggedright}p{1.5in}|>{\raggedright}p{1.5in}|>{\raggedright}p{1.25in}p{0mm}}
			\hline
			Riemann's Moduli Theory	&	Isomorphism class of a Riemann Surface (i.e. a point of Riemann's moduli space) \cite[Chapter 1]{imayoshi-book}	& Isomorphism class of a Number Field  & Dedekind-Weber Analogy \cite{dedekind-weber} & \\ \hline 
			Teichm\"uller Theory &		Riemann surface + a quasi-conformal mapping (i.e. a point of Teichm\"uller space) \cite[Chapter 1]{imayoshi-book} & Number field + extra arithmetic-topological data (precisely codified in  \constrtwoh{Definition }{def:arithmeticoid})
			[\iut\ rests on the existence of such arithmetic Teichm\"uller data.] & Extension of the Dedekind-Weber Analogy due to Joshi-Mochizuki \cite{joshi-teich-def} &\\ 
			\hline
		\end{tabular}
\end{tcolorbox}

\newcommand{\syflv}{\sY_{\lvbht,L_v}}
\newcommand{\syflvmax}{\sY_{\lvbhtmax,L_v}}
\newcommand{\sxflv}{\sX_{\lvbht,L_v}}
\newcommand{\sxflvmax}{\sX_{\lvbhtmax,L_v}}
\newcommand{\bsY}{\mathscrbf{Y}}
\newcommand{\bsX}{\mathscrbf{X}}
\newcommand{\yadl}{\bsY_L}

\begin{tcolorbox}[label=tab:analogies2,grow to right by = 1cm, colback=white,coltitle=black,colframe=white!25!olive,fonttitle=\bfseries,title={Table of Analogies 2: Classical and Arithmetic Teichm\"uller Theory}]
\centering
		\begin{tabular}{>{\raggedright}p{2.5in}|>{\raggedright}p{2.5in}p{0mm}}
			\hline
			{Classical Teichm\"uller Theory}	& {Arithmetic Teichm\"uller Theory} & \\ \hline
{Approach of \cite{nakai59} via Banach space structure of Royden Algebras of Riemann Surfaces} &
			{(This paper) Banach spaces of Berkovich $p$-adic analytic functions } \Cref{thm:main,thm:main2} & \\	\hline
Virasoro action \cite{beilinson88} & \Cref{th:main3}{\bf(5)} Arithmetic (analog of) Virasoro action & \\ \hline
Teichm\"uller (metric) Space $\sT_\Sigma$ of a Riemann Surface $\Sigma$ & Arithmetic Teichm\"uller Space $\yadl$ of a Number Field $L$ (the metrisable topological space $\yadl$ is constructed and described in detail in \cite{joshi-teich-def}) &\\
			\hline
				Each point of $\sT_{\Sigma}$ corresponds to a deformation of the complex structure of $\Sigma$ using quasi-conformal mappings& Each point of $\yadl$ corresponds to a topological deformation of the arithmetic of $L$ \constrtwoh{Theorem }{th:consequence-inequivalent-arithmeticoids} & \\
				\hline
			Every point of $\sT_{\Sigma}$ provides an isomorph of $\pi_1^{top}(\Sigma)$ & Every point of $\yadl$ provides an isomorph of $\pi_1^{et}(\Spec(L))=\gal(\bL/L)$ & \\ \hline
		\end{tabular}
	\end{tcolorbox}

\clearpage
\subsection{The goals of Arithmetic Teichm\"uller Theory}\label{ss:gatt}
The goal of this paper (Constructions I) and \cite{joshi-teich-def} (Constructions I(1/2)) is to provide for each number field $L$ (assumed for simplicity to have no real embeddings) and for each geometrically connected, smooth, quasi-projective variety $X/L$  over a number field $L$, and  for each place $v\in\vl$, the following
\benumlab
\item 
a $v$-adic analog of the Classical Teichm\"uller Space, more precisely a category $\fjxlv$ (\Cref{def:arith-hol-space-local} for $v\in\vlnon$ and \Cref{def:arith-hol-space-local-arch} for $v\in\vlarch$; in dimension one, for $v\in\vlarch$, one obtains classical Teichm\"uller spaces), of $X/L_v$, such that objects of $\fjxlv$ can be distinguished from each other (by means of their \textit{Arithmetic Holomorphic Structures} (\Cref{def:arith-hol-strs}) and in general the $v$-adic Berkovich analytic space structure of functions changes in a manner similar to classical Teichm\"uller space (\Cref{th:variation-of-arith-hol-strs}); for $v\in\vlarch$, the Banach structure of rings of functions varies (\Cref{th:nakai1}, \Cref{th:nakai2}). But all objects of $\fjxlv$ have isomorphic tempered fundamental groups (\Cref{pr:fund-group-same} for $v\in\vlnon$ and \Cref{th:nakai2} for $v\in \vlarch$) and hence these, and other results proved here, justify the moniker \textit{`Arithmetic Teichm\"uller space of $X/L_v$.'}
\item The construction of {\bf(1)} provides an Adelic Arithmetic Teichm\"uller Space i.e. an adelic category $\tfjxl$ (\Cref{th:main4}) built from the local Arithmetic Teichm\"uller Spaces $\fjxlv$ (for each $v\in\vl$) together with global compatibilities, forced by the constructions of {\bf(3)}, which justify the moniker \textit{`Adelic Arithmetic Teichm\"uller space of $X/L$.'}
\item As a special case, taking $X=\Spec(L)$, the construction of {\bf(1,2)} provides a topological space, with a certain list of properties, which justify the moniker \textit{`Arithmetic Teichm\"uller Space of the number field $L$.'} [This construction is detailed in  \cite{joshi-teich-def}.]
\eenum

In the texts on Teichm\"uller Theory, Classical Teichm\"uller Spaces  are constructed by considering quasi-conformal mappings  between Riemann surfaces instead of conformal mappings. A key property of quasi-conformal mappings is that they are not necessarily holomomorphic mappings, and secondly under such mappings, locally i.e. in a neighborhood of a point, metrics are scaled or distorted (\cite{ahlfors-quasiconf-book}). In the context of ({\bf(1), (2)}) above, there is an analogous notion of $v$-adic quasi-conformal mappings and under such mappings, local i.e. $v$-adic metrics are distorted or scaled relative to each other (for each $v\in\vl$). This is proved here in \Cref{th:distortions},  \Cref{cor:distortions}.

In the adelic context, global arithmetic of a number field enters through the validity of the product formula of \cite{artin45}. In the present case, the global compatibility of {\bf(2,3)} is forced (in \cite{joshi-teich-def}) by requiring that the (global) product formula holds for normalized $v$-adic metrics (for all $v\in\vl$). [This strategy of using the product formula to globalize the adelic theory is also employed in \iutthr.] However, a key feature of the theory is that one cannot uniformly fix the global normalization at all points of the Arithmetic Teichm\"uller Space i.e. the choice of  normalized $v$-adic metrics (for all $v\in\vl$) cannot be uniformly fixed for all objects parameterized by the Arithmetic Teichm\"uller Space $\tfjxl$ [Mochizuki's assertion to this effect is 
\iutf{IUT1, Example 3.5, Remark 3.5.1}].  

In particular, after \cite{joshi-teich-def}, it is possible to quantify the phrase ``there are many distinct \textit{avatars} or versions  of the arithmetic of a fixed number field $L$,'' (this is claimed in \iutf{IUT1, Page 24}) and also to precisely quantify that in any two avatars, relative to each other, $L$ appears stretched or shrunk in a precise topological sense (at  some or all of its primes). Since number theorists have always viewed any given number field as a rigid object on which the drama of Diophantine Geometry unfolds, the existence of such stretchings or shrinkings may seem quite tantalizing but is mathematically quite robust \cite{joshi-teich-def} (and this is an important feature common with Classical Teichm\"uller Theory). 

In particular, after the present paper and \cite{joshi-teich-def}, the arithmetic of a fixed number field $L$ itself (in the precise sense of \cite{joshi-teich-def}) can be treated as a dynamical variable, and this opens the possibility of applying techniques of Arithmetic Dynamics \cite{silverman-dynamics-book} to this dynamical variable.

\subsubsection{The genesis of the theory}
As noted earlier, the possibility that an Arithmetic Teichm\"uller Theory of the type {\Cref{ss:gatt}\bf(1,2,3)} exists is suggested by Shinichi Mochizuki for the special case of elliptic curves in \iutf{IUT1-IUT3}. His proof of the $abc$-conjecture in \iutf{IUT4} remarkably uses the idea of averaging over such Teichm\"uller data (for each fixed number field $L$).  From the standpoint of this paper, Mochizuki's strategy of proof of the $abc$-conjecture may be understood as averaging over the adelic Arithmetic Teichm\"uller Space or a suitably chosen subset of it. [The central and foundational difficulty with \iut\ has been that of unambiguously ascertaining the claimed existence of Teichm\"uller data. A detailed discussion of this can be found in \cite{joshi-final}, \constrfour{\ssep}{se:intro}.] It should be noted that averaging over Classical Teichm\"uller spaces is a  well-established paradigm \cite{wright19}.  

\subsubsection{Mochizuki's Key Principle of Inter-Universality}
It should be noted that Mochizuki neither suggests nor constructs Arithmetic Teichm\"uller Spaces in his work.  Generally speaking, as well as in substantantive details, Mochizuki's approach is  purely group-theoretic and this makes it estremely difficult to establish the existence of Arithmetic Holomorphic Structures. Mochizuki also recognized this problem. Mochizuki proposed to solve this problem via his key \emph{Principle of Inter-Universality} \iutf{\ssep I3, Pages 24-27}. This is the founding principle of his theory and it gives his theory its name. This principle suggests  the use of distinct Grothendieck universes to serve as proxies for the missing Arithmetic Holomorphic Structures i.e. the missing Teichm\"uller data. 

Unfortunately, this principle is stated without proof by Mochizuki (\iutf{\ssep I3, Pages 24-27}).  My treatment of this Key Principle is given in \ssep\ref{appendix-univ}.  Especially, \Cref{thm:zfc}, which provides a precise formulation and proof of Mochizuki's key principle, and notably of Mochizuki's assertion that each geometric base-point for tempered fundamental groups provides a label for a Grothendieck Universe which is distinguished by its label and which, as Mochizuki repeatedly reminds us, serves as a place-holder for the missing Arithmetic Holomorphic Structure.  Thus, the results of this paper provide a way to view Mochizuki's Theory with greater clarity and precision.

However,  it should be noted that my approach to establishing the Theory of Arithmetic Teichm\"uller Spaces is based on entirely different set of ideas from those espoused by Mochizuki and  is grounded in holomorphic function theory (as detailed in \Cref{se:grothendieck-conj}), wherein group theory manifests itself as symmetries of spaces on which the function live. Thus, my approach puts the theory on par with classical Teichm\"uller Theory from the very beginning. 

\subsubsection{The three roads leading to Arithmetic Teichm\"uller Theory}
On the other hand, let me also say that, regardless of how one wants to approach the issue of establishing the existence of Arithmetic Teichm\"uller Spaces using either 
\benumlab
\item \present, or

\item \textit{Mochizuki's Key Principle of Inter-Universality} \iutf{IUT1, \ssep I3, Pages 25--26} (which underlies his theory, but treated correctly as discussed in \ssep\ref{ss:key-principle}),  or
\item the theory of tempered fundamental groupoids (this is essentially {\bf(2)}), or
\item the Theory of Diamonds \cite{scholze-diamonds} (\Cref{re:diamonds}),
\eenum
all point to the existence of a unique (Arithmetic) Teichm\"uller Theory which is Arithmetic-Geometric in origin, and is as established in \present.

\subsection{Main innovations of this paper (Constructions I) and Constructions I(1/2)}
In this subsection, and all sections of this paper except for Section 3, the term `local constructions' will mean constructions over $p$-adic fields, while `global constructions' will mean constructions over number fields. In both the cases, varieties will \emph{always}  assumed to be quasi-projective and smooth over the relevant fields.  Thus, for instance construction of the local Arithmetic Teichm\"uller Spaces will mean construction of Arithmetic Teichm\"uller Spaces over $p$-adic fields etc.

I emphasize that the order of reading  (Constructions I, I(1/2) which deal with Arithmetic Holomorphic Structrues and Arithmetic Teichm\"uller Spaces). These two papers are independent of the remaining papers in this series of papers. [For a very brief discussion of the remaining papers in this series (The Geometric Case, Constructions II, IIIA, IIIB, IV) see \Cref{ss:remaining}.]

Along with (Constructions I, I(1/2)), readers will also find my Anabelomorphy paper \cite{joshi-anabelomorphy} useful as it provides concrete examples of how familiar notions and constructions of arithmetic geometry behave if the absolute Galois group is fixed, but not the base field. My approach to deformations of additive structures is detailed in my Formal Groups paper \cite{joshi-formal-groups}. The Anabelomorphy and Formal Groups papers can be read independently of the present series of papers.

Let me give a brief overview of the key innovations which are carried out in each of the present series of papers:
\subsubsection{Constructions I (this paper): Arithmetic Holomorphic Structures on Algebraic Varieties over Number Fields} The main innovations in this paper are as follows:
\begin{itemize}
	\item A robust local definition of Arithmetic Holomorphic structures \Cref{def:arith-hol-strs} via \Cref{thm:main2}, preserving anabelian properties.
	\item The establishment of the local and global arithmetic Teichm\"uller spaces (local \Cref{def:arith-hol-space-local}, \Cref{def:arith-hol-space-local-arch}), (global \Cref{th:main4}).
	\item Proofs of the key assertions,  tacitly used in \iut, that symmetries of certain monoids i.e.  operations performed on certain monoids, give rise to changes of Arithmetic Holomorphic Structures. A typical result obtained here is \Cref{th:main3}--a result reminiscent of the classical Virasoro uniformization of \cite{beilinson88}.
	\item \Cref{pr:local-diamond} establishes the relationship between the local Arithmetic Teichm\"uller spaces constructed here and diamond $X^\Diamond_{E}$ associated to any quasi-projective variety $X$ over $p$-adic field $E$ in \cite[Section 15]{scholze-diamonds}. 
	\item In Classical Teichm\"uller theory one encounters variation of periods of Riemann surfaces while the biholomorphism class of the Riemann surface remains fixed. A typical example of this type is this: the complex upper half-plane $\fH\subset \C$ is the classical Teichm\"uller space in genus one, and two distinct points of the classical Teichm\"uller space $\tau\in\fH$ and $-\frac{1}{\tau}\in\fH$ represent the biholomorphism class of the same Riemann surface (namely the elliptic curve with period lattice $[1,\tau]$), and at the same time if $|\tau|$ is small, then $|-1/\tau|$ is large and vice versa. \iut\ rests on Mochizuki's assertion  that this archimedean phenomenon has  a $p$-adic analog, and loc. cit.  rests on exploiting such a variation of $p$-adic periods of elliptic curves \cite[Page 3]{mochizuki-gaussian}. The existence of this $p$-adic variation is established here in \Cref{th:distortions}.
	\item One may think of (non-archimedean) \Cref{thm:main2} as providing a variation/deformation of the Banach ring structure of analytic functions. In the archimedean case, and especially, for Riemann surfaces, classical work of \cite{nakai59,nakai60,nakai1972} shows that Classical Teichm\"uller Theory admits a similar interpretation as a variation/deformation of Banach ring structure of Royden algebras of functions. The archimedean case is detailed in \Cref{app:archimedean-theory2}, \Cref{th:nakai1} (for Riemann surfaces) and in general in \Cref{def:arith-hol-space-local-arch}, \Cref{th:nakai2}. This not only puts the non-archimedean Teichm\"uller Theory of the present paper and (archimedean) Classical Teichm\"uller  Theory on a parallel footing, but it also allows us to include quasi-conformal mappings in the archimedean theory.
	\item In \Cref{appendix-univ}, I provide a proof of \Cref{thm:zfc} which is asserted without proof in \iut, but which underlies \emph{Mochizuki's Key Principle of Inter-Universality}. \end{itemize}  
It should be noted that the approach developed here may be termed Construction of Arithmetic Holomorphic Structures as opposed to Anabelian Reconstruction approach of \iut.

\subsubsection{Constructions I(1/2) \cite{joshi-teich-def}: Arithmetic Teichm\"uller Spaces of Number Fields} As stated in the Introduction, a remarkable and beautiful idea due to Mochizuki in \iut, is that just as a fixed Riemann surface has many distinct Teichm\"uller avatars (quantified by points of the relevant Teichm\"uller space), a Number Field $L$ has Teichm\"uller avatars (called \emph{arithmeticoids of $L$} in \cite{joshi-teich-def}). However, the arithmetic case is particularly difficult to quantify using the methods, and the framework, presented in \iut. 
\begin{itemize}
	\item In \cite[\ssep 5.1]{joshi-teich-def}, using the results of the present paper, I provide a canonical construction of distinct Teichm\"uller avatars (called \emph{arithmeticoids} in \cite[{Definition~\ref{II5-def:arithmeticoid}}]{joshi-teich-def}) of a fixed number field as asserted in \iutf{IUT1, \ssep I3, Page 25}.
	\item Going beyond Mochizuki's ideas, \cite[{Theorems~\ref{II5-th:distance-bet-arith}, \ref{II5-th:parameter-space-of-deformations}}]{joshi-teich-def} constructs a metrisable,   arithmetic Teichm\"uller space $\yadl$,  of such avatars. In this space, properties asserted by Mochizuki (such as the $\log$-Link mentioned above) appear as symmetry properties.
	\item For a smooth, quasi-projective variety $X/L$, \constrtwoh{Theorem }{th:consequence-inequivalent-arithmeticoids} shows how the choice of arithmeticoids changes the  Berkovich analytic structure  for all of any non-empty set of primes of $L$. An immediate consequence is that the theory underlying \iut\ is highly non-vacuous--a point which has been difficult to demonstrate by purely using methods of \iut--leading to some doubts about its claims.
	\item A key innovation of \cite{joshi-teich-def} is the construction of global Frobenius morphism of a Number Field. The reader is referred to \constrtwoh{\ssep }{ss:intro-frob} for a detailed discussion of this notion and its properties (which are established in \constrtwoh{\ssep}{se:global-frobenius-num-fld}). 

\item A second key innovation of \cite{joshi-teich-def} is the construction of the action $L^*\act\yadl$ which may be terms as a global Arithmetic Virasoro action. This is detailed in \constrtwoh{Theorem }{th:galois-action-on-adelic-ff}.
\end{itemize}  
This completes the circle of ideas originating in the Dedekind-Weber analogy between Riemann surfaces and Number Fields.

\subsubsection{The remaining papers in this series}\label{ss:remaining}
The other papers in the series (The Geometric Case, Constructions II, IIIA, IIIB, IV) deal with application to the $abc$-conjecture following Mochizuki's ideas.  
\begin{tcolorbox}[label=tab:analogies3,grow to right by = 1cm, colback=white,coltitle=black,colframe=white!25!olive,fonttitle=\bfseries,title={Brief guide to the contents of remaining papers in this series}]
\begin{description}
	\item[\hyperlink{cite.joshi-teich-estimates}{Constructions II}] provides a canonical local construction of Mochizuki's $\Theta$-Link as a correspondence (on the Arithmetic Teichmuller Space) similar to a Hecke correspondence.
	\item[\hyperlink{cite.joshi-teich-IIIB-rosetta}{Constructions IIIA}] provides  
	\begin{itemize}
		\item a `Rosetta Stone' \constrthrB{\ssep }{se:intro-rosetta-stone} for a parallel reading of this work and \iut,
		\item a refinement of Mochizuki's Anabelian Reconstruction Theory \cite{mochizuki-topics3}  which allows one to reconstruct Arithmetic Holomorphic Structures (in the sense of this paper) from the tempered fundamental groups.  This refinement appears in \constrthrB{\ssep }{ss:supplement-anab-recon} and provides a pathway to apply Mochizuki's Anabelian Reconstruction ideas for disambiguating assertions of \iutf{IUT3, Theorem 3.11 and Corollary 3.12}.
		\item  A global construction of Mochizuki's $\Theta$-values locus (referred to as \emph{Mochizuki's Ansatz} in \cite{joshi-teich-estimates,joshi-teich-IIIB-rosetta,joshi-teich-IIIC-moccor}) and which is the key object of interest in \iutf{IUT3, Theorem 3.11 and Corollary 3.12}.
	\end{itemize} 
	\item[\hyperlink{cite.joshi-teich-IIIC-moccor}{Constructions IIIB} and \hyperlink{cite.joshi-teich-abc-conj}{Constructions IV}]  deal with \moccor\ and the $abc$-conjecture following Mochizuki's strategy.   More details can be found in the introductions of these papers.
\end{description}
\tcblower
\emph{Constructions I, I(1/2), II, IIIA are adequate for arriving at a completely unambiguous and a canonical reading of \iutf{IUT1-IUT3}.}
\end{tcolorbox}
\subsection{Acknowledgments} 
The theory presented here is  inspired by Mochizuki's  work on anabelian geometry \cite{mochizuki96-gconj,mochizuki99}, $p$-adic Teichm\"uller Theory \cite{mochizuki-teichmuller-book} and Inter-Universal Teichm\"uller Theory \iut.  So,  my intellectual debt to Shinichi Mochizuki cannot be overstated. \emph{Notably this paper could not have existed without his work.}  

I thank Peter Scholze for some correspondence regarding \Cref{thm:main,thm:main2}, the definition Arithmetic Holomorphic Structures,  and for discussions about their role in \iut\ in May-June 2024. I thank Kiran Kedlaya for correspondence, comments and encouragement. One well-known expert (on rigid analytic geometry, Berkovich spaces, arithmetic geometry) with whom I corresponded while writing of this paper (especially \Cref{se:grothendieck-conj}), has chosen to remain anonymous, but  our correspondence is acknowledged here with pleasure.

Thanks are due to Yuichiro Hoshi  for answering many of my elementary questions about tempered fundamental groups and \iut; also to Taylor Dupuy, Anton Hilado for some comments and conversations regarding the 2020 versions  of this paper.  I thank Jacob Stix for providing a correction to \ssep\ref{pa:temp-fun-grp-vs-etale-fun-grp} and for answering a question regarding \cite{stix16}. I also thank Emmanuel Lepage for some correspondence in the context of \cite{lepage-thesis} and \Cref{thm:main,thm:main2}.

\newcommand{\ebh}{\widehat{\bE}}
\newcommand{\ebhx}[1][x]{\widehat{\bE^{#1}}}
\newcommand{\bdr}{B_{dR}}
\newcommand{\bdre}{{\bdr}_{,E}}
\newcommand{\bdrep}{B^+_{dR,E}}
\newcommand{\kbh}{\widehat{\bK}}
\newcommand{\sR}{\mathscr{R}}
\newcommand{\rs}{\sR_\Sigma}
\newcommand{\rsp}{\sR_{\Sigma'}}

\section{Classical Teichm\"uller Theory}\label{se:classical-teich-thry}
\addtocontents{toc}{\protect\setcounter{tocdepth}{1}}
\nwss
\subsection{Classical Teichm\"uller Theory in a Nutshell}

Since the principal result of this paper is the demonstration of the existence of a $p$-adic Arithmetic Teichm\"uller Theory, it will be useful to recall the salient theorems  of Classical Teichm\"uller Theory (these are organized in \Cref{th:main-classical} stated below). This list of results, and others, be gleaned from standard texts on Classical Teichm\"uller Theory (\cite{lehto-book}, \cite{nag-book}, \cite{imayoshi-book}). This theorem is the one which my work imitates in the arithmetic case. 

\bthm\label{th:main-classical}
Let $\Sigma$ be a connected, compact Riemann surface of topological type $g\geq 0$ and let $\sT_{g}=\sT(\Sigma)$ be the Teichm\"uller space of $\Sigma$ and let $\sM_{g}$ be the moduli of Riemann surfaces. Then 
\benumlab
\item $\sT_g$ parameterizes, up to a suitable equivalence hereafter dubbed \text{Teichm\"uller equivalence}, pairs $(\Sigma',f)$ consisting of a compact Riemann surface $\Sigma'$ and a non-constant quasi-conformal mapping $f:\Sigma'\to \Sigma$  \cite[Chapter V, 2.1]{lehto-book}. Note that, in general, $f$ need not be a holomorphic mapping (and hence may be a non-algebraic mapping), but $f$ is a homeomorphism, so $\Sigma'$ has the same topological type as $\Sigma$. Importantly, $\sT_g$ is not an algebraic variety or a scheme.]
\item The space $\sT_g=\sT(\Sigma)$ can also be  described as the set of equivalence classes (the equivalence relation is not recalled here) of pairs  $(\Sigma,ds^2)$ consisting of the Riemannian manifold $\Sigma$  equipped with an arbitrary Riemannian metric $ds^2$ (\cite[Theorem 1.8]{imayoshi-book}).
\item For $(\Sigma',f)\in \sT_g$, let $[\Sigma']$ denote the bi-holomorphism class i.e. the isomorphism class of the Riemann surface $\Sigma'$. Then one has a natural, surjective mapping $\sT_g\to \sM_{g}$ given by $(\Sigma',f)\mapsto [\Sigma']$.
\item There exists an infinite discrete group, denoted $\Gamma_g$ and called the Teichm\"uller modular group or the mapping class group, which acts freely on $\sT_g$ and $\sT_g\to M_g$ is the quotient mapping
$\sT_g\to \sT_g/\Gamma_g=\sM_g$ (\cite[Chapter 6]{imayoshi-book}).
\item For $g=1$ (elliptic curves), $\sT_1=\fH$ is the upper half-plane, $\Gamma_1=\PSL_2(\Z)$, and $\sM_1=\fH/\PSL_2(\Z)=\C$ \cite[Chapter 1]{imayoshi-book}.
\item The fiber of the natural morphism $\sT_{g}\to \sM_{g}$  over the isomorphism class of $[\Sigma]\in\sM_{g}$ provides a collection of pairs $(\Sigma',f)\in \sT_g$ for which $[\Sigma']=[\Sigma]$.
\item There exists a Virasoro symmetry i.e.  the Virasoro Uniformization Theorem (of \cite{kontsevich87}, \cite{beilinson88}) holds (see \cite[Theorem 17.3.2]{frenkel01-book} for an accessible exposition). [For more on this Virasoro symmetry see \cref{ss:virasoro-symmetry}.]
\item For every $(\Sigma',f)\in \sT_g$, one has an isomorphism $\pi_1(\Sigma')\mapright{\isom} \pi(\Sigma)$  between the fundamental groups of $\Sigma'$ and $\Sigma$ given by the homeomorphism underlying  the quasi-conformal mapping $f:\Sigma'\to \Sigma$.
\eenum
\ethm

\brem It should be noted that \Cref{th:main-classical}{\bf(7)} is generally not a part of Classical Teichm\"uller Theory textbooks \cite{lehto-book,nag-book,imayoshi-book}. On the other hand,
one of the important aspect of  this paper is the demonstration of the existence of the analog of the Virasoro action \Cref{th:main-classical}{\bf(7)} in the theory of Arithmetic Teichm\"uller Spaces. [While, the underlying property (which I prove) appeared under the name of \textit{`Indeterminacy Ind2'} in Mochizuki's work \iutf{IUT3, Theorem 3.11 and Corollary 3.12}, I recognized it as being the analog of the Virasoro Uniformization Theorem, and gave a precise formulation and proof in \Cref{th:main3}.] 
\erem

\subsection{Classical Teichm\"uller Theory as a functor in Banach Algebras}
\label{app:archimedean-theory2}
Classical Teichm\"uller Spaces admits description as a functor in infinite dimensional Banach algebras (see Theorem~\ref{th:nakai1} below). This rests on a remarkable Banach ring theoretic characterization of quasi-conformal and conformal mappings of Riemann surfaces due to \cite{nakai59,nakai60,nakai-book}. 

In particular, this implies that the $p$-adic theory of this paper based on the idea that the $p$-adic Banach algebra structure can be made to vary (\Cref{thm:main2}, \Cref{th:variation-of-arith-hol-strs}), has an archimedean counterpart in classical Teichm\"uller Theory. This puts the $p$-adic Arithmetic Teichm\"uller Theory of this paper  on par with classical Teichm\"uller Theory. [I discovered Mitsuru Nakai's aforementioned works after I had established \cite{joshi-untilts-2020}.]

\para Let $\Sigma$ be a connected Riemann surface (open or closed). Recall (from \cite{royden1953}, \cite{nakai60}) that the \emph{Royden algebra $\rs$} of a connected Riemann surface $\Sigma$ is the $\C$-algebra of all complex valued functions $f:\Sigma\to \C$ satisfying  the following three properties:
\benumlab
\item $f$ is absolutely continuous on $\Sigma$ in the Tonelli sense (see \cite[Chap. III, 1A, Page 147]{nakai-book} for the definition),
\item $|f|$ is bounded on $\Sigma$, and
\item the Dirichlet integral $$D[f]=\int\int_\Sigma \left(\left|\frac{\partial^2f}{\partial x^2}\right|+\left|\frac{\partial^2f}{\partial y^2}\right|\right)dx\cdot dy <\infty.$$
\eenum
For $f\in\rs$, define $$\abs{f}_\Sigma=\displaystyle{\sup_\Sigma \abs{f}}+ \sqrt{D[f]}.$$
The following is proved in \cite{royden1953} \cite{nakai59}. 
\bthm
Let $\Sigma$ be any connected Riemann surface. 
\benumlab
\item The mapping $\abs{-}_\Sigma:\rs \to \R$ given by $$f\in\rs \mapsto \abs{f}_\Sigma$$ defines  a norm on $\rs$, and
\item  the normed algebra $(\rs,\abs{-}_\Sigma)$  is a commutative, complex Banach algebra.
\eenum
\ethm

\brem\  
\benumlab
\item Note that $\rs$ is an infinite dimensional, commutative Banach algebra even if $\Sigma$ is a compact Riemann surface.
\item Note that $\rs$ is also equipped with several different topologies other than its norm-topology. The norm topology is the strongest among these topologies (\cite{nakai60}). 
\eenum 
\erem

\para 
\newcommand{\sTs}{\mathscr{T}_\Sigma}
Let $T_\Sigma$ be the Teichm\"uller space of $\Sigma$ and let $\sTs$ be the  category whose objects are the points of $T_\Sigma$ i.e. it objects are Riemann surfaces $(\Sigma',f:\Sigma' \to \Sigma)\in T_\Sigma$. A morphism between $(\Sigma',f:\Sigma' \to \Sigma)\to (\Sigma'',g:\Sigma'' \to \Sigma)$ is a quasi-conformal mapping $h:\Sigma'\to \Sigma''$. 

\para The next theorem is my reformulation of  main theorem of \cite{nakai59,nakai60} it provides a remarkable algebraic characterization of quasi-conformal mappings (the theory of \cite{clausen-scholze-lectures} espouses an algebraization of analysis in greater generality):
\bthm\label{th:nakai1}
Let  $\Sigma,\Sigma'$ be two connected two Riemann surfaces. Then 
\benumlab
\item $\Sigma,\Sigma'$ are quasi-conformal if and only if one has a $\C$-algebra isomorphism $$\rs\isom \rsp$$  between their Royden algebras. 
\item Any $\C$-algebra isomorphism $\rs\isom \rsp$ is a homeomorphism for the respective norm-topologies.
\item Moreover, $\Sigma$ and $\Sigma'$ are conformally equivalent (i.e. the Riemann surfaces are isomorphic) if and only if $$(\rs,\abs{-}_\Sigma)\isom(\rsp,\abs{-}_{\Sigma'})$$ is an isometry of Banach algebras. 
\item The construction of $(\rs,\abs{-}_\Sigma)$ is functorial for quasi-conformal mappings between Riemann surfaces.
\item In particular one has a functor $$(\Sigma',f:\Sigma'\to \Sigma)\mapsto (\rsp,\abs{-}_{\Sigma'})$$ from the Teichm\"uller category $\sTs$  to the category of infinite dimensional complex Banach algebras.
\eenum
\ethm
\bp 
The first assertion is \cite[Theorem III.1]{nakai60}, the second assertion is \cite[Theorem III.4]{nakai60}, the third assertion is \cite[Theorem III.4(i)]{nakai60}. The fourth assertion is \cite[Theorem III.3(i)]{nakai60}. The last assertion should be clear from the previous four assertions.
\ep

\brem\  
\benumlab
\item According to \constranab{Lemma }{le:anab-arch-prop}, a quasi-conformal mapping between Riemann surfaces is an anabelomorphism of Riemann surfaces.
\item Thus,  quasi-conformality of Riemann surfaces has a purely algebraic characterization and conformality of Riemann surfaces has a purely Banach algebra theoretic characterization.
\item Note that by \Cref{th:nakai1}{\bf(5)}, $\sT_{\Sigma}$ may be considered as an archimedean version of a Banach-Colmez Space of \cite{colmez2000}. 
\eenum
\erem

\brem  Now fix a connected Riemann surface $\Sigma$ and let $\sR=\rs$ be the Royden algebra of $\Sigma$. Let $f:\Sigma'\to \Sigma$ be a quasi-conformal mapping. Then using the isomorphism $\rsp\to \rs$ induced by $f$ one can pull-back the norm $\abs{-}_{\Sigma'}$ of $\rsp$ and  write $(\rs,\abs{-}_{\Sigma'})$ for the norm on $\rs$ given using $(\Sigma',f)$. Thus for each $(\Sigma',f)\in T_\Sigma$, one has a Banach algebra norm on the fixed topological ring $\sR=\sR_\Sigma$. In particular classical Teichm\"uller theory provides a distinguished collection 
$$\left\{ (\sR, \abs{-}_{\Sigma'}): \Sigma'\in T_{\Sigma} \right\}$$
of Banach algebra norms  on the fixed topological ring $\sR$. Thus one may  view classical Teichm\"uller Theory as deformations of Banach algebra structure  on the fixed $\C$-algebra $\sR$ i.e. classical Teichm\"uller Theory is a variation of Banach ring structures. Notably, this observation brings the classical theory closer to the $p$-adic theory of Banach-Colmez Spaces \cite{colmez2000}.

The constructions of the present paper can be similarly viewed: as changes in the Banach ring structure of functions as one deforms the valued field $\C_p$.
\erem

\section{The Absolute Grothendieck Conjecture is false for Berkovich Spaces}\label{se:grothendieck-conj}
\nwss
By \Cref{th:main-classical}{\bf(8)}, one knows that each $(\Sigma',f)\in\sT_{\Sigma}$ has topological fundamental group isomorphic to that of $\Sigma$, and hence the fundamental group does not determine the datum $(\Sigma',f)\in\sT_{\Sigma}$ nor does it determine the complex structure of $\Sigma'$. One can say, therefore, that the Absolute Grothendeick Conjecture is trivially false for Riemann surfaces because of the existence of distinct Teichm\"uller data $(\Sigma',f)$ and this non-validity is a manifestation of the existence of distinct complex analytic structures. 

The $p$-adic approach (to Arithmetic Teichm\"uller Spaces) of this paper began with a demonstration that the Absolute Grothendieck Conjecture is false in the category of Berkovich analytic spaces while the tempered fundamental groups remain isomorphic. This is established in  \Cref{thm:main2}. [A brief version of this result appeared in \cite{joshi-untilts-2020}.] This leads to a precise definition of Arithmetic Holomorphic Structures in \Cref{se:arith-hol-strs} and to the construction of Arithmetic Teichm\"uller Spaces in \Cref{se:construct-att}.

\para\label{se:tempered-groups} \emph{All valuations on base fields considered in this paper will be rank one valuations.} For the theory of tempered fundamental groups see \cite{andre-book,andre03} or \cite{lepage-thesis}. As is noted in \cite{andre-book}, tempered fundamental groups are natural in the $p$-adic analytic contexts because they capture  finite \'etale coverings \emph{and} discrete coverings such as those arising from Tate or Mumford Uniformization available in the $p$-adic contexts. Berkovich spaces  (see \cite{berkovich-book} and \cite{berkovich93}) will be strictly analytic (and mostly will arise as analytifications of geometrically connected smooth quasi-projective varieties).

\para In what follows I will work with algebraically closed, perfectoid fields of characteristic zero. A typical example of such a field is the completed algebraic closure $\C_p$ of $\Q_p$. Such fields can also be characterized in many different ways.  For the convenience of the readers unfamiliar with perfectoid fields, the following lemma (immediate from \cite[Definition 3.1]{scholze12-perfectoid-ihes}), provides a translation of this condition into  more familiar hypothesis.

\blem\label{lem:perfectoid}
Let $K$ be a valued field and let $R\subset K$ be the valuation ring and assume that $\abs{p}_K<1$. The following conditions are equivalent:
\benumlab
\item $K$ is an algebraically closed field, complete with respect to a rank one non-archimedean valuation with residue characteristic $p>0$.
\item $K$ is an algebraically closed, perfectoid field.
\eenum
\elem

\bp A perfectoid field has residue characteristic $p>0$ and is complete with respect to a rank one valuation. So  (2)$\implies$(1) is trivial. So it is enough to prove that (1)$\implies$(2).   I claim that  Frobenius $\phi:R/pR\to R/pR$ is surjective. let $\bar{x}\in R/pR$ and suppose $x\in R$ is an arbitrary lift of $\bar{x}$. Then as $K$ is algebraically closed, $x^{1/p}\mod{pR}$ provides a lift of $\bar{x}$. As $K$ is complete with respect to a rank one valuation and Frobenius is surjective on $R/pR$, so  $K$ is perfectoid by \cite[Definition 3.1]{scholze12-perfectoid-ihes} and by my hypothesis $K$ also algebraically closed. This proves (1)$\implies$(2).
\ep

\para For a perfectoid algebraically closed field $K$ as above, one has naturally associated field $K^\flat$, algebraically closed, perfectoid of characteristic $p>0$, called the \emph{tilt of $K$} and $K$ is called an untilt of $K^\flat$ (see \cite[Lemma 3.4]{scholze12-perfectoid-ihes}).

\para Fix an  algebraically closed  field, perfectoid $F$ of characteristic $p>0$ (see \cite{scholze12-perfectoid-ihes}). For example readers can simply assume, without any loss of generality, that $F= \C_p^\flat$ as this case is quite adequate for my purposes.

\para By an $E$-\emph{untilt} of $F$, I will mean the data $(E\into K, K^\flat\isom F)$ where $E\into K$ is an isometric embedding of a $p$-adic field $E$ into a perfectoid field $K$, of characteristic zero, with an isometry $K^\flat\isom F$. If $E=\Q_p$, I will simply say an untilt of $F$ instead of $\Q_p$-untilt of $F$. Note that by \cite[Proposition 3.8]{scholze12-perfectoid-ihes}  $K$ is algebraically closed as its tilt $K^\flat=F$ is algebraically closed (by my hypothesis). If $F=\cpt$  then $K^\flat$ is isometric with $\cpt$. By the theory of \cite{fargues-fontaine} untilts $K$ of $F$ exist and are parametrized by Fargues-Fontaine curves.

\para Let $E$ be a $p$-adic field which is fixed for the present discussion. I will work with untilts $K$, of $F$, equipped with continuous embeddings $E\into K$ with the valuation of $K$ providing a valuation on $E$ which is equivalent to the natural $p$-adic valuation on $E$. By \cite{fargues-fontaine} for a given pair $(F,E)$,  such fields $K\hookleftarrow E$, exist and are parametrized by Fargues-Fontaine curves  (denoted here by $\sxfe$). Without further mention, all untilts $K$ will be assumed to be of this type (for our chosen $p$-adic field $E$).

\para Crucial point for this paper is that \emph{there exist untilts of $\cpt$ which are not topologically isomorphic.} This is the main result of \cite[Theorem~1.3]{kedlaya18} (also see \cite{matignon84}). Note that all characteristic zero untilts of $\cpt$ have the cardinality of $\C_p$ and are complete and algebraically closed fields and hence are abstractly isomorphic fields but may not be topologically isomorphic after \cite[Theorem 1.3]{kedlaya18} (also see \cite{matignon84}).

\para Now fix a geometrically connected, smooth quasi-projective variety $X/E$, with $E$ a $p$-adic field. Let $\xan/E$ be the strictly analytic space associated to $X/E$. Let \be\pit{X/E}=\pi_1^{temp}(\xan/E)\ee be the tempered fundamental group of the strictly $E$-analytic space associated to $X/E$ in the sense of \cite{andre03} or \cite{andre-book}. 

(\emph{Note that my notation $\pit{X/E}$ suppresses the passage to the   analytification $\xan/E$  for simplicity of notation. The theory of (tempered) fundamental groups also requires a choice of base point which will be suppressed from my notation for the moment.})

\para\label{pa:temp-fun-grp-vs-etale-fun-grp} The relationship between the tempered fundamental groups and the  \'etale fundamental groups of geometrically connected, smooth quasi-projective varieties as follows. 

Let $X/E$ be a geometrically connected, smooth, quasi-projective variety over a $p$-adic field $E$. Then one has a natural homomorphism (\cite[Proposition 4.4.1]{andre03}, \cite[Section 2.1.4]{andre-book}):
$$\pit{X/E}\to\pi_1^{et}(X/E),$$ 
which is injective if $\dim(X)=1$ (\cite[Chapter III, Prop. 2.1.6]{andre-book}). In any dimension,
the image of the above homomorphism is dense and moreover  $\pi_1^{et}(X/E)$ is the profinite completion (\cite[Proposition 4.4.1]{andre03}):
$$\widehat{\pit{X/E}}=\pi_1^{et}(X/E).$$ 

\para Let $E'/E$ be a finite extension of $E$ with a continuous embedding $E'\into K$ (as $K$ is algebraically closed, valued field containing $E$, such $E'$ exists).   One can consider $X_{E'}=X\times_{E}E'$ (similarly $X_K=X\times_EK$). Then one has an exact sequence by \cite[Prop. 2.1.8]{andre-book}
$$1\to \pit{X_{E'}/E'}\to \pit{X/E}\to \gal(E'/E)\to 1.$$

Let $\bE\subseteq K$ be the algebraic closure of $E$ contained in $K$.

By  varying $E'$ over all finite extensions of $E\into K$  one obtains (see \cite[Section 5.1]{andre03}) an exact sequence of topological groups:
$$1\to  \ilim_{E'/E}\pit{X_{E'}/E'}\to\pit{X/E} \to \gal(\bE/E)\to 1.$$

\newcommand{\Et}{\widetilde{E}}

\bthm\label{thm:main}
Let $F$ be an algebraically closed perfectoid field of characteristic $p>0$ (for example $F=\cpt$). Let $E$ be a $p$-adic field. Let $K,K_1,K_2$ be arbitrary untilts of $F$ with continuous embedding $E\into K$ (resp. into $K_1$ and $K_2$). Let $\bE$ (resp. $\bE_1,\bE_2$) be the algebraic closure of $E$ in $K$ (resp. in $K_1,K_2$). Let $X/E$ be a geometrically connected, smooth, quasi-projective variety over $E$. Then one has the following:
\benumlab
\item\label{thm:main-1} a continuous isomorphism $$\pit{X/K}\isom \ilim_{E'/E}\pit{X/E'},$$
where the inverse limit is over all finite extensions $E'$ of $E$ contained in $K$, and a
\item\label{thm:main-2} a short exact sequence of topological groups $$1\to \pit{X/K}\to \pit{X/E}\to G_E\to 1,$$ and
\item\label{thm:main-4} In particular for any two untilts $K_1,K_2$ of $F$, one has a continuous isomorphism $$\pit{X/K_1}\isom \pit{X/K_2}.$$
\eenum
\ethm
\bp 
The assertion \ref{thm:main-1} is true assuming only that $K$ is a complete algebraically closed field containing $E$ isometrically and is due to \cite{lepage-thesis}. My own proof of \ref{thm:main-1}, before I found the assertion in \cite{lepage-thesis},  was by reworking of \cite[Prop. 5.1.1]{andre03} for any $K$ algebraically closed perfectoid field, and I was interested in proving \ref{thm:main-1} because  I wanted to establish \ref{thm:main-4} (whose importance  will become clear in Theorem~\ref{thm:main2} below).   Here I provide an approach to the proof of \ref{thm:main-1} via the reduction to the principle of invariance of fundamental groups under extension of algebraically closed fields (also due to \cite{lepage-thesis}), for completeness. \emph{So  \ref{thm:main-4} is the new and important  observation here--from the point of view of Theorem~\ref{thm:main2} below.}

Let me remind the reader that my hypothesis on  $K,K_1,K_2$ imply that $K,K_1,K_2$ are algebraically closed and complete with respect to a rank one valuation.  

Let me prove \ref{thm:main-1}, this will also lead to \ref{thm:main-2}. Since $K$ is  algebraically closed, it follows that $K$ contains an algebraic closure $\bE$ of $E$.  Let $\Et\subseteq K$ be the closure  (with respect to valuation topology of $K$) of $\bE$. 

It is clear that $\Et\supset \bE$ is complete and algebraically closed field and $\Et$ contains the algebraic closure $\bE\subset K$ of $E$ contained in $K$ as a dense subfield. In particular $\Et$ is the completion of $\bE$ with respect to the induced valuation. In other words $\Et$ is a copy of the completion of an algebraic closure of $E$ (usually denoted $\ebh$) equipped with an isometric embedding $\iota:\ebh\into K$ with $\iota(\ebh)=\Et$.
Hence  $K/\Et$ is an isometric extension of algebraically closed, complete valued fields (with rank one valuations). 

Now one can apply  the principle of invariance  of  fundamental groups under passage to extensions of algebraically closed fields. This principle is well-known for \'etale fundamental groups of proper varieties (see \cite[Expos\'e X, Corollaire 1.8]{sga1}). For  tempered fundamental groups (and $X$ not necessarily proper) this principle is proved in \cite[Proposition 2.3.2]{lepage-thesis}. Thus  applying \cite[Proposition 2.3.2]{lepage-thesis} to the extension $K/\Et$ one has an isomorphism of topological groups $$\pit{X_K}\isom \pit{X_{\Et}}.$$

On the other hand by \cite[Proposition 5.1.1]{andre03}, as $\Et$ is the completion of the algebraic closure of $\bE\subset K$ of $E$, one has an isomorphism
\be \pit{X/\Et}\isom \ilim_{E'/E}\pit{X_{E'}/E'}
\ee
and  an exact sequence of topological groups 
$$1\to \pit{X/K} \isom \ilim_{ E'/E}\pit{X_{E'}/E'}\to\pit{X/E} \to \gal(\bE/E)\to 1.$$
 This proves the assertions \ref{thm:main-1}, \ref{thm:main-2} as claimed.
 
Let me now prove \ref{thm:main-4}. The claimed isomorphism $\pit{X/K_1}\isom \pit{X/K_2}$ follows from the fact that both the groups can be identified with $\ilim_{E'/E}\pit{X/E'}$ where the inverse limit is over all finite extensions of $E'/E$ contained in $K_1$ (resp. $K_2$) and the fact that there is an equivalence between categories of finite extensions of $E$ contained in $K_1$ and the category of finite extensions of $E$ contained in $K_2$, since finite extensions of $E$ are given by adjoining roots of polynomials with coefficients in $E$ and this data is independent of the embedding of $E$ in $K_1$ or $K_2$ and moreover any abstract isomorphism of  finite extensions of a complete discretely valued field is in fact an isometry--i.e given a finite extension of $E$,   $E'\into K_1$ contained in $K_1$, there is an isometry $E'\into K_2$ and vice versa.
\ep

\para\label{pa:vertical-horizontal-variations1} The importance of working with algebraically closed perfectoid fields $K_1,K_2$ with isometric tilts $K_1^\flat \isom K_2^{\flat}$ (i.e. with untilts of a fixed algebraically closed perfectoid field of characteristic $p>0$) will become clear from Theorem~\ref{th:main3} which will be proved later.  Note that if $K_1,K_2$ are arbitrary algebraically closed perfectoid fields, then $K_1^\flat$ and $K_2^\flat$ need not be isometric. A simple example of this is given as follows. Let $\cpmax\supset\C_p$ be a maximally (i.e. spherically) complete field (i.e. a maximal immediate extension of $\bQ_p$) (see \cite{kaplansky42}, \cite{poonen93} for the construction of $\cpmax$), then $\cpmax$ and $\C_p$ do not have isometric tilts other wise $\cpt\isom {\cpmax}^\flat$ is also spherically (i.e. maximally) complete, which is certainly not the case. [The field extension $\cpmax\supset\C_p$ is of uncountable transcendence degree (both the fields have the same cardinality).]

\para\label{pa:vertical-horizontal-variations2} The passage from $\C_p$ to $\cpmax$ should be considered as a \emph{``vertical variation''} of the algebraically closed perfectoid field because it involves extension of their tilts $\C_p^{max\flat}\supset\cpt$ (because this is also of uncountable transcendence degree). On the other hand \cite{kedlaya18} shows that a there is also a \emph{``horizontal''} or an \emph{``iso-tilted or equi-tilted variation''} possible in which the tilts stay fixed isometrically.

\para{\label{pa:E-isomorphisms-analytic-spaces}} Let $E$ be a $p$-adic field.  In \cite[Section 2.3, Section 3.1]{berkovich-book},  Berkovich constructs the category of \emph{analytic spaces over $E$} (or more simply the category of \emph{Berkovich spaces} over $E$) (a similar theory is also sketched  in \cite{berkovich93}).  While I refer the reader to these references for the general case, let me recall what this means in the context I will use.  
By \cite[Section 3.1]{berkovich-book} an \emph{analytic space over $E$} is a $K$-analytic space for some  valued field  $K\supseteq E$ (with a rank one valuation) and equipped with an isometric embedding $E\into K$. Let $X/E$ be a quasi-projective variety and let $\xan$ denote analytification  of $X/E$ (in the sense of \cite{berkovich-book}).   Thus the  $K$-analytic space  $\xan_K=\xan\times_{E}K$ is an analytic space  over $E$.  Let $K_1,K_2$ be two valued fields (with rank one valuation) containing $E$ (isometrically).  Let $\xan_{K_1}=\xan\times_E K_1$ and similarly define $\xan_{K_2}$. So one has two analytic spaces over $E$.  By \cite[Section 2.3, Section 3.1]{berkovich-book}  one can consider the notion of  (iso)morphisms  $\xan_{K_1}\mapright{\isom} \xan_{K_2}$ of analytic spaces over $E$.  Specifically, this reduces to defining the notion of (iso)morphisms between affinoid spaces  over $E$. This is done as follows (in the notation of \cite[Section 1.2]{berkovich-book}): if $(\mathscr{M}(A_i),A_i)$ are $K_i$-affinoid spaces over $E$, for $i=1,2$, then an (iso)morphism between them is given by a bounded (iso)morphism, $f:A_1\mapright{\isom}A_2$, of  Banach rings compatible with their structure as normed algebras over the valued field $E$ (and the corresponding continuous (iso)morphism between the semi-norm spectra $\mathscr{M}(A_i)$).  In particular if $E\supseteq \Q_p$ then one can consider  (iso)morphisms of analytic spaces over $\Q_p$.   Thus an isomorphism $\xan_{K_1}\mapright{\isom} \xan_{K_2}$ of analytic spaces over $\Q_p$ makes sense and is the $p$-adic analytic analog, of the notion of  isomorphisms of schemes over $\Z$. [Note that this can be obviously formulated more generally, without assuming that $X/E$ is quasi-projective, but I have restricted myself to the case I will use in Theorem~\ref{thm:main2} given below.]

\newcommand{\xrig}{X^{rig}}

\para Let $X/E$ be a geometrically connected, projective variety over a $p$-adic field and $E\into K$ an isometric embedding into a complete valued field with a rank one valuation. Then one has the (projective) analytic spaces $\xan/E$ and $\xan/K$. Projectivity (though not essential for my argument) ensures,  by \cite[Chap 3]{berkovich-book}, that a number of adjectives which may be applied to an analytic space, can be applied to  both of these spaces: both are proper (hence separated, so quasi-separated), strict, good, compact (hence quasi-compact) i.e. covered by a finite number of affinoid open subsets and the construction below applies to analytic spaces which enjoy some of these properties (but not necessarily projectivity). 

By definition of an analytic space, $\xan_K$ is equipped with an  atlas of affinoid opens. This data can be used to equip $\xan_K$ with a sheaf of Banach algebras $\O_{\xan_{K}}$ (to be precise this means that for any quasi-compact open, $U\subset \xan_{K}$, the algebra $\O_{ \xan_{K}}(U)$ is a Banach algebra which is functorial in such $U$ with the following properties: (1) if $U=\sM(A)$ is an affinoid open then $\O_{ \xan_K}(U)=A$ and (2)  if $U$ is any quasi-compact open with $U=\cup_i U_i$  a finite cover by affinoids then $\O_{ \xan_{K}}(U)\to \prod_i \O_{ \xan_{K}}(U_i)$ is a closed embedding of Banach algebras (note that the Banach norms provided in this construction are not claimed to be unique (locally) but equivalent). This construction is detailed  in \cite[Section 3.3.2, Section 4.1.2]{temkin15}. \emph{The important point  here is not the sheaf itself, but the fact that the spaces of  local analytic functions acquire a Banach space structure, which agrees with the norm on constant functions i.e. on our field $K$, in a manner that is compatible with gluing of local analytic functions and independent of the  gluing data.} The most succinct way of expressing all this is to say that one has a sheaf of Banach algebras $\O_{ \xan_{K}}$ on $\xan_{K}$ for a suitable Grothendieck topology on $\xan_K$. \emph{This implies, in particular, that the ring $\Gamma(\xan_K,\O_{\xan_K})$ of global analytic functions on $\xan_K$ is naturally a Banach algebra. }

\para Let me briefly sketch a proof of the above claims. Readers familiar with the construction of such a sheaf may skip this paragraph. By definition, an analytic space is equipped with an atlas of affinoid open subsets and some gluing data and the analytic space can be equipped with a Grothendieck topology using this data i.e. one restricts the notion of open subsets for the purpose of constructing sheaves. The construction of $\O_{\xan_{K}}$ uses this datum. The key tool in the construction of $\O_{\xan_K}$ is Tate's Acyclicity Theorem \cite[3.3.2.1]{temkin15}.
The following general facts about Banach algebras and Banach modules over Banach algebras will be useful to remember: 
\benumlab
\item If $A$ is any $K$-affinoid algebra then $A$ is equipped with a  norm (and even a power multiplicative norm, if one assumes additionally that $A$ is reduced, which is certainly true in the case which I am concerned with here \cite[6.2.4, Theorem 1]{bosch-non-arch-analysis}, but the existence of some norm on $A$ can always be inferred from the Gauss norm), equipping $A$ with a structure of  a Banach algebra (i.e.  $A$ is complete with respect to this norm), and any two norms on $A$ are equivalent and moreover the  restriction of this norm  to  $K\into A$ is the valuation norm on $K$. If $K\supseteq E$ is a complete valued subfield, then one can think of a $K$-affinoid algebra as an $E$-Banach algebra.
\item If $A,B$ are Banach $E$-algebras then any $E$-linear homomorphism $f:A\to B$ is continuous if and only if it is bounded \cite[2.1.8, Proposition 2]{bosch-non-arch-analysis}. 
\item Product of  $E$-Banach algebras $A_1,\ldots, A_n$, is also a $E$-Banach algebra with the obvious definition of a norm.
\eenum
The sheaf $\O_{\xan_K}$ is constructed as follows (see \cite[Section 3.3.2, Section 4.1.2]{temkin15}). First consider the case of an affinoid open set.  If $U=\sM(A)\subset \xan_K$ is an affinoid open subset then $\O_{\xan_K}(U)=A$ is evidently a Banach algebra over $E$.  If $U$ is  covered by a finite number of affinoid opens $U=\cup_{i=1}^nU_i$, with $U_i=\sM(A_i)$ and $U_i\cap U_j=\sM(A_{i,j})$ then using Tate's Acyclicity Theorem \cite[8.2.1, Theorem 1]{bosch-non-arch-analysis} or \cite[3.3.2.1]{temkin15}) one obtains the equality of Banach algebras $A=\ker(\prod A_i \to \prod A_{i,j})$ so one can indeed define the sheaf $\O_{\xan_{K}}$ using the rule $\O_{\xan_K}(U)=A$ on affinoids (provided in the atlas and the net of affinoids defining the analytic space $\xan_{K}$). Moreover, Tate's acyclicty theorem also shows, that this gives  a sheaf of Banach algebras on $U=\sM(A)$ which  is independent of the choice of the covering. 

Now suppose   $U\subseteq \xan_{K}$ is an arbitrary quasi-compact open subset of $\xan_K$.  Choose a finite covering of $U=\{U_i\}_{i=1}^n$  by affinoids with $U_i=\sM(A_i)$, with $U_i\cap U_j=\sM(A_{i,j})$. Then  $\O_{\xan_{K}}(U)$ is the equalizer of the two restriction  arrows
 $$\xymatrix{
\prod A_i	\ar@<0.5ex>[r]\ar@<-0.5ex>[r] &  \prod_{i,j} A_{i,j}.
}
$$
So $\O_{\xan_{K}}(U)$ is closed in the product and hence carries a natural structure of Banach algebra. This is independent of the choice of the covering: any two such covers of $U$ have a common refinement and  provide isomorphisms between the three possible Banach structures on $\O_{\xan_K}(U)$  via  \cite[2.8.1 Banach's Open Mapping Theorem]{bosch-non-arch-analysis}. So one gets independence of coverings and also natural compatibility of the Banach norms on $\O_{\xan_K}(U)$.

Since $\xan_K$ is covered by a finite number of affinoids, this constructs  $\O_{X^{an}_K}$ as a sheaf of $E$-Banach algebras for the Grothendieck topology of $\xan_{K}$ given by the net of compact analytic domains in $X$. Moreover, affinoid locally on $\xan_K$, the norm on $\O_{\xan_K}$, on constant functions $K$, coincides with the valuation.   The construction  of the sheaf of Banach algebras $\O_{\xan_{K}}$ is functorial for morphisms of analytic spaces described above.  Moreover one also sees from this local description that the ring  of global functions  $\Gamma(\xan_K,\O_{\xan_K})=H^0(\xan_K,\O_{\xan_K})$, on $\xan_K$,  is a  Banach $E$-algebra and on the constant functions, this  norm is equivalent to the one given by the valuation.

\para Now let me prove the following important observation:
\bthm\label{thm:main2}
Let $X/E$ be a geometrically connected, smooth projective variety. Let $K_1,K_2$ be two untilts of $\cpt$ which contain $E$. Suppose that $K_1,K_2$ are not topologically isomorphic. Then
\benumlab
\item\label{thm:main2-1} one has  an isomorphism of topological groups $$\pit{X/K_1}\isom \pit{X/K_2},$$
\item\label{thm:main2-2} but the analytic spaces  $\xan/{K_1}$ and $\xan/{K_2}$  are not isomorphic  as analytic spaces over $\Q_p$ (in the sense of \ssep\ref{pa:E-isomorphisms-analytic-spaces}).
\item\label{thm:main2-3} In particular the Absolute Grothendieck Conjecture fails in the category of Berkovich spaces over perfectoid fields of characteristic zero.
\eenum
\ethm
\bp 
After Theorem~\ref{thm:main},
only \ref{thm:main2-2}  needs to be proved as \ref{thm:main2-2} $\implies$ \ref{thm:main2-3}. The hypothesis of Theorem~\ref{thm:main2} are non-vacuous--by \cite{kedlaya18}, fields $K_1,K_2$ exist. 

Assume that $X/E, K_1,K_2$ are as in my hypothesis and that $X$ is geometrically connected, smooth and projective over $E$. Suppose, if possible, that $\xan/K_1$ and $\xan/K_2$ are isomorphic as analytic spaces over $\Q_p$. Then one has a bounded isomorphism of Banach rings $$K_1\isom H^0(\xan/K_1,\O_{\xan/K_1}) \isom H^0(\xan/K_2,\O_{\xan/K_2})\isom K_2.$$ 
I claim that this is in fact an isomorphism of valued fields. Write $K_1^\circ$ (resp. $K_2^\circ$) for the respective subrings of power bounded elements of $K_1$ (resp. $K_2$) (for the respective norms). By \cite[1.2.5, Proposition 4]{bosch-non-arch-analysis}, the above isomorphism induces an  isomorphism of $K_1^\circ\mapright{\isom} K_2^\circ$. Further, as the norm $\abs{-}_{K_i}$ on  $K_i$ arises from the valuation of these fields, so one sees that the norms are power-multiplicative. This implies, by \cite[1.3.1, Proposition 4]{bosch-non-arch-analysis}, that one has the equality $$K_i^\circ=\{x\in K_i: \abs{x}_{K_i} \leq 1\},$$ i.e. $K_i^\circ$ is the valuation subring of $K_i$ and thus the valued fields $K_i$ have isomorphic valuation rings and hence $K_1$ and $K_2$ are therefore isomorphic as valued fields. 
Thus one has arrived at a contradiction.
\ep

\brem 
For the existence of uncountably many topologically non-isomorphic perfectoid fields see the construction of inequivalent valuations given by \cite{schmidt33} (sketched here in \Cref{pa:existence-perfectoid}) or \cite[Theorem 8]{kaplansky42}.
\erem

\brem\label{re:uncountable-copies}
As an aside let me remark that the proof of \cite{kedlaya18} (also see \cite[Th\'eor\`eme 2 and \S3 Remarque 2]{matignon84}) provides an uncountable collection of perfectoid fields $K_1,K_2$ with tilts isometric to $\cpt$ and such that  $K_1,K_2$ are not topologically isomorphic.
\erem

\para\label{pa:anabelomorphisms} Let me introduce some terminology from \cite{joshi-anabelomorphy}. I will say that two geometrically connected varieties $X/E$ and $X'/E'$ over fields $E,E'$ are \emph{anabelomorphic} (resp. \emph{tempered anabelomorphic}) if one has a topological isomorphism of their \'etale fundamental groups (resp. tempered fundamental groups if $E,E'$ are $p$-adic fields):
$$\alpha:\pi_1(X/E)\isom \pi_1(X'/E') \text{ resp. } \alpha:\pit{X/E}\isom \pit{X'/E'}$$
and in this situation I will write $\alpha:X/E\anabelmap X'/E'$ for this \emph{anabelomorphism} (resp. \emph{tempered anabelomorphism}). I will say that an anabelomorphism (resp. \emph{tempered anabelomorphism}) $\alpha:X/E\anabelmap X'/E'$ is a \emph{strict anabelomorphism} (resp. \emph{strict tempered anabelomorphism}) if $X/E, X'/E'$ are anabelomorphic but not isomorphic (resp. anabelomorphic but not isomorphic  analytic spaces). 

In this terminology, Theorem~\ref{thm:main2} says there exist perfectoid fields $K,K'\supset E$ and a strict tempered anabelomorphism  \be 
\xan/K\anabelmap \xan/K'.\ee

\emph{Note that anabelomorphism defines an equivalence relation on geometrically connected, smooth, quasi-projective varieties and it makes perfect sense to talk about the anabelomorphism class of a variety.}

Let me remark that in \cite{joshi-anabelomorphy} I show that anabelomorphy of $p$-adic fields changes important invariants of $p$-adic fields such as discriminants and more importantly it also impacts geometric invariants of varieties such as minimal discriminants of elliptic curves. 

Another notion introduced in \cite{joshi-anabelomorphy} is that of amphoricity: a quantity, a property or an algebraic structure associated to $X/E$ is said to be \emph{amphoric} if it is an invariant of the anabelomorphism class of $X/E$.

\para\label{pa:question-about-filtrations-on-Pi}
As is well-known from \cite{mochizuki-local-gro}, a $p$-adic field $E$ is not amphoric \cite{mochizuki-local-gro} i.e. $G_E$ does not determine the isomorphism class of $E$, but as was shown in \cite{mochizuki-local-gro}, $G_E$ equipped with its upper numbering ramification filtration $G_E^\mydot$ determines $E$. Now let $X/E$ be a geometrically connected, smooth, quasi-projective variety. Let $\Pi=\pit{X/E}$ and let $\Pi\supset\overline{\Pi}$ be the geometric tempered fundamental group. Mochizuki has shown that the quotient $\Pi_{X/E} \to G_E$ is amphoric \cite[Lemma 1.3.8]{mochizuki04} and hence its kernel $\overline{\Pi}\subset \Pi$  is an amphoric subgroup i.e. determined by the isomorphism class of $\Pi=\pit{X/E}$. In the light of this and  Theorem~\ref{thm:main} and Theorem~\ref{thm:main2} one can ask the following question:
\begin{question}
	Let $K$ be a complete, algebraically closed valued field containing an isometric embedding of $E$.  Is there some filtration by normal subgroups $\Pi^\mydot_{X}\subset \Pi_X$ which determines the pair of analytic spaces $(\xan_E, \xan_K)$ up to an isomorphism? 
\end{question}
Mochizuki's Theorem that $\overline{\Pi}$ is amphoric should be considered as the analog of the assertion (of \cite{mochizuki-local-gro}) that the inertia subgroup $I_E\subset G_E$ is amphoric. If $G\isom G_E$ is an isomorph of $G_E$ then one may equip $G$ with many different inertia filtrations corresponding to anabelomorphisms $G_E\isom G\isom G_{E'}$. Similarly if $\Pi$ is an isomorph of $\Pi_{X/E}$ then above suggests remarkably that there are many different filtrations $\Pi^\mydot$ each corresponding to an algebraically closed, complete valued field $K\supset E$. At least when $X/E$ is a hyperbolic curve this question is quite reasonable.

\numberwithin{equation}{subsection}

\section{Arithmetic Holomorphic Structures}\label{se:arith-hol-strs}
The results of Section~\ref{se:grothendieck-conj} namely, Theorem~\ref{thm:main} and Theorem~\ref{thm:main2} provide us with a way of defining Arithmetic Holomorphic Structures (see \Cref{def:arith-hol-strs}) and allow one to assert with impunity that this notion of Arithmetic Holomorphic Structures is fully compatible with our conventional way of  thinking about distinct classical holomorphic structures on Riemann surfaces. As I note in Theorem~\ref{thm:arith-hol-strs}, my definition of arithmetic holomorphic structures also provides arithmetic holomorphic structures in sense this term is used in \iut. In particular, in the $p$-adic contexts one can retain both:  our ideas from classical theory complex holomorphic functions and also preserve ideas from  anabelian geometry. Thus this work provides a natural solution to the problem (raised in the context of \iut) of producing \emph{geometrically, and canonically labeled copies of the tempered fundamental groups} in a manner compatible with the classical case (\Cref{th:main-classical}). This aspect is detailed in \Cref{se:untilts-of-Pi}.

While the term Arithmetic Holomorphic Structures was coined and used by Mochizuki in \iut,	 my definition \Cref{def:arith-hol-strs} is far more precise and Theorem~\ref{thm:main2}(2) shows that  it is a holomorphic structure in the sense of holomorphic functions, secondly this holomorphic structure   is dependent on the arithmetic and topological properties of the coefficient field and hence the nomenclature `Arithmetic Holomorphic Structure' is quite apt.   Arithmetic Holomorphic Structures defined here also provide arithmetic holomorphic structures in Mochizuki's sense. \textit{This is the starting point of the Teichm\"uller Theory presented in the \present.} [Mochizuki does not take this Berkovich analytic space approach which I have taken here.]

Once one has defined Arithmetic Holomorphic Structures, by imitating the classical theory, one defines the Arithmetic Teichm\"uller Space (\Cref{def:arith-hol-space-local}) as the category of all Arithmetic Holomorphic Structures with a fixed tempered fundamental group. 

Let me also point out that the approach taken here for the $p$-adic case involves deforming (infinite dimensional) Banach spaces of (Berkovich) analytic functions. As should be clear from \Cref{app:archimedean-theory2}, especially \Cref{th:nakai1}, that  Banach space approach  also exists in the classical case (but not so well-known--I discovered Nakai's works after \cite{joshi-untilts-2020} was posted on the arxiv). At any rate, one now has a complete parallel between the $p$-adic and archimedean theory.

\subsection{Definitions}
Let $\sM(A)$ be the Berkovich spectrum \cite{berkovich-book} of any Banach algebra $A$, and let $\xan_E,\xan_K$ be the analytic spaces arising from $X/E$ etc. (as before $X/E$ is a geometrically connected, smooth, quasi-projective variety over a $p$-adic field $E$). It will be useful to recall the definition of geometric base-points in the context of \cite{andre-book}.

\begin{defn}\label{def:arith-hol-strs}
	Let $X/E$ be a geometrically connected, smooth, quasi-projective variety over a $p$-adic field $E$.	
	A (pointed) \textit{arithmetic holomorphic structure} $$(X/E,(K\supset E,K^\flat\isom F),*_K:\sM(K)\to\xan_E)$$ on $X/E$ consists of
\benumlab
\item  a choice of an untilt $(K\supset E,K^\flat\isom F)$ for some algebraically closed, perfectoid field $F$ of characteristic $p>0$, and 
\item a choice of a morphism of Berkovich analytic spaces $*_K:\sM(K)\to\xan_E$. 
\eenum
\textit{Morphisms of arithmetic holomorphic structures} are be defined as morphisms of the data $(X/E,(K\supset E,K^\flat\isom F),*_K:\sM(K)\to\xan_E)$ of the arithmetic holomorphic structure.
\end{defn}

\begin{defn}\label{de:geom-base-points} 
Following see \cite[Chapter III, 1.2.2]{andre-book}, a \textit{$K$-geometric base-point} (or less precisely, a \textit{geometric base-point}) of the analytic space $\xan_E$   is an algebraically closed, complete (rank one) valued field $K$ and an isometric extension  of valued fields $K/E$ and a morphism of analytic spaces $$*_K:\sM(K)\to \xan_E$$. 
\end{defn}

\brem\ 
\benumlab 
\item Let me remark that   \cite{andre-book} works in the context of  \cite{berkovich-book}, and hence any geometric base-point as defined here and in \cite[Chapter III, 1.2.2]{andre-book} provides an algebraically closed, valued field which is complete with respect to a rank one valuation. Hence, especially that, a $K$-geometric base-point in \Cref{de:geom-base-points}  provides an algebraically closed, complete (rank one) valued field $K$ containing $\Q_p$ isometrically. Thus, by Lemma~\ref{lem:perfectoid}, a geometric base-point provides  an algebraically closed, perfectoid  field $K\supset E\supset \Q_p$.  
\item The discussion of \iutf{IUT1, \S I3, Page 21} makes it clear that \iut\ requires arbitrary geometric base-points (and so the valued field $K\supset E$ must be arbitrary). \textit{This is how algebraically closed, perfectoid fields enter \iut\     (see \cite{joshi-teich-IIIB-rosetta}).} 
\eenum
\erem  

\blem\label{le:geom-base-point} 
Each arithmetic holomorphic structure $$(X/E,(K\supset E,K^\flat\isom F),*_K:\sM(K)\to\xan_E)$$ is equipped with a $K$-geometric base-point, namely the morphism 
$*_K:\sM(K)\to\xan_E$.
\elem  
\bp 
This is clear from the definitions.
\ep

Here is one immediate consequence of \Cref{def:arith-hol-strs} and \Cref{le:geom-base-point}.
\bpro\label{pr:galois-seq-untilt} Let $(X/E; (K\supset E,K^\flat\isom F),*_K:\sM(K)\to \xan_E)$ be an arithmetic holomorphic structure on $X/E$. Then the geometric base-point provided by \Cref{le:geom-base-point} provides the exact sequence of topological groups 
\be\label{eq:galois-seq-untilt} 1\to \pit{X/K} \to \pit{X/E}\to G_{E;K}\to 1.\ee
The isomorphism class each of the groups and the isomorphism class of this exact sequence of topological groups is independent of the choice of the arithmetic holomorphic structure on $X/E$ and also independent of the choice geometric base-point. 
\epro 
\bp 
This is immediate from \Cref{thm:main} and properties of tempered fundamental groups established in \cite{andre03}.
\ep

Let me record the following which makes it clear that distinct i.e. non-isomorphic arithmetic holomorphic structures exist.

\bpro\label{le:non-isom-arith-hol-strs}
Let $F$ be any algebraically closed perfectoid field of characteristic $p>0$. Let $((K_1\supset E,K_1^\flat\isom F),*_{K_1}:\sM(K_1)\to\xan_E)$, $((K_2\supset E,K_2^\flat\isom F),*_{K_2}:\sM(K_2)\to\xan_E)$ be two arithmetic holomorphic structures. Then the following assertions hold:
\benumlab
\item The valued fields $K_1,K_2$ have characteristic zero, the same cardinality, the same residue fields and the same value groups as that of $F$.
\item There exist arithmetic holomorphic structures $((K_1\supset E,K_1^\flat\isom F),*_{K_1}:\sM(K_1)\to\xan_E)$, $((K_2\supset E,K_2^\flat\isom F),*_{K_2}:\sM(K_2)\to\xan_E)$ which are not isomorphic.
\eenum 
Furthermore, if $F=\cpt$, then $K_1,K_2$ need not even be topologically isomorphic fields (and hence need not be isomorphic as valued fields) and, if $X/E$ is projective, then by \Cref{thm:main2}, $\xan_{K_1}$ and $\xan_{K_2}$ need not be $\Q_p$-isomorphic analytic spaces.
\epro
\bp
The first assertion is an immediate consequence of \cite{fargues-fontaine}, \cite{scholze12-perfectoid-ihes}. By \cite[Corollaire 2.2.22]{fargues-fontaine}, one sees that for any given perfectoid field $F$ of characteristic $p>0$ and any given $p$-adic field $E$, there exist many non-isomorphic untilts and this together with \Cref{def:arith-hol-strs}, implies that there exist many non-isomorphic arithmetic holomorphic structures. The last assertion for $F=\cpt$ is a consequence of \Cref{thm:main2} and  the fundamental theorem of \cite{kedlaya18}.
\ep

\brem 
Theorem~\ref{thm:main} makes it clear that one is dealing with holomorphic structures in the familiar sense of complex analytic function on the other hand, this theorem and \Cref{le:non-isom-arith-hol-strs} also makes it clear that this phenomenon arises from topological changes in arithmetic of the fields one is dealing with here. This is the reason why the term Arithmetic Holomorphic Structures is appropriate here.
\erem

\brem 
Let me remark importantly that by \cite{kedlaya18} that there exist  algebraically closed, perfectoid fields with an isometrically embedded $\Q_p$ and with tilts isometric to $\cpt$, but such that $\bQ_p$ is not necessarily dense in these fields. So one cannot hope to work (solely) over $\bQ_p$ once one recognizes that arbitrary algebraically closed perfectoid fields are required in \iut\ via its requirement of arbitrary geometric base-points. 
\erem

\bpro\label{pr:various-functors}
From the category whose objects are of the form $(X/E,(K\supset E,K^\flat\isom F),*_K:\sM(K)\to\xan_E)$, there are natural forgetful functors, to various categories, given as follows:
$$
\begin{tikzcd}
	& \pi_1^{temp}(\xan_E,*_K: \sM(K)\to\xan_E)&&\\
	& (X_E^{an}, *_K:\sM(K)\to\xan_E)\arrow[u,mapsto]&&\\
	&\arrow[dl,mapsto] \arrow[d,mapsto]\arrow[u,mapsto] (X/E,(K\supset E,K^\flat\isom F),*_K:\sM(K)\to\xan_E)\arrow[dr,mapsto]&&\\
	X/E& (K\supset E,K^\flat\isom F)\arrow[dl,mapsto]\arrow[dr,mapsto] & X_K^{an}\to X_E^{an} \\
	 K &&  K^\flat
\end{tikzcd}
$$
\epro
\bp 
The functor
$$(X/E,(K\supset E,K^\flat\isom F),*_K:\sM(K)\to\xan_E)\mapsto \left( X_K^{an}\to X_E^{an}\right)$$ is obtained by considering the pair of spaces $(\xan_E,\xan_K)$ and the morphism $\xan_K\to\xan_E$ obtained by base change via  the given isometric inclusion $E\into K$. 

The morphism $*_K:\sM(K)\to\xan_E$ is a geometric base-point of the analytic space $\xan_E$ and the functor 
$$(X_E^{an}, *_K:\sM(K)\to\xan_E)\mapsto \pi_1^{temp}(\xan_E,*_K: \sM(K)\to\xan_E)$$ is the tempered fundamental group functor constructed in \cite{andre03}   using this geometric base-point given in the datum of the object $(X/E,(K\supset E,K^\flat\isom F),*_K:\sM(K)\to\xan_E)$. 

The functor  $$(X_E^{an}, *_K:\sM(K)\to\xan_E)\mapsto  \pi_1^{temp}(\xan_E,*_K: \sM(K)\to\xan_E)$$ is thus  the composite functor obtained as
$$(X_E^{an}, *_K:\sM(K)\to\xan_E)\mapsto  (\xan_E, *_K:\sM(K)\to\xan_E) \mapsto \pi_1^{temp}(\xan_E,*_K: \sM(K)\to\xan_E).$$

The other functors are clear.
\ep

\brem
Suppose one is given any $K$-geometric base point $*_K:\sM(K)\to\xan_K$, then
by composition with the base-change morphism $\xan_K\to\xan_E$ provided by the given isometric embedding $E\into K$, one obtains  a $K$-geometric base-point of $*_K:\sM(K)\to \xan_E$ (note the conflation of notation). In particular, in defining $\fjxe$, one can work instead with the datum $$((K\supset E,K^\flat\isom F), *_K:\sM(K)\to \xan_K)$$ instead of $$((K\supset E,K^\flat\isom F), *_K:\sM(K)\to \xan_E).$$
\erem

\subsection{Notational Contractions} In this paper, for simplicity of notation, I may contract the notation 
$$(X/E; (K\supset E,K^\flat\isom F),*_K:\sM(K)\to \xan_E)$$ to $$(X/E; (K\supset E,K^\flat\isom F)).$$
This may give the impression that the geometric base-point datum is not) is unimportant, but let me caution the reader that this is certainly not the case. I may also further contract the notation to
 $$(X/E; E\into K).$$
   This may also give the impression that the tilting data (i.e. the isometry $K^\flat\isom F$ and the geometric base-point datum is not) is unimportant, but let me caution the reader that this is certainly not the case. 
   
\brem\  
\benumlab
\item   The tilting datum $K^\flat\isom F$ provided by an untilt allows us to compute arithmetic degrees in one fixed location namely the value group of $F$ even as the untilt moves (see \cite{joshi-teich-estimates} for how this gets used and its relevance for \iutf{IUT3, Corollary 3.12}). 

\item As is established in \constrthrB{\ssep }{ss:mochizuki-dichotomy}, in the parlance of \iut, the perfectoid field $K$ provides ``\'etale-Picture''  given by group theoretic data of the isomorphs $\pit{X/E;K}$ and $G_{E;K}$ of the topological groups $\pit{X/E}$ and $G_E$ respectively.
While the tilting data is ``Frobenius-Picture.'' In fact, after \Cref{def:arith-hol-strs},  Mochizuki's term ``Frobenius-like'' can be understood in the literal sense as  the construction of $K^\flat$ requires $p^{th}$-powers i.e. Frobenius morphism.
\eenum
\erem

\subsection{Relationship with the Diamond $X_E^{\Diamond}$}\label{re:diamonds} Let $X/E$ be a geometrically connected, smooth, quasi-projective variety over a $p$-adic field $E$. Write $X_E^{\Diamond}=(\xan_E)^{\Diamond}$ for the diamond associated to the analytic space $\xan_E$ by \cite[Section 15]{scholze-diamonds}. 
In this section I want to record an important observation which makes it clear that the theory of the \present, the theory of tempered fundamental groups (and hence \iut) are inseparably linked to Scholze's Theory of Diamonds described in \cite{scholze-diamonds}. This remark was first pointed out in \cite{joshi-untilts}.

\bpro\label{pr:local-diamond} 
There is a natural functor from the category of objects of the form $$(X/E,(K\supset E,K^\flat\isom F),*_K:\sM(K)\to\xan_E)$$ to the diamond $X^\Diamond_E$ given by the rule

\begin{tikzcd}
 (X/E,(K\supset E,K^\flat\isom F),*_K:\sM(K)\to\xan_E)\arrow[d,mapsto]\\  ((K\supset E,K^\flat\isom F),*_K:{\rm Spa}(K,\O_K)\to X^{ad}_E)\in X^\Diamond_E(F).
\end{tikzcd}

\epro
\bp 
This is immediate from the results of \cite{huber94}, \cite[Section 15]{scholze-diamonds}:  the $K$-geometric base-point i.e. the morphism $*_K:\sM(K)\to\xan_E$ of Berkovich analytic spaces provides, by the means of the fully-faithful functor from Berkovich analytic spaces  to adic spaces constructed by \cite[Proposition 4.3 and Proposition 4.5]{huber94}, to a morphism $*_K:{\rm Spa}(K,\O_K)\to X^{ad}_E$ and conversely. 
\ep

\section{Construction of Arithmetic Teichm\"uller Spaces}\label{se:construct-att}
\subsection{Definition in the non-archimdean case} Classical Teichm\"uller, \Cref{thm:main}, the definition of Arithmetic Holomorphic Structures  on $X/E$ \Cref{def:arith-hol-strs} and \Cref{thm:arith-hol-strs} suggests a natural definition of Arithmetic Teichm\"uller Space of $X/E$: namely, this space should be the category whose objects are arithmetic holomorphic structures all of which have tempered fundamental groups isomorphic to the tempered fundamental group of $X/E$.   Here is a more precise version:

\begin{defn}\label{def:arith-hol-space-local}
	Let $E$ be a $p$-adic field and let $X/E$ be a geometrically connected, smooth, quasi-projective variety over $E$. Then the \textit{Arithmetic Teichm\"uller Space $\fjxe$ of $X/E$} is a category whose objects are
	arithmetic holomorphic structures
	$$	(Y/E',(E'\into K, K^\flat \isom F), *_K:\sM(K)\to \yan_{E'})$$ 
	where 
	\benumlab 
	\item $E'$ is  a $p$-adic field,
	\item  $F$ is some algebraically closed perfectoid field of characteristic $p>0$,
    \item  $(E'\into K, K^\flat \isom F)$ is an   $E'$-untilt of $F$,
	\item	$Y/E'$ is a geometrically connected, smooth, quasi-projective variety such that
	\item one has	$$\dim(Y)=\dim(X),$$ 
	\item and
	\begin{enumerate}
		\item 	$Y/E'$ is tempered anabelomorphic to  $X/E$ i.e. one has an anabelomorphism
			$$	\pit{Y/E'}\isom\pit{X/E},$$ 
		\item and any such anabelomorphism induces an anabelomorphism of their geometric tempered fundamental subgroups
		      $$	\pit{Y/K}\isom\pit{X/\C_p}.$$
	\end{enumerate}
\eenum
Morphisms in $\fjxe$ are morphisms of the data $$(Y/E',(E'\into K, K^\flat \isom F), *_K:\sM(K)\to \yan_{E'}).$$
\end{defn}

\brem
In the principal case of hyperbolic curves over $p$-adic fields which occurs in the Anabelian Geometry literature, the requirement {\bf(4)}(b) is implied by {\bf(4)}(a) (by \cite[Lemma 1.3.8]{mochizuki04}). I do not know if this is still true in the general case considered here.
\erem

\brem\label{pa:kedlaya-remark} In recent correspondence, Kiran Kedlaya pointed out to me the following consequence of \cite[Proposition 8.8.9]{kedlaya-liu19}: \emph{every deformation of an analytic space (arising from the analytification of a  quasi-projective variety) over a perfectoid field arises from a deformation of the perfectoid field}. So the idea of moving the algebraically closed perfectoid field (considered here) is, in a rather precise sense (of \cite{kedlaya-liu19}), optimal. \erem

\subsection{Definition in the Archimedean Case} \textcolor{red}{This section still needs some cleanup in case $\dim(X)>1$.} This section may be skipped by readers who are interested only in the case of smooth, hyperbolic curves i.e. the dimension one case: in the dimension one case one simply takes classical Teichm\"uller space \cite{imayoshi-book} as the Archimedean Teichm\"uller Space.

Now suppose that $E$ is an algebraically closed, archimedean complete local field. So $E\isom \C$ and  think  $X$ is a complex, smooth, quasi-projective variety. I will refer to this as the archimedean case. I will write $\xan_\C$ for the complex analytic space i.e. the complex manifold given by the scheme $X$.

In the archimedean case, if the complex variety $X/\C$ has $\dim(X/\C)=1$, then classical Teichm\"uller Theory provides a robust theory, but if $\dim(X/\C)>1$ then there is no good definition of (archimedean) Teichm\"uller space and hence one has some latitude in defining what it should be. 

Here I provide a reasonable working definition of (archimedean) Teichm\"uller space which works for all $\dim(X/\C)\geq 1$ and which coincides with the classical case of $\dim(X/\C)=1$ (Riemann surfaces). Further work is, of course needed, if $\dim(X/\C)>1$. 

The notion of quasi-conformal mappings extends to higher dimensional Riemannian manifolds, but on the other hand for Riemannian manifolds of dimension greater than two, the notion of quasi-conformal mappings is not adequate for a satisfactory theory. So an alternative approach is needed to build the theory.

My approach here is based on \cite{nakai59}, \cite{nakai60}, \cite{nakai1972} and \cite[Appendix: Higher Dimensions]{nakai-book}.  As is detailed in \Cref{app:archimedean-theory2}, for $\dim(X/\C)=1$, the results of \cite{nakai59}, \cite{nakai60} using Royden algebras of functions are closest to the Banach space approach of this paper and allow us to view classical Teichm\"uller Theory as a functor in infinite dimensional Banach algebras (as is done here in the $p$-adic case). For $\dim(X/\C)\geq 2$,  \cite{nakai1972} suggests a natural generalization (again using results of \cite{nakai1972} on Royden algebras) which is compatible with $\dim(X/\C)=1$ case. In particular, for all dimensions $\dim(X/\C)\geq1$ one has a good category equipped with a functor to infinite dimensional Banach algebras.

To adopt \cite{nakai1972} approach one must work with Riemannian manifolds. Write $\xan=\xan_{\C}$ for the associated complex manifold. Since $X/\C$ is a  complex, smooth, quasi-projective variety, $\xan$ can be equipped with a Riemannian metric. For example, using any embedding $X\into \P^N$ in projective space and considering the restriction of the Fubini-Study metric on $\P^N$ to $\xan$, one sees that $\xan$ is a Riemannian manifold. In this subsection, I will view $\xan$ for any complex, smooth, quasi-projective variety as a Riemannian manifold.

A \textit{Nakai quasi-isometry} between Riemannian manifolds is a quasi-isometry of Riemannian manifold as defined in \cite[Paragraph 2]{nakai1972}. A Nakai quasi-isometry is a quasi-conformal mapping and also a homeomorphism. By \cite{nakai1972}. If $\dim(X/\C)=1$    (i.e. the underlying real manifold $\xan_{\C}$ is of dimension two) then one sees that a Nakai quasi-isometry is the same as a quasi-conformal mapping. Hence one recovers the classical Teichm\"uller space if $\dim(X/\C)=1$. For $\dim(X/\C)\geq 1$, the notion of Nakai quasi-isometry appears to be the correct generalization of quasi-conformal mappings which is compatible with the Banach space approach of this paper by the results of \cite{nakai1972}.

\begin{defn}\label{def:arith-hol-space-local-arch}
	Let $X/\C$ be a complex, smooth, quasi-projective variety, considered as being equipped with a Riemannian metric. The archimedean \textit{pre-Teichm\"uller Space  of $X/\C$} consists of pairs $(\yan,f:\yan\to \xan)$ where $Y$ is a connected, smooth, quasi-projective variety over complex numbers and $f:\yan\to \xan$ is a Nakai quasi-isometry of Riemannian manifolds. Note that this includes pairs of the form $(\xan,f:\xan\to \xan)$ and hence the set of such pairs is not vacuous as soon as $X\neq \emptyset$. 
\end{defn}
The Teichm\"uller space $\fJ(X,\C)$ of $X$ is then the pre-Teichm\"uller Space of $X/\C$ subject to a suitable equivalence. 
\begin{defn}
Two pairs $(\yan,f:\yan\to \xan)$ and $(Z^{an},g:Z^{an}\to \xan)$ as in \Cref{def:arith-hol-space-local-arch} are equivalent if $f\circ g^{-1}$ is homotopic to a mapping which induces an isometry between Royden algebras of $\yan$ and $Z^{an}$. The archimedean Teichm\"uller Space of $X/\C$ is the set of equivalence classes of pairs $(\yan,f:\yan\to \xan)$ in the pre-Teichm\"uller Space of \Cref{def:arith-hol-space-local-arch}.
\end{defn}

\brem
By \cite{nakai60}, if $\dim(X/\C)=1$, then two pairs $(\yan,f:\yan\to \xan)$ and $(Z^{an},g:Z^{an}\to \xan)$ as in \Cref{def:arith-hol-space-local-arch} are equivalent if $f \circ g^{-1}$ is homotopic to a conformal mapping. So $\fJ(X,\C)$ is the classical Teichm\"uller Space if $\dim(X/\C)=1$.
\erem

The following result is immediate from \cite{nakai1972}:
\bpro 
Suppose $f:\yan\to \xan$ is a Nakai quasi-isometry of Riemannian manifolds. Then the real manifolds $\yan,\xan$ are homeomorphic and hence have the same real dimension and hence one has $\dim(X/\C)=\dim(Y/
\C)$ and have isomorphic topological fundamental groups.
\epro

The next theorem is my reformulation of  main theorem of \cite{nakai1972} it provides an algebraic characterization of Nakai quasi-isometries similar to the case of Riemann surfaces discussed in \Cref{app:archimedean-theory2}:
\bthm\label{th:nakai2}
Let  $X,Y$ be connected, smooth, quasi-projective varieties over $\C$. Let $\Sigma=\xan_{\C}$, $\Sigma'=\yan_{\C}$ be the underlying complex manifolds considered as being equipped with a Riemannian metrics. Then 
\benumlab
\item A homeomorphism $\Sigma'\to \Sigma$ is a Nakai quasi-isometry if and only if one has a $\C$-algebra isomorphism $$\rs\isom \rsp$$  between their Royden algebras. 
\item Any $\C$-algebra isomorphism $\rs\isom \rsp$ is a homeomorphism for the respective norm-topologies.
\item The construction of $(\rs,\abs{-}_\Sigma)$ is functorial for Nakai quasi-isometries between Riemann manifolds.
\item In particular one has a functor $$(\Sigma',f:\Sigma'\to \Sigma)\mapsto (\rsp,\abs{-}_{\Sigma'})$$ from the Teichm\"uller category $\sTs$  to the category of infinite dimensional complex Banach algebras.
\item Moreover, if $\Sigma,\Sigma'$ are Nakai quasi-isometric then
\begin{enumerate}
	\item  $\dim(X/\C)=\dim(Y/\C)$ as complex varieties, and
	\item $\pi_1^{top}(\Sigma')\isom \pi_1^{top}(\Sigma)$  and hence
	\item $\pi_1^{et}(Y)\isom \pi_1^{et}(X)$ is an isomorphism of topological groups, and
	\item hence $X,Y$ are anabelomorphic complex varieties.
\end{enumerate}
\eenum
\ethm
\bp 
The assertions {\bf(1)--(4)} are immediate from \cite{nakai1972} and \cite[Appendix: Higher Dimensions]{nakai-book}. The assertion {\bf(5)}(a,b) follows from the fact that the Riemannian manifolds $\Sigma'$ and $\Sigma$ are Nakai quasi-isometric and hence homeomorphic. The assertion {\bf(5)}(c) follows from a well-known theorem of Grothendieck \cite{sga1}. The assertion {\bf(5)}(d) is immediate from the definition of anabelomorphic varieties.
\ep

\brem 
Let me remark that in the context of Arakelov Theory (and hence in Diophantine geometry) one works with the existence of Green's functions. An important reason why one works with Nakai quasi-isometries is that Nakai quasi-isometries (i.e. quasi-conformal mappings if $\dim(X/\C)=1$) preserve the property of the possession of non-constant Green's functions i.e. if $Y\to X$ is a Nakai quasi-isometry then either both $X,Y$ possess a non-constant Green's function or neither one does. This is established in \cite{nakai1972}, \cite{nakai-book}.
\erem

\subsection{Variation of Berkovich Analytic Structure}
As observed in \Cref{th:nakai1}, \Cref{th:nakai2}, in the archimedean case, a variation of Banach spaces of functions underlies classical Teichm\"uller spaces. For $E$ non-archimedean, the category $\fjxe$ shares this important property with the classical Teichm\"uller Space. Namely, the Berkovich analytic space structure provided by $(Y/E',(E'\into K, K^\flat \isom F), *_K:\sM(K)\to \yan_{E'})$ varies as  $(Y/E',(E'\into K, K^\flat \isom F), *_K:\sM(K)\to \yan_{E'})$ varies over objects of $\fjxe$. This is
recorded in the following consequence of \Cref{thm:main2}. This is one of the properties of $\fjxe$ which justifies the term Arithmetic Teichm\"uller Space $\fjxe$.

\bthm\label{th:variation-of-arith-hol-strs}\ 
Let $X/E$ be a geometrically connected, smooth, projective variety over a $p$-adic field $E$.
\benumlab
\item Each $(Y/E',(E'\into K, K^\flat \isom F), *_K:\sM(K)\to \yan_{E'})\in \fjxe$ provides the analytic spaces $\yan_{E'}, \yan_{K}$ and the morphism  $\yan_K\to \yan_{E'}$ of analytic spaces provided by the given isometric inclusion $E'\into K$. 
\item Two objects $(Y_1/E_1',(E_1'\into K_1, K_1^\flat \isom F), *_{K_1}:\sM(K_1)\to \yan_{E_1'})$ and 
$(Y_2/E_2',(E_2'\into K_2, K_2^\flat \isom F), *_{K_2}:\sM(K_2)\to \yan_{E_2'})$ of $\fjxe$ need not be isomorphic even if $Y_1/E_1=Y_2/E_2$.
\item If $Y_1/E_1$ and $Y_2/E_2$ are both projective then
 $(\yan_{E_1'},\yan_{K_1},\yan_{K_1}\to \yan_{E_1'})$ need not be isomorphic to $(\yan_{E_2'},\yan_{K_2},\yan_{K_2}\to \yan_{E_2'})$ even if $Y/E_1'=Y/E_2'$.
\eenum
\ethm
\bp 
The first assertion is clear. The second assertion is immediate from the fact that the two given untilts, $(E_1'\into K_1, K_1^\flat \isom F)$, $(E_2'\into K_2, K_2^\flat \isom F)$ of $F$ need not be isomorphic. The third assertion is immediate from \Cref{thm:main2}.
\ep

\brem 
One way of understanding the above result is that the Berkovich analytic space structure i.e. the analytic space datum provided by  objects of $\fjxe$ varies over $\fjxe$.
\erem

\subsection{Existence of distortions of $p$-adic metrics}
A key property classical Teichm\"uller spaces, grounded in the basic theory of quasi-conformal mappings  of Riemann surfaces is that metrics appear scaled relative to each other (these scalings are encoded in the Beltrami parameter) and local geometries are distorted relative to each other \cite[Chapter 1 and Chapter 4]{imayoshi-book}. 

In the context of Berkovich analytic spaces \cite{berkovich-book}, there is a natural notion of metrics namely, the local Banach norms provided in the datum of a Berkovich analytic space. The following theorem demonstrates that for this notion of Banach norms, objects of $\fjxe$ share this property of local metric scalings with the objects of classical Teichm\"uller Space (\cite[Theorem 1.7]{imayoshi-book}).

\bthm\label{th:distortions} 
Let $X/E$ be geometrically connected, smooth, quasi-projective variety over a $p$-adic field $E$. 
Let $(Y/E',(E'\into K_1, K_1^\flat \isom F), *_{K_1}:\sM(K_1)\to \yan_{E'})$ and $(Y/E',(E'\into K_2, K_2^\flat \isom F), *_{K_2}:\sM(K_2)\to \yan_{E'})$ be two objects of $\fjxe$.
Let $(E',\abs{-}_{K_1})$ $(E',\abs{-}_{K_2})$ be the restriction of the absolute value of $K_1$ to $E'$ (resp. of $K_2$ to $E'$). Then there exists a real number $\alpha\in\R$ such that
$$(E',\abs{-}_{K_2}) =(E',\abs{-}^\alpha_{K_1}).$$  
However, there is no uniform choice of $\alpha$ for all pairs of objects of $\fjxe$.
\ethm
\bp 
One has two embeddings $\iota_1:E'\into K_1$ and $\iota_2:E'\into K_2$. Write $E_1'=\iota_1(E')\subset K_1$ and similarly define $E_2'$ i.e. view $E'$ as a (discretely) valued subfield of $K_j$ via $\iota_j$.  Abstractly, $E'\isom E_1'\isom E_2'\isom E'$ and any  abstract  isomorphism  $E_1'\isom E_2'$ is in fact an isomorphism of discretely valued fields, and the two valuations on $E'$ induced from $K_1,K_2$ on $E'$  are equivalent. Hence existence of $\alpha\in\R$ as asserted is clear.  
The second assertion is now immediate from the existence of Fargues-Fontaine curves \cite{fargues-fontaine}. 
\ep

\bcor\label{cor:distortions} 
In particular, in the notation of \Cref{th:distortions} and its proof, the $p$-adic metrics (i.e. the local Banach norms) of $\yan_{E_2'}$ are scaled relative to that of $\yan_{E_1'}$. However, there is no uniform choice of $\alpha$ for all pairs of objects of $\fjxe$.
\ecor
\bp 
Any isomorphism of valued fields $E_1'\isom E_2'$ (as in the above proof) also gives isomorphism of analytic spaces $\yan_{E'_1}\isom \yan_{E'_2}$. The theorem asserts that valuation of $E_1'$ is scaled relative to $E_2'$. 
The remaining assertion is immediate from the above theorem and the existence of Fargues-Fontaine curve $\sY_{F,E'}$ established in \cite[Chapitre 2, 6]{fargues-fontaine}.
\ep

\brem 
In \iutf{IUT1, Remark 3.9.3}, Mochizuki asserts the existence of local distortions \Cref{cor:distortions} as a key point in his theory.
\erem

\brem 
In the archimedean case, the existence of distortions of metrics is built into the theory of Nakai quasi-isometries in all dimensions by \cite{nakai1972}.
\erem

\subsection{Invariance of the tempered fundamental groups in Arithmetic Teichm\"uller Spaces} Here is an important but trivial consequence of the definition of $\fjxe$: 
\bpro\label{pr:fund-group-same} Let $X/E$ be a geometrically connected, smooth, quasi-projective variety over a field $E$.
\benumlab
\item If $E$ is a $p$-adic field, then for every triple $(Y/E',E'\into K)\in\fjxe$ one has  an anabelomorphism  $$\pit{Y/E'}\isom \pit{X/E},$$
\item inducing an isomorphism of their geometric tempered fundamental subgroups
$$	\pit{Y/K}\isom\pit{X/\C_p},$$ 
and hence an isomorphism 
$$G_{E';K}\isom G_{E;\C_p}$$ of absolute Galois groups (i.e. $E'$ and $E$ are anabelomorphic $p$-adic fields).
\item In the archimedean case (i.e. $E=\C$), for every object $(Y,f:Y\to X)$, on has an isomorphism of topological fundamental groups
$$\pi_1^{top}(Y)\isom \pi_1^{top}(X)$$
and hence of their profinite completions
$$\pi_1^{et}(Y)\isom \pi_1^{et}(X).$$
\eenum
\epro
\bp 
The proof is clear from the definition of $\fjxe$ in both the cases.
\ep
\brem 
This together with \Cref{th:variation-of-arith-hol-strs}, \Cref{th:distortions}, \Cref{cor:distortions} shows why $\fjxe$ has properties similar to classical Teichm\"uller Spaces.
\erem

\subsection{The Subcategories $\fjxe_F$ of $\fjxe$} Since $F$ is an arbitrary, algebraically closed perfectoid field of characteristic $p>0$, the category $\fjxe$ is too big in practice. Hence one works with the following full subcategories of $\fjxe$. 
\bpro
There is a natural functor 
$$(Y/E',E'\into K,\iota:K^\flat\isom F)\mapsto F$$
from $\fjxe$ to the category of perfectoid fields of characteristic $p>0$, with isometries of such fields as morphisms.
\epro
\bp 
The proof is clear.
\ep

\begin{defn}\label{def:subcat-F}
Let $E$ be $p$-adic field and fix $F$ to be an algebraically closed perfectoid field of characteristic $p>0$. Then let $$\fJ(X,E)_F$$ be the full subcategory of $\fjxe$ whose objects are
$$\left\lbrace (Y/E',E'\into K,\iota:K^\flat\isom F') \in \fjxe: F'\text{ isometric to }F\right\rbrace.$$
Thus $\fjxe_F$ is the fiber, over $F$, of the above functor.
\end{defn}

\subsection{Arithmetic Teichm\"uller Spaces With Auxillary Conditions} For greater flexibility and with a view to applications, it is useful to consider variants of the construction of $\fJ(X,E)$. Let $\Sigma$ be a (finite) set of geometric or arithmetic conditions one can impose on the data $(Y/E',E\into K)$ (all objects of $\fjxe$ need not satisfy the conditions $\Sigma$). 

Let 
$$\fJ^\Sigma(X,E)$$ be the full subcategory of $\fjxe$ consisting of objects of $\fjxe$ which satisfy $\Sigma$:
\be \fJ^\Sigma(X,E)=\left\{ (Y/E',E'\into K)\in \fJ(X,E) \text{ and } (Y/E',E'\into K) \text{ satisfies } \Sigma \right\}.
\ee

\newcommand{\fJhypxe}{\fJ^{hyp}(X,E)}
\para Two important examples of $\Sigma$ which are useful from the point of view of \iut, are the following:
\benumlab
\item Let
$$\Sigma_{hyp}=\{ Y/E' \text{ is  hyperbolic}\}.$$
\item 
If $\dim(X)=1$, then let
$$\Sigma_{Belyi}=\{ Y \text{ hyperbolic and also of strict Belyi Type}\}.$$
Note that  $\dim(X)=1$ implies (by \Cref{def:arith-hol-space-local}) that $\dim(Y)=1$; for the definition of Strict Belyi Type see \cite[Definition 3.5]{mochizuki-topics2}.
\eenum
For $\Sigma_{hyp}$, I will write $\fJhypxe$ instead of $\fJ^{\Sigma_{hyp}}(X,E)$. Hence
$$\fJ^{hyp}(X,E)=\left\{ (Y/E',E'\into K)\in \fJ(X,E):  Y/E' \text{ is  hyperbolic variety} \right\},$$
and for $\Sigma_{Belyi}$, I will write
$$\fJ^{SB}(X,E)=\left\{ (Y/E',E'\into K)\in \fJ(X,E): Y/E' \text{ is  of Strict Belyi Type} \right\}.$$

\bpro\label{pr:top-types}
Let $X/E$ be a geometrically connected, smooth, hyperbolic curve over a $p$-adic field $E$ and let
$(Y/E',(E'\into K, K^\flat \isom F), *_K:\sM(K)\to \yan_{E'})\in\fJ^{hyp}(X/E)$ be any object. Then 
\benumlab
\item  $Y/E'$ has the same topological type=(genus, number of punctures) as that of $X/E$.
\item any anabelomorphism
$$\pit{Y/E'}\mapright{\isom} \pit{X/E}$$
provides an anabelomorphism
$$G_{E'} \isom G_E$$
i.e. the $p$-adic fields $E'$ and $E$ are anabelomorphic.
\eenum
\epro
\bp 
The first assertion is \cite[Lemma~1.3.9]{mochizuki04}. The second assertion is immediate from \cite[Lemma~1.3.8]{mochizuki04}. 
\ep

\bpro\label{pro:sb-type-hyp} 
Let $X/E$ be a geometrically connected, smooth, hyperbolic curve over a $p$-adic field $E$.
Suppose $X/E$ is also of Strict Belyi Type. Then $$\fJ^{SB}(X,E)=\fJhypxe.$$
\epro
\bp 
To prove this it is enough to show that if $X/E$ is of Strict Belyi Type, then every hyperbolic curve $Y/E'$ occurring as a part of the datum of any object of $\fJhypxe$ is also of Strict Belyi Type. This is immediate from \Cref{pr:fund-group-same}, the proof of \cite[Theorem 2.10]{mochizuki07-cuspidalizations-proper-hyperbolic}, and the definition of Strict Belyi Type hyperbolic curves \cite[Definition 2.9]{mochizuki07-cuspidalizations-proper-hyperbolic}.
\ep

\brem 
Thus for a geometrically connected, smooth, hyperbolic curve, $\fJ^{hyp}(X/E)$ has properties similar to classical Teichm\"uller Spaces: namely local scalings (\Cref{th:distortions}, \Cref{cor:distortions}), the topological types of all objects and  the tempered fundamental group (and hence \'etale fundamental group) is fixed by \Cref{pr:top-types}.
\erem

\para For hyperbolic curves of Strict Belyi Type one has the following:
\bpro\label{pr:belyi}
Let $E$ be a $p$-adic field and $X/E$ be a geometrically connected, smooth, quasi-projective, hyperbolic curve  of Strict Belyi Type over $E$. Then
\benumlab
\item For every $(Y/E',E'\into K)\in \fJ^{SB}(X,E)$ one has an isomorphism of schemes (over $\Z$) $$Y\isom X.$$
\item Hence one has a natural action of ${\rm Aut}(\Pi)$ on $\fJ^{SB}(X,E)$ via its action on
$${\rm Isom}(Y,X)\isom {\rm Isom}^{\rm out}(\pi_1(Y),\pi_1(X)).$$ 
\item However  two such triples $(Y/E',E'\into K),(Y''/E'',E''\into K')\in\fJ^{SB}(X,E)$ may  not be isomorphic in general.
\eenum
\epro
\bp 
The first two assertions are immediate from Mochizuki's proof of the Absolute Grothendieck Conjecture  for hyperbolic curves over $p$-adic fields of Strict Belyi Type (see \cite[Corollary 2.12]{mochizuki07-cuspidalizations-proper-hyperbolic}) and the construction of $\fJ^{SB}(X,E)$. The second assertion is immediate from \Cref{le:non-isom-arith-hol-strs}.
\ep

\para Let me also record the following proposition in the special case of hyperbolic curves of strict Belyi Type.
\bpro\label{cor:relation-with-iut} 
Let $E$ be a $p$-adic field. Let $X/E$ be a geometrically connected, smooth, quasi-projective, hyperbolic curve over $E$ of strict Belyi type. Then there is a natural action of $Out(\Pi)$ where $\Pi=\pit{X/E}$ on $\fJ_{SB}(X,E)$.
\epro
\bp 
This is immediate from Proposition~\ref{pr:belyi}: for any $(Y/E',E'\into K)$ one has an isomorphism of schemes $Y\isom X$. 
\ep

\subsection{Existence of Virasoro Type Uniformization}\label{ss:virasoro-symmetry}
\newcommand{\syfep}{\mathscr{Y}_{F,E'}}

The data $$(Y/E',(E'\into K, K^\flat \isom F), *_K:\sM(K)\to \yan_E)\in\fjxe$$ is quite rigid but the data $$(Y/E',(E'\into K, K^\flat \isom F))$$ has symmetries stemming from symmetries of both $Y/E'$ and $(E'\into K, K^\flat \isom F)$. Hence along with the data $(Y/E',(E'\into K, K^\flat \isom F), *_K:\sM(K)\to \yan_E)\in\fjxe$ it is also convenient to work with the category of data $(Y/E',(E'\into K, K^\flat \isom F))$  which omits  the geometric base-point $*_K:\sM(K)\to \yan_E$.  This leads to the following definition:
\begin{defn}\label{def:arith-hol-space-local-detached}
	Let $E$ be a $p$-adic field and let $X/E$ be a geometrically connected, smooth, quasi-projective variety over $E$. Let  $\ufjxe$  be the  category whose objects are	pairs
	$$	(Y/E',(E'\into K, K^\flat \isom F))$$
	where 
	\benumlab
	\item $E'$ is a $p$-adic field, and 
	\item $F$ is some perfectoid field of characteristic $p>0$.
	\item $(E'\into K, K^\flat \isom F)$ is an $E$-untilt of $F$.	
	\item $Y/E'$ is anabelomorphic to $X/E$,
	\item $\dim(Y)=\dim(X)$; and
	\eenum
	Morphisms in $\ufjxe$ are morphisms of the pairs $(Y/E',(E'\into K, K^\flat \isom F))$.
\end{defn}

\bpro\ 
\benumlab
\item There is a forgetful functor 
$$\fjxe\to \ufjxe$$ 
given by forgetting the geometric base-point:
$$(Y/E',(E'\into K, K^\flat \isom F), *_K:\sM(K)\to \yan_E)\mapsto (Y/E,(E\into K, K^\flat \isom F)).$$
\item Given any $$(Y/E',(E'\into K, K^\flat \isom F))\in\ufjxe,$$ each choice of a geometric base-point $$*_K:\sM(K)\to \yan_{E'}$$ provides an object $$(Y/E',(E'\into K, K^\flat \isom F), *_K:\sM(K)\to \yan_{E'})\in\fjxe.$$
\item There is a functor 
$$\ufjxe\ni(Y/E',(E'\into K, K^\flat \isom F))\mapsto (\yan_{K}\to \yan_{E'})$$
from $\ufjxe$ to the category of morphisms of analytic spaces $\Q_p$-analytic spaces (with morphisms defined in the obvious way).
\item Given $(Y/E',(E'\into K, K^\flat \isom F))\in \ufjxe$, any choice of $K$-geometric base-point $*_K:\sM(K)\to \yan_{E'}$ gives tempered fundamental group
$$\pit{Y/E'} \quad (\isom \pit{X/E}).$$
\eenum
\epro
\bp 
All the assertions are clear.
\ep

\begin{defn}\ 
\benumlab
\item I will refer to $(\yan_{K}\to \yan_{E'})$ as the arithmetic holomorphic structure associated to $$(Y/E',(E'\into K, K^\flat \isom F))\in\ufjxe.$$ 
\item Let $\ufjxe_{E}$ be the full subcategory of $\ufjxe$ consisting of the objects of the form $(Y/E,(E\into K, K^\flat \isom F))\in\ufjxe$. Evidently $(X/E, (E\into K, K^\flat \isom F))\in\ufjxe_E$ and hence $\ufjxe_E$ is not an empty category.
\eenum
\end{defn}

\para The goal of the next few paragraphs is to demonstrate that there are natural symmetries which transform the arithmetic holomorphic structures associated to $(Y/E,(E\into K,K^\flat\isom F))\in\ufjxe_{E}$. These symmetries are the arithmetic analogs fo the Virasoro action on certain extended moduli spaces of curves and Teichm\"uller spaces established in \cite{kontsevich87}, \cite{beilinson88} (for an accessible account see \cite[Virasoro Uniformization Theorem]{frenkel01-book}). The existence of this property in the present context reinforces my point that one is dealing with a Teichm\"uller Theory of the arithmetic sort in which arithmetic holomorphic structures change under certain Virosoro action type symmetries. 

The key point in the demonstration is that an $E$-untilt  $(E\into K,K^\flat\isom F)$ of $F$ can be viewed as a closed classical point on the Fargues-Fontaine curve $\syfe$  constructed \cite{fargues-fontaine}. A less precise version of this symmetry occurs, in the form of ``Indeterminacy Ind2'' in \iut, where it plays a central role in \iutf{IUT3, Theorem 3.11 and Corollary 3.12}.

\para Let $E$ be a $p$-adic field and let $\pi$ be a uniformizer for its ring of integers, let $q$ be the cardinality of the residue field $\O_{E}/\pi$ of $E$. By a \emph{Lubin-Tate formal group $\sG$ over $\O_E$}, I mean a formal group constructed by \cite{lubin65}, using some polynomial $Q(T)\in\O_E[T]$ satisfying the following two hypothesis of \cite{lubin65}:
\benumlab
\item $Q(T)=\pi T+O(T^2)$, and
\item $Q(T)=T^q\bmod{\pi}$.
\eenum
By \cite{lubin65} the formal groups determined by two such polynomials is naturally isomorphic. By the results of \cite[Chapter 1, 2]{fargues-fontaine} especially \cite[Proposition 2.1.7]{fargues-fontaine}, the $\O_E$-algebra $W_{\O_E}(\O_F)$ required in the construction of the Fargues-Fontaine curves is independent of the choice of the Lubin-Tate polynomial $Q(T)$ used to define $\sG$ and in particular these constructions are independent of the choice of the Lubin-Tate group $\sG$.

\para Associated to a Lubin-Tate $\O_E$-formal group $\sG$ over $\O_E$, is a $\pi$-divisible group over $\O_E$,  and its special fibre (over the residue field of $\O_E$) and also a $\pi$-divisible formal $\O_E$-module. I will pass between these objects whenever needed (to invoke results of \cite{fargues-fontaine}), but \emph{beware that I will notationally conflate all of these objects as $\sG$.} Hopefully there will be no confusion.

\para I will also use $\sG$ for the special fiber of $\sG/\O_E$, I hope that readers will be able to unravel the usage from the context (in \cite{fargues-fontaine}, the special fiber is denoted by $\sG_k$ where $k$ is the residue field of $F$). This means for example where \citeauthor{fargues-fontaine} write $\sG_k(\O_F)$, I will simply write $\sG(\O_F)$. By \cite[Proposition 4.4.1]{fargues-fontaine} $\sG(\O_F)$ is naturally a Banach space over $E$.

\para Suppose $E$ is a $p$-adic field and $\O_E$ its ring of $p$-adic integers, let $K$ be an algebraically closed perfectoid field and let $\sG$ be an Lubin-Tate group or a $\O_E$-formal group equipped with an homomorphism $\O_E\to\End_{\O_K}(\sG)$. Let $\pi\in\O_E$ be a uniformizer. Let 
$$\widetilde{\sG(\O_K)}=\invlim_{\text{mult. by }\pi}\sG(\O_K),$$
this is naturally an $E$-vector space (and hence also an $\O_E$-module).  Note that $\sgtok$   is denoted as $X(\sG)(\O_K)$ \cite[Chap IV]{fargues-fontaine}.

\para Let me explicate this for $E=\Q_p,F=\cpt$ and $K=\C_p$, here $\sG(\cpt)\isom\sgtocp$. In this case $\sG(\O_{\C_p})=\mathfrak{m}_{\O_{\C_p}}$ is considered as a $\Z_p$-module via the Lubin-Tate action in which the endomorphism $p$ acts on $\sG(\O_{\C_p})$ by  the endomorphism $z\mapsto z^p+pz$ of $\sG(\O_{\C_p})$ and $\sgtocp$ is the $\Q_p$-vector space  obtained by formally inverting the Lubin-Tate action of $p$ on the group  $\sG(\O_{\C_p})\isom \mathfrak{m}_{\O_{\C_p}}$.

\para Let $E$ be a $p$-adic field and let $F$ be an algebraically closed perfectoid field of characteristic $p>0$. Recall that in \cite[Chapter 2]{fargues-fontaine} the construction of the Fargues-Fontaine curve proceeds via the construction of an auxiliary curve, denoted  in loc. cit by $Y_{F,E}$ (resp. $X_{F,E}$) (or simply by $Y$ if the choice of $F,E$ is unambiguous) and denoted here by $\syfe$ (resp. $\sxfe$). More precisely $\syfe$ is constructed as an adic space  (but I will not use this fact here).

\para Of particular interest to us are the  sets of closed points of $\abs{\syfe}$ (resp. $\abs{\sxfe}$) of closed (classical points) of $\syfe$ (resp. closed points of $\sxfe$). By \cite[Corollaire 2.2.2 and D\'efinition 2.3.1]{fargues-fontaine}, the set of closed classical points of $\syfe$ is in bijection $\abs{\syfe}\ni y\mapsto (E\into K_y,K_y^\flat\isom F)$ with the set of untilts of $F$ (up to isomorphisms of untilts) and where $K_y$ is naturally identified with the residue field of $y\in\syfe$. This will be used in what follows. The curve $\syfe$ is equipped with a natural Frobenius morphism $\varphi:\syfe \to \syfe$ and $\sxfe=\syfe/\varphi^\Z$  is the (adic space) quotient by the infinite cyclic group generated by the Frobenius morphism $\varphi$. In the discussion which follows, I will habitually conflate $\sxfe$ and $\syfe$ with $\abs{\sxfe}$ (resp. $\abs{\syfe}$). 

\para By \cite[Th\'eor\`eme 6.5.2(4)]{fargues-fontaine}, one has a canonical identification 
$$\abs{\syfe}/\varphi^\Z\isom \abs{\sxfe},$$
given by $y\mapsto \{\varphi^n(y):n\in\Z\}$. It is standard that $\syfe/\varphi^\Z\to \sxfe$ is in fact a morphism of adic spaces which provides the above identification on points.

The following theorem demonstrates how anabelomorphy of $p$-adic fields enters the theory of Fargues-Fontaine curves $\syfe$ and $\sxfe$.
\bthm\label{th:anab-ff-curves} 
Let $E,E'$ be anabelomorphic $p$-adic fields and let $F$ be algebraically closed perfectoid field of characteristic $p>0$. Then any anabelomorphism $\alpha:G_E\to G_{E'}$ induces a homeomorphism of topological spaces compatible  with the respective topological group actions
$$G_E\act \abs{\syfe} \isom \abs{\syfep}\curvearrowleft G_{E'}$$ and 
$$G_E\act\abs{\sxfe} \isom \abs{\sxfep}\curvearrowleft G_{E'}.$$
\ethm
\bp 
One has the natural identifications (\cite[Proposition 2.1.10]{fargues-fontaine})
\be\label{eq:ff-curve-y}\abs{\syfe}=\left(\sG(\O_F)-\{0\}\right)/\O_E^*,\ee
and (\cite[Theorem 6.5.2]{fargues-fontaine})
\be\label{eq:ff-curve-x}\abs{\sxfe}=\left(\sG(\O_F)-\{0\}\right)/E^*,\ee
and similarly for $E'$. The topological space $\left(\sG(\O_F)-\{0\}\right)$ is independent of $E,E'$ and the anabelomorphism $\alpha:G_{E}\to G_{E'}$ induces  isomorphisms of the topological groups
$$\O_E^*\isom \O_{E'}^*$$
and 
$$E^*\isom {E'}^*.$$
Hence all the assertions are immediate.
\ep
\brem
On the other hand, I have shown in \cite{joshi-gconj} that, if $E,E'$ are strictly anabelomorphic $p$-adic fields, then schemes $\sxfe$ and $\sxfep$ are anabelomorphic but not  isomorphic as schemes.
\erem

\para With this preparation let me prove the following:
\bthm\label{th:main3} 
Let the notation and assumptions be as in the previous paragraph. Assume $K$ is an algebraically closed perfectoid field with $K^\flat=F$. 
\benumlab
\item\label{th:main3-1} The isomorphism class of the topological $\O_E$-module $\sgtok$ is independent of $K$ more precisely, there is a natural homeomorphism of $\O_E$-modules
$$\sgtok\isom \sG(\O_F),$$
(in fact this is an isomorphism of Banach $E$-vector spaces).
\item\label{th:main3-2} The isomorphism class of the topological $\O_E$-module $\sG(\O_F)$ is independent of the choice of the Lubin-Tate $\O_E$-formal group $\sG$.
\item\label{th:main3-3} There is a natural action of the group  ${\rm Aut_{\O_E}}(\sG(\O_F))$, of topological automorphisms of the $\O_E$-module $\sG(\O_F)\isom\sgtok$, on the set of closed points of degree one of the Fargues-Fontaine curve $\syfe$, arising from the natural identification  $$\abs{\syfe}=\left(\sG(\O_F)-\{0\}\right)/\O_E^*.$$
\item\label{th:main3-4} Let $y\in\syfe$ be a closed point of degree one.  Then one has an action of ${\rm Aut_{\O_E}}(\sG(\O_F))$ on closed points of degree one of  $\sxfe$ via mapping $$\{\varphi^n(y):n\in\Z\}\mapsto \{\varphi^n(\sigma(y)):n\in\Z\}.$$
\item\label{th:main3-5}
Thus given any topological $\O_E$-linear automorphism  $$\sigma:\sG(\O_F)\mapright{\isom}\sG(\O_F),$$ and a closed point $y\in\syfe$ of degree one, giving an untilt $y=(E\into K_y,K_y^\flat\isom F)$ of $F$ with residue field  $K_y$, there is an untilt $\sigma(y)=(E\into K_{\sigma(y)},K_{\sigma(y)}^\flat\isom F)$ given via \eqref{eq:ff-curve-y} with residue (also perfectoid, algebraically closed field) $K_{\sigma(y)}$ and with isometries $K_y^\flat\isom F\isom K_{\sigma(y)}^\flat$ and an embedding $E\into K_{\sigma(y)}$.
\item\label{th:main3-6} In particular ${\rm Aut}_{\O_E}(\sG(\O_F))$ acts naturally on $\ufjxe_F$ via $$(Y,y=(E\into K,K^\flat\isom F))\mapsto (Y, \sigma(y)=(E\into K_{\sigma(y)},K_{\sigma(y)}^\flat\isom F))$$ for all  $\sigma\in {\rm Aut}_{\O_E}(\sG(\O_F))$.
\eenum
\ethm
\bp 
Before proceeding to the proofs let me remark that items \ref{th:main3-1} and \ref{th:main3-2} are due to \cite{fargues-fontaine} and I include them here for completeness.  The assertion \ref{th:main3-1} is \cite[Proposition~4.5.11]{fargues-fontaine} (what I have denoted  as $\sgtok$ is denoted by $X(\sG)(\O_K)$ in loc. cit.). The independence from the choice of the $\O_E$-formal Lubin-Tate group $\sG$ is clear from \cite{lubin65} as the Lubin-Tate $\O_E$-formal group over $\O_E$ is unique up to isomorphism by \cite{lubin65}.

The identification of $\abs{\syfe}$ with the  $\left(\sG(\O_F)-\{0\}\right)/\O_E^*$ is \cite[Proposition 2.1.10]{fargues-fontaine} and hence for any $\sigma\in {\rm Aut_{\O_E}}(\sG(\O_F))$, $\sigma$ is evidently a bijection on $\left(\sG(\O_F)-\{0\}\right)/\O_E^*$. Thus the  claim \ref{th:main3-3} is immediate. 

The proof of \ref{th:main3-4} is now clear now that \ref{th:main3-3} has been established. 

The assertion \ref{th:main3-5} is self-evident and it implies \ref{th:main3-6}. This completes the proof.
\ep

\para\label{pa:virasoro-unif-remark} For readers familiar with the Geometric Langlands Program \cite{beilinson00b} over $\C$, let me remark that the action of ${\rm Aut}_{\Z_p}(\sG(\O_{\cpt}))$ considered in \Cref{th:main3} \emph{is the $p$-adic analog of the action of the Virasoro Algebra on moduli spaces of marked Riemann surfaces described in the Virasoro uniformization Theorem \cite[Section 4]{beilinson88}, \cite{beilinson00b}, \cite{frenkel01-book}}. In the Geometric Langlands setting of \cite{beilinson00b},  the Virasoro algebra plays a fundamental role and manifests itself via the action of the group scheme $\C\subset R\mapsto {\rm Aut}_{cont}(R(\kern-0.5mm(T)\kern-0.5mm))$. As is described in \cite[Section 4
]{beilinson88} or \cite[Theorem~17.3.2]{frenkel01-book}, this action  also changes  complex structures of marked Riemann surfaces (in general). \emph{As has been noted above, ${\rm Aut}_{\Z_p}(\sG(\O_{\cpt}))$ acts by changing the analytic structure   of $(X\times_EK)^{an}$ and hence must be considered as the $p$-adic analog of the Virasoro action in the complex setting.} As was remarked in \cite{beilinson88}, the Virasoro uniformization Theorem complements the Teichm\"uller Uniformization.

\para An important consequence of this is that topological $\O_E$-linear automorphisms of $\sG(\O_F)$ can be used to change arithmetic holomorphic structures in the sense of Theorem~\ref{thm:main2}. In \iut, Mochizuki claims this assertion in the guise of ``Indeterminacy Ind2'' \iutf{IUT3, Theorem 3.11}.

\newcommand{\syFqp}{\sY_{F,\Q_p}}
\para An important consequence of Theorem~\ref{th:main3} and \ref{le:artin-hasse-exp} is the following:
\bthm\label{co:action-on-cpt-cat}
Let $F$ be an algebraically closed, perfectoid field of characteristic $p>0$. There is a natural action of ${\rm Aut}_{\Z_p}(\sG(\O_{F}))\isom {\rm Aut}_{\Z_p}(\hgm(\O_{F}))$ on the closed points $\abs{\sY_{F,\Q_p}}$ of $\sY_{F,\Q_p}$, which provides an action of ${\rm Aut}_{\Z_p}(\hgm(\O_{F}))$ by correspondences on $\fJ(X,E)_{F}$. Explicitly this is given as follows. 
\benumlab
\item Let 
$$\sG(\O_{F})) \mapright{\sigma} \sG(\O_{F}))$$
equivalently, let $$  $$ 
$$\hgm(\O_{F}) \mapright{\sigma} \hgm(\O_{F})$$
be a $\Z_p$-linear isomorphism of topological groups.
\item 
 Let $\tau: \abs{\syfep}\to\abs{\syFqp}$ be  the finite morphism  given by mapping an $E'$-untilt $y\in\abs{\syfep}$ of $F$ to the underlying $\Q_p$-untilt $\tau(y)=x\in\abs{\syFqp}$ of $F$ \cite[Proposition 2.3.20]{fargues-fontaine}.
\item  Let $y_1,\ldots,y_m\in\abs{\syfep}$ be all the points lying over $\sigma(x)\in\syFqp$ i.e. $\tau^{-1}(\sigma(x))=\{ y_1,\ldots,y_m\}\subset\abs{\syfep}$.
\eenum
Then the divisorial correspondence (i.e. one-to-many mapping between object class)
\be  
\ufjxe_F \dashrightarrow \ufjxe_F
\ee
is given by the rule
$$(Y/E',y=(E'\into K_y,K_y^\flat\isom F))\mapsto ((Y/E',y_j=(E'\into K_{y_j},K_{y_j}^\flat\isom F)))_{j=1,\ldots,m}.$$ 
Moreover, the number $m=[E':\Q_p]=[E:\Q_p]$ depends only on $E$ and not on the chosen object of $\ufjxe$.
\ethm

\bp 
Let assertion $m=[E':\Q_p]$ is \cite[Proposition 2.3.20]{fargues-fontaine}. The assertion $[E':\Q_p]=[E:\Q_p]$ is immediate from the fact $E'$ and $E$ are anabelomorphic $p$-adic fields \Cref{def:arith-hol-space-local}{\bf(4)}, \Cref{pr:top-types}. The rest of the proof is clear from the method of proof of \Cref{th:main3}.
\ep

\para The results of the preceding paragraphs can be assembled into the following theorem: 
\bthm\label{th:main4.5-const-local-teich}
Let $E$ be a $p$-adic field, let $X/E$ be a geometrically connected, smooth, quasi-projective variety over $E$. Then there exists a category $\fjxe$, called the \emph{Arithmetic Teichm\"uller Space associated to $X/E$}, with the following properties:
\benumlab
\item  objects of $\fjxe$ are $(Y/E',(E'\into K, K^\flat \isom F), *_K:\sM(K)\to \yan_{E'})$ (\Cref{def:arith-hol-space-local}).  Morphisms between these objects are defined in the obvious way. 
\item The category $\fjxe$ is an  \emph{anabelian variation providing $\Pi=\pit{X/E}$} (see \ssep\ref{pa:geometric-anabelian-var}, \Cref{th:anab-var}) i.e.
for any  $(Y/E',(E'\into K, K^\flat \isom F), *_K:\sM(K)\to \yan_{E'})\in\fjxe$, one has an isomorphism of topological groups $$\pit{Y/E'}\isom  \pit{X/E},$$
inducing an isomorphism of their geometric tempered fundamental subgroups
$$\pit{Y/K}\isom  \pit{X/K},$$
\item There are forgetful functors (see \Cref{pr:various-functors}):
\benum[label={{\bf(\roman{*})}}]
\item $(Y/E',(E'\into K, K^\flat \isom F), *_K:\sM(K)\to \yan_{E'}) \mapsto Y/\Z$ (i.e. to Schemes/$\Z$).
\item $(Y/E',(E'\into K, K^\flat \isom F), *_K:\sM(K)\to \yan_{E'}) \mapsto E'$ (i.e. to $p$-adic fields).
\item $(Y/E',(E'\into K, K^\flat \isom F), *_K:\sM(K)\to \yan_{E'}) \mapsto K$ (i.e. to algebraically closed perfectoid fields of characteristic zero).
\item $(Y/E',(E'\into K, K^\flat \isom F), *_K:\sM(K)\to \yan_{E'}) \mapsto K^\flat$  (i.e. to algebraically closed perfectoid fields of characteristic  $p>0$).
\eenum
\item There are functors to analytic spaces (see \Cref{pr:various-functors}) $$(Y/E',(E'\into K, K^\flat \isom F), *_K:\sM(K)\to \yan_{E'})\mapsto (\yan_{E'}),$$
and
$$(Y/E',(E'\into K, K^\flat \isom F), *_K:\sM(K)\to \yan_{E'})\mapsto (\yan_K),$$
and
$$(Y/E',(E'\into K, K^\flat \isom F), *_K:\sM(K)\to \yan_{E'})\mapsto (\yan_K\to \yan_{E'}).$$
\item There are functors to Mochizuki's anabelian landscape (see \ssep\ref{pa:anab-functors}): (one uses the given perfectoid field to compute algebraic closures) $$(Y/E',(E'\into K, K^\flat \isom F), *_K:\sM(K)\to \yan_{E'})\mapsto \pit{Y/E'} \act\O_{\bE}^*\subset \O_K,$$
and also
$$(Y/E',(E'\into K, K^\flat \isom F), *_K:\sM(K)\to \yan_{E'})\mapsto \pit{Y/E'} \act\otmu\subset \O_{K}^*/\mu(K),$$ and similarly $$(Y/E',(E'\into K, K^\flat \isom F), *_K:\sM(K)\to \yan_{E'})\mapsto \pit{Y/E'} \act\O_{ \bE}^\triangleright\subset \O_{K}^\triangleright.$$
\item If $\dim(X)=1$ and $X$ is of Strict Belyi Type (this condition is defined in \cite[Definition 3.5]{mochizuki-topics2}) then one has an action of ${\rm Aut}(\Pi)$ on $\fJ(X,E)$ (Proposition~\ref{pr:belyi}).
\item For a fixed algebraically closed, perfectoid field $F$ of characteristic $p>0$, there are categories $\fJ(X,E)_F$ consisting of $(Y/E',(E'\into K, K^\flat \isom F), *_K:\sM(K)\to \yan_{E'})$ such that $K^\flat=F$.
\item Now fix an algebraically closed perfectoid field $F$ of characteristic $p>0$, a uniformizer $\pi$ for $E$ and let $\sG/\O_E$ be the Lubin-Tate formal group.  Then there is a natural action of ${\rm Aut}_{\O_E}(\sG(\O_{F}))$ on $\fJ(X,E)_F$ (Theorem~\ref{th:main3}). Notably for $F=\cpt$ one has a natural action (Corollary~\ref{co:action-on-cpt-cat}) of $${\rm Aut}_{\Z_p}(\sG(\O_{\cpt})) \isom {\rm Aut}_{\Z_p}(\sgtocp)$$ on $\fJ(X,E)_{\cpt}$.
\item The category  $\fJ(X,E)_{\cpt}$ is self-similar (Theorem~\ref{th:main6-fundamental-domain}).
\eenum
\ethm

\bp 
The only assertion which remains to be proved is the last claim that $\fJ(X,E)_{\cpt}$ is self-similar and this is Theorem~\ref{th:main6-fundamental-domain} and will be proved in the next section.
\ep

\section{Connectedness of Arithmetic Teichm\"uller Spaces and the abs. Grothendieck Conjecture}\label{se:connectedness}
\newcommand{\lmod}{L_{\rm mod}}
\nws
\renewcommand{\tripod}{\P^1-\{0,1,\infty\}}
In this section I want to demonstrate that there is a natural notion of connectedness for $\fjxe$ which is in fact related to the validity of the Absolute Grothendieck Conjecture for hyperbolic curves over $p$-adic fields. In particular, as Mochizuki has established this validity for curves of Strict Belyi Type, one obtains connectedness of $\fJhypxe$ for $X/E$ of strict Belyi Type.

\subsection{The Absolute Grothendieck Conjecture}\label{ss:abs-groth-conj}
\nwss
My formulation of the Absolute Grothendieck Conjecture over $p$-adic fields is based on \cite[Corollary 2.3]{mochizuki07-cuspidalizations-proper-hyperbolic}. This conjecture asserts the following. Assume that $E$ is a $p$-adic field, let $X/E$, (resp. $Y/E'$) be a geometrically connected, smooth quasi-projective variety over a $p$-adic field $E$ (resp. $E'$).  Then any isomorphism of topological groups
$$\pit{X/E}\mapright{\isom} \pit{Y/E'}$$
arises from an isomorphism of $\Z$-schemes $$X\mapright{\isom} Y.$$
As stated this is, of course, \textit{a very optimistic formulation}, and one expects that some additional hypothesis are required for the validity of this assertion. The only case  in which this conjecture is presently known is the case of geometrically connected, smooth, hyperbolic curves of Strict Belyi Type over $p$-adic fields. This case was established by Mochizuki in \cite[Corollary 2.12]{mochizuki07-cuspidalizations-proper-hyperbolic}. In \cite{joshi-gconj}, I have shown that this conjecture is false for Fargues-Fontaine curves (these are not of finite type over $p$-adic fields).

\newcommand{\fjhypxe}{\fJhypxe}
\subsection{Connectedness of $\fjxe$ and the validity of the Absolute Grothendieck Conjecture}
There is a natural notion of connectedness in the theory of Arithmetic Teichm\"uller Spaces:
\begin{defn}\label{def:connectedness}
I will say that $\fjxe$ (resp. any full subcategory $\fJ$ of $\fjxe$--for example, if $\dim(X)=1$, then one can take $\fJ=\fjhypxe$) is a \textit{connected Arithmetic Teichm\"uller Space} (resp. is connected) if for every object $(Y/E',(E'\into K, K^\flat \isom F), *_K:\sM(K)\to \yan_{E'})\in\fjxe$ (resp. $(Y/E',(E'\into K, K^\flat \isom F), *_K:\sM(K)\to \yan_{E'})\in \fJ$), one has an isomorphism  $X\isom Y$ of schemes over $\Z$.
\end{defn}

\brem 
Since isomorphism of $\Z$-schemes is an equivalence relation on schemes, one sees that the class of objects of $\fjxe$ (resp. $\fjhypxe$) can be partitioned, by this equivalence relation, into a disjoint union of connected components.
\erem

The relationship between connectedness of Arithmetic Teichm\"uller Spaces and the Absolute Grothendieck Conjecture is summarized in the following:

\bthm\label{th:groth-and-connectedness} Let $E$ be a $p$-adic field and let $X/E$ be a geometrically connected, smooth, quasi-projective variety over $E$. Then 
\benumlab
\item  If  the Absolute Grothendieck Conjecture (see \cref{ss:abs-groth-conj}) holds true for any pair of objects of $\fjxe$, then $\fjxe$ is a connected Arithmetic Teichm\"uller Space.
\item In particular if $X/E$ is additionally an hyperbolic curve of strict Belyi Type  \cite{mochizuki07} then the Arithmetic Teichm\"uller Space $\fJ^{SB}(X,E)=\fjhypxe$ is connected.
\eenum
\ethm
\bp

The first assertion is immediate from  the formulation of the Absolute Grothendieck Conjecture in \cref{ss:abs-groth-conj}. ones therefore sees that $X\isom Y$ as $\Z$-schemes and hence $\fjxe$ is connected. This proves the first assertion.

To prove the second assertion, note that one has $\fJ^{SB}(X,E)=\fjhypxe$ by \Cref{pro:sb-type-hyp}. By \cite[Corollary 2.12]{mochizuki07-cuspidalizations-proper-hyperbolic}, it is shown that the Absolute Grothendieck Conjecture holds for any pair of smooth, hyperbolic, anabelomorphic curves of Strict Belyi Type over any pair of $p$-adic fields. Hence $\fjhypxe$ is connected and the second assertion follows.
\ep

\subsection{Existence of a field of definition of $\fjhypxe$ for $X/E$ of Strict Belyi Type}\label{se:absolute-grothendieck-conj}
Assume for this section that $X/E$ is a geometrically connected, smooth hyperbolic curve of Strict Belyi Type. Then by \Cref{th:groth-and-connectedness}, the Arithmetic Teichm\"uller Space $\fjhypxe$ is connected.
In this situation, I observe below that every object of $\fjxe$ has a common field of definition. 

In \cite[Theorem 1.9]{mochizuki-topics3} Mochizuki provides a anabelian reconstruction theoretic proof of the existence of this common field over which every object of $\fjhypxe$ is defined. The existence of this common field of definition is analogous to role played by the reflex field in the context of Shimura varieties.

\bthm\label{th:field-of-def} Let $X/E$ be a geometrically connected, smooth hyperbolic curve over a $p$-adic field $E$. Assume that $X$ is of Strict Belyi Type (this implies that $X$ is definable over some number field). Let $\lmod$ be the field of moduli of $X/E$. Then for every $(Y/E',Y/K)\in \fjhypxe$,
\benumlab
\item  $Y/E'$ is of strict Belyi Type,
\item and  $Y/E'$ is definable over the number field $\lmod$.
\eenum
In particular, there is a common number field $\lmod$ which is determined by every object of $\fjhypxe$.
\ethm

\brem 
Note that  in the situation of Theorem~\ref{th:field-of-def}, by \cite[Theorem 5.4(2)]{joshi-untilts}, one sees that $\fjxe$ is a connected Arithmetic Teichm\"uller Space and hence this common number field $\lmod$ should be considered as the  field of definition of the Arithmetic Teichm\"uller Space $\fjxe$ in a manner entirely analogous to the notion of the reflex field of a Shimura variety datum (see \cite{deligne79-shimura}). 
\erem

\bp 
Let $(Y/E',E'\into K)\in\fjhypxe$. This condition means that one has an anabelomorphism
\be  
\pit{Y/E';K}\isom \pit{X/E,\C_p}.
\ee
The first assertion is proved as follows. The isomorphism of tempered fundamental groups also forces an isomorphism of the corresponding \'etale fundamental groups (see \cite[Proposition 4.4.1]{andre03}). One knows from   the proof of \cite[Proposition 2.4(iv)]{mochizuki07-cuspidalizations-proper-hyperbolic} that the property of being of strict Belyi Type is an amphoric property and hence $Y/E'$ is also of strict Belyi Type. In particular $Y/E'$ is also definable over a number field.

Now  to deduce second assertion,  one uses Mochizuki's proof of the validity of the Absolute Grothendieck Conjecture for curves of strict Belyi Type (by  \cite[Corollary 2.12]{mochizuki07-cuspidalizations-proper-hyperbolic} which was referred above in \cref{ss:abs-groth-conj}). This furnishes us with an isomorphism of absolute schemes $Y\isom X$ (over $\Z$). One also sees from \Cref{pr:top-types} that $X/E$ and $Y/E'$ have the same genus and the same number of punctures (with respect to their respective, smooth proper models).

Since $X/E$ and $Y/E'$ are both of strict Belyi type, and hence both are definable over number fields, say $L$ and $L'$ respectively.  The  isomorphism $Y\isom X$ as $\Z$-schemes provides the same moduli point corresponding to the common isomorphism class of these curves in the moduli. Let $\lmod$ be the field of moduli of $X/E$, then $L\supseteq \lmod$ and one has that  $Y/E'$ also has the field of moduli $\lmod$ (so $L'\supseteq \lmod$). So the second assertion is proved. 
\ep

\section{Untilts of absolute Galois Groups and Fundamental Groups}\label{se:untilts-of-Pi}
\subsection{Untilts of fundamental groups of Riemann surfaces}\label{se:riemann-surfaces}
It will be useful to begin with a discussion of the classical case of Riemann surfaces.
Let  $\Sigma$ be a connected Riemann surface, which one assumes to be hyperbolic to avoid trivialities. $\Pi=\pi_1^{top}(\Sigma)$ be the topological fundamental group of a connected Riemann surface $\Sigma$.
Let $\sT_{\Sigma}$ be the Teichm\"uller space of $\Sigma$.

\para The assignment $$\sT_{\Sigma}\ni (\Sigma',f) \mapsto \pi_1^{top}(\Sigma')\mapright{\isom} \Pi=\pi_1^{top}(\Sigma)$$ provides a function from the Teichm\"uller space $\sT_{\Sigma}$ to the isomorphism class of the group $\Pi$. Hence $\pi_1^{top}(\Sigma')$ is an isomorph of $\Pi$ with $(\Sigma',f)\in\sT_{\Sigma}$ serving as a geometrically distinguishable feature of this copy of $\Pi$. 

I will write $\Pi_{(\Sigma',f)}$ for the isomorph of $\Pi=\Pi_{\Sigma}$ provided by the Teichm\"uller datum $(\Sigma',f)$ and I will say \textit{$\Pi_{(\Sigma',f)}$ is an untilt of  the fundamental group $\Pi=\Pi_{\Sigma}$ of the Riemann surface $\Sigma$.}

For the anabelian contexts of \iut, the construction of the Arithmetic Teichm\"uller Space $\fjxe$ provides  similar, intrinsic arithmetic-geometric labels for isomorphs of the tempered fundamental group.

\subsection{Untilts of absolute Galois groups of $p$-adic fields}
Let $E$ be a $p$-adic field. Recall that to define the absolute Galois group $G_E$ of $E$, one requires an algebraic closure of $E$.

\begin{defn}\label{def:untilt-galois}
Let $F$ be an algebraically closed, perfectoid field of characteristic $p>0$. Let $(K\supset E,K^\flat\isom F)$ be a characteristic zero $E$-untilt of $F$. Let $\bE\subset K$ be the algebraic closure of $E$ contained in $K$ (note that $K$ is algebraically closed, perfectoid of characteristic zero, and $\bE$ is equipped with the valuation induced from $K$). I will refer to $\bE$ as \textit{the algebraic closure of $E$ provided by the given $E$-untilt $(K\supset E,K^\flat\isom F)$ of $F$.} Let 
$$G_{E,(K\supset E,K^\flat\isom F)}=\gal(\bE/E) \isom G_E$$
be the absolute Galois group of $E$ computed using the algebraic closure provided by the untilt $E$-untilt. 
I will refer to $G_{E,(K\supset E,K^\flat\isom F)}=\gal(\bE/E)$ as the \textit{untilt of the absolute Galois group $G_E$ of $E$ provided by the $E$-untilt $(K\supset E,K^\flat\isom F)$ of $F$}.  
\end{defn}

\brem 
For brevity, I will habitually shorten the notation $$G_{E,(K\supset E,K^\flat\isom F)}$$ to
$$G_{E;K}.$$
\erem
\newcommand{\sGalE}{\mathcal{Gal}_E}
The following is immediate from \Cref{def:untilt-galois}.
\bpro
Let $E$ be a $p$-adic field. Then 
\benumlab
\item Any pair of untilts of $G_E$ are isomorphic as topological groups, but there is no canonical isomorphism between them.
\item There is a category, denoted $\sGalE$, of all untilts of $G_E$ which has for its objects the collection of untilts of $G_E$ given by
$$\left\{G_{E,(K\supset E,K^\flat\isom F)} : F \text{ is some algebraically closed perfectoid field of } char(F)=p>0.\right\}$$
and whose morphisms are isomorphisms of topological groups.
\item
The category $\sGalE$ comes equipped with functors
$$\begin{tikzcd}
&G_{E,(K\supset E,K^\flat\isom F)}\arrow[dl,mapsto]\arrow[d,mapsto]\arrow[dr,mapsto]&{}\\
{} F & K	&{} (K\supset E,K^\flat\isom F) 
\end{tikzcd}
$$
to the category of perfectoid fields of characteristic $p>0$, the category of perfectoid fields of characteristic zero, the category of $E$-untilts of some perfectoid field of characteristic $p>0$ respectively.
\eenum
\epro
\bp 
The proof is clear from the construction of untilts of $G_E$.
\ep
\subsection{Untilts of tempered fundamental groups}
Now I am ready to define untilts of tempered fundamental groups.
\begin{defn}\label{def:untilts-def} Let $X/E$ be a geometrically connected, smooth, quasi-projective variety over a $p$-adic field $E$. Let $(X/E; (K\supset E,K^\flat\isom F),*_K:\sM(K)\to \xan_E)$ be an arithmetic holomorphic structure on $X/E$.
	An \emph{untilt of the tempered fundamental group $\Pi=\pit{X/E}$} is the isomorph 
	\be\label{eq:untilt-isomorph} \pit{(X/E; (K\supset E,K^\flat\isom F),*_K:\sM(K)\to \xan_E)}=\pi_1^{temp}(\xan_E,*_K: \sM(K)\to\xan_E)\ee
	of the tempered fundamental group $\pit{X/E}$ computed using the geometric base-point provided by the given  arithmetic holomorphic structure on $X/E$.
\end{defn}

To understand the Anabelian consequences of the theory of Arithmetic Holomorphic Structures, let me begin with the following assertion.

\bthm\label{thm:arith-hol-strs}
Let $X/E$ be a geometrically connected, smooth, quasi-projective variety over a $p$-adic field $E$. Let $\C_p$ be the completion of a fixed algebraic closure of $\Q_p$. Let $\pit{X/E}$ be the tempered fundamental group of $X/E$ computed using any geometric base-point $\sM(\C_p)\to \xan_E$. Then
\benumlab
\item Non-isomorphic arithmetic holomorphic structures on $X/E$ (Definition~\ref{def:arith-hol-strs}) exist (by Lemma~\ref{le:non-isom-arith-hol-strs}).
\item For any algebraically closed  perfectoid field $F$ of characteristic $p>0$, the topological isomorphism class of the pro-discrete group $$\pit{(X/E; (K\supset E,K^\flat\isom F),*_K:\sM(K)\to \xan_E)}$$ is independent of the choice of  an arithmetic holomorphic structure $$(X/E; (K\supset E,K^\flat\isom F),*_K:\sM(K)\to \xan_E)$$  on $X/E$. 
\item Hence each arithmetic holomorphic structure provides an isomorph \eqref{eq:untilt-isomorph} of  the topological groups $\pit{X/E}\supset \pit{X/\C_p}$  and an exact sequence \eqref{pr:galois-seq-untilt}.
\item Any choice of an arithmetic holomorphic structure (in the sense of Definition~\ref{def:arith-hol-strs}) provides an arithmetic holomorphic structure   in the sense of Mochizuki \iutf{IUT1, Example 1.8}:
$$\pit{(X/E; (K\supset E,K^\flat\isom F),*_K:\sM(K)\to \xan_E)}\onto G_{E;K}$$
of i.e. an isomorph of the  surjection  $$\pit{X/E}\onto G_{E}.$$ 
\item If $F$ is an algebraically closed, perfectoid field of characteristic $p>0$, then one has continuous families of isomorphs of $\pit{X/E}$ parameterized by Fargues-Fontaine curves $\syfe$ and $\sxfe$.
\item In particular one has continuous families of arithmetic holomorphic structures (in the sense of Definition~\ref{def:arith-hol-strs}) and also in the sense of \iut\ which are parameterized by Fargues-Fontaine curves $\syfe$ and $\sxfe$ for any arbitrary algebraically closed perfectoid field $F$ of characteristic $p>0$.
\eenum
\ethm

\bp 
The assertion {\bf(1)} follows from \Cref{le:non-isom-arith-hol-strs}. The assertions {\bf(2,3)} is immediate from the independence of the tempered fundamental groups on geometric base-points \cite{andre03}. The assertion {\bf(4)} is immediate from the exact sequence \eqref{pr:galois-seq-untilt} and \iutf{IUT2, Example 1.8}. Now {\bf(5)} is immediate from \Cref{def:arith-hol-strs} and the constructions of $\syfe$ \cite{fargues-fontaine}. The last assertion {\bf(6)} is immediate from {\bf(4)}.
\ep

\brem 
\Cref{thm:arith-hol-strs}, especially \Cref{thm:arith-hol-strs}{\bf(6)}, represents a substantial strengthening of the anabelian philosophy espoused by Mochizuki in \iut, as the existence of continuous parameter spaces for isomorphs of tempered fundamental groups is not suggested by \iut.
\erem

\brem One may think of an untilt of $\pit{X/E}$ as an isomorph equipped an arithmetic-geometric label (namely the datum of the arithmetic holomorphic structure giving birth to the said isomorph. This is similar to labeling procedure discussed in the classical case  \Cref{se:riemann-surfaces}.
\erem

\bcor\label{cor:untilted-fun-grp}
Let $X/E$ be a geometrically connected, smooth, quasi-projective variety over a $p$-adic field $E$. Then the class of untilts of the tempered fundamental group of $X/E$ is the class of (isomorphic) topological groups
\be
\left\{ \pit{(X/E; (K\supset E,K^\flat\isom F),*_K:\sM(K)\to \xan_E)}\right\}\ee  where  $F$  is some algebraically closed, perfectoid field of  characteristic  $p>0$,  and  where  $(X/E; (K\supset E,K^\flat\isom F),*_K:\sM(K)\to \xan_E)$ is some arithmetic holomorphic structure on  $X/E$.
\ecor

\brem 
The above corollary is the precise version of Mochizuki's Key Principle of Inter-Universality \iutf{IUT1, \ssep I3, Pages 25--26} which according to Mochizuki lies the foundation of the anabelian theory described in \iut.
\erem

\brem In this paper, for simplicity of notation, I will habitually contract the notation 
$$\pit{(X/E; (K\supset E,K^\flat\isom F),*_K:\sM(K)\to \xan_E)}$$ to $$\pit{X/E; K}.$$  
\erem

\section{Galois cohomology and Arithmetic Teichm\"uller Spaces}\label{ss:galois-group-cohomology}
\newcommand{\Zh}{\widehat{\Z}}
In this section I want to outline some consequences, for Galois cohomology, of the existence of Arithmetic Teichm\"uller Space $\fjxe$ constructed in \ssep\ref{se:construct-att}. As an  application one obtains a robust version of Mochizuki's idea, in \iutf{IUT3, Section 3}, of collating Galois cohomology classes arising from distinct arithmetic holomorphic structures. This collation process is central for the formulation of \moccor\ and for the volume computations of \iutf{IUT4}, \cite{joshi-teich-abc-conj}. However, the present section is devoted to developing the local formalism which facilitates this collation of Galois cohomology classes, while the adelic case (which is needed for arithmetic applications) is detailed in \cite{joshi-teich-def}, \cite{joshi-teich-IIIB-rosetta}. This allows me to independently and explicitly demonstrate many assertions of \iut\ dealing with this collation process--notably \moccor.

\subsection{} Let $X/E$ be a geometrically connected, smooth, quasi-projective variety over a $p$-adic field. Then by \Cref{def:arith-hol-space-local} one has the arithmetic Teichm\"uller space $\fjxe$ associated with  $X/E$. For applications to \moccor, one assumes additionally that $X/E$ is a hyperbolic curve.
Choose an arithmetic holomorphic structure $(X/E,(\C_p\supset E,\cpt\mapright{id} \cpt),*_{\C_p}:\sM(\C_p)\to\xan_E)$ where $\cpt \mapright{id} \cpt$ is the identity isometry of $\cpt$.  The $E$-untilt $(\C_p\supset E,\cpt\mapright{id} \cpt)$ is a closed classical point of $\syfqpe{\cpt}$ lying over the canonical point of $\sxQp$ for the composite morphism $\syfqpe{\cpt}\to \sX_{\cpt,E}\to\sxQp$. I will refer to this as the \textit{standard arithmetic holomorphic structure on $X/E$}.

\subsection{} Let $(Y/E',(E'\into K, K^\flat \isom F), *_K:\sM(K)\to \yan_{E'})\in\fjxe$. This comes equipped with a preferred algebraic closure of $\bE'$ of $E'$ namely the algebraic closure $\bE'\subset K$ contained in the given algebraically closed perfectoid field $K$. Hence one also has $G_{E';K}=\gal(\bE'/E')$, of course (as abstract topological groups one has $G_{E';K}\isom G_{E'}$ but I am emphasizing the choice of the preferred algebraic closure $\bE'$). By \Cref{pr:fund-group-same} one has an isomorphism of tempered fundamental groups $\pit{Y/E';K}\isom \pit{X/E;\C_p}$ and also of geometric tempered fundamental subgroups
$\pit{Y/K}\isom \pit{X/\C_p}$ and an isomorphism $G_{E'}\isom G_E$ (i.e. $E'$ and $E$ are anabelomorphic $p$-adic fields).

\subsection{}\label{ss:distinct-kummer-theories} Let  $\mu(\bE')\subset \bE'\subset K$ the (multiplicative) group of all the roots of unity contained in $\bE'$ and let $\mu_p(\bE')\subset \bE'\subset K$ the (multiplicative) group of all the $p$-power roots of unity contained in $\bE'$. Write $\Zh(1)_K\subset K$ (resp. $\Z_p(1)_K$) for the $G_{E';K}$-module of roots of unity (resp. the $G_{E';K}$-module of $p$-power roots of unity) contained in $\bE'_K\subset K$. The notation  $\Zh(1)_K\subset K$ is to emphasize that this $G_{E';K}$-module is contained in the $G_{E';K}$-module $\bE'_K$. Similarly one has the $G_{E';K}$-modules $\Q(1)_K$ and $\Q_p(1)_K$. One may consider these $G_{E';K}$-modules as $\pit{Y/E';K}$-modules through the quotient homomorphism $\pit{Y/E';K}\onto G_{E';K}$. Write $$G_{E';K}\act\Z(1)_K,\ G_{E';K}\act\Z_p(1)_K,\ \pit{X/E';K}\act\Z(1)_K,\  \pit{X/E';K}\act\Z_p(1)_K$$ etc. for these modules.
These modules in turn provide continuous cohomology groups
$$H^i(G_{E';K}, \Z(1)_K),\ H^i(G_{E';K},\ \Z_p(1)_K),\ H^i(\pit{X/E';K},\Z(1)_K), \text{etc.}$$
I will refer to these modules as the \emph{Kummer Theory and Galois cohomology for the $p$-adic field $E'$ provided by the arithmetic holomorphic structure $(Y/E',(E'\into K, K^\flat \isom F), *_K:\sM(K)\to \yan_{E'})\in\fjxe$.} These considerations may also be applied to various subspaces of these cohomology groups $H^1_e(G_{E'},-),H^1_f(G_{E'},-),H^1_g(G_{E'},-)$ considered in \cite{bloch90}.

Now let $(X/E,(\C_p\supset E,\cpt\mapright{id} \cpt),*_{\C_p}:\sM(\C_p)\to\xan_E)$ be the standard arithmetic holomorphic structure chosen above. Then I will omit the reference to the algebraically closed, perfectoid field from the above notations and simply write $G_E=\gal(\bE/E)$, and $\Zh(1)\subset \C_p$ (resp. $\Z_p(1)\subset \C_p$) and $H^i(G_E,\Z(1))$, $H^i(G_E,\Z_p(1))$ etc. and refer to this as the \textit{standard Kummer Theory and Galois cohomology} etc.  Thus each arithmetic holomorphic structure in $\fjxe$ provides us with a distinguished Kummer Theory. To summarize:

\bpro Given any arithmetic holomorphic structure $(Y/E',(E'\into K, K^\flat \isom F), *_K:\sM(K)\to \yan_{E'})\in\fjxe$ one has, for each integer $i\geq 0$,  a Kummer Theory and continuous Galois cohomology given by the assignment:
\be 
G_{E';K}\act\Z(1)_K,\ G_{E';K}\act\Z_p(1)_K,\ \pit{X/E';K}\act\Z(1)_K,\ \pit{X/E';K}\act\Z_p(1)_K
\ee
\be 
(Y/E',\yan/K)\mapsto H^i(\pit{Y/E';K},\Zh(1)_K)
\ee
and also
\be 
(Y/E',\yan/K)\mapsto H^i(G_{E';K},\Zh(1)_K),
\ee
and a similar construction for the subspaces $H^1_e(G_{E'},-),H^1_f(G_{E'},-),H^1_g(G_{E'},-)$.
\epro

\subsection{} One has the following:
\bthm\label{pr:galois-cohomology-and-teich} 
Let $X/E$ be a geometrically connected, smooth, quasi-projective variety over a $p$-adic field $E$. Let $(Y/E',(E'\into K, K^\flat \isom F), *_K:\sM(K)\to \yan_{E'})\in\fjxe$ and let $(X/E,(\C_p\supset E,\cpt\mapright{id} \cpt),*_{\C_p}:\sM(\C_p)\to\xan_E)$ be the standard arithmetic holomorphic structure on $X/E$. Then
\benumlab
\item Any anabelomorphism $\pit{X/E;K}\isom \pit{X/E;\C_p}$ provides, for all integers $i\geq 0$, an isomorphism of continuous cohomology groups
$$
H^i(\pit{Y/E';K},\Zh(1)_K) \isom H^i(\pit{X/E;\C_p},\Zh(1))
$$
and also an isomorphism
$$
H^i(G_{E';K},\Zh(1)_K) \isom H^i(G_{E},\Zh(1))
$$
and also an isomorphism (of $\Q_p$-vector spaces)
$$
H^i(G_{E';K},\Q_p(1)_K) \isom H^i(G_{E},\Q_p(1)), 
$$
\item and a similar isomorphism holds for the subgroups (resp. subspaces) $H_e^1,H^1_f,H^1_g$ of these  groups (resp. subspaces), for example:
$$
H^i_f(G_{E';K},\Q_p(1)_K) \isom H^i_f(G_{E},\Q_p(1)).
$$
\eenum
\ethm
\bp
It is enough to remark that by \cite[Proposition 4.2(iv)]{hoshi-mono} (or \cite[Proposition 1.2.1]{mochizuki04})  the anabelomorphism $G_{E';K}\isom G_E$ induces an isomorphism of $\Zh(1)_K\isom \Zh(1)$ which is compatible with $G_{E';K}$ (resp. $G_E$-action) on either side. This proves the first assertion. The second assertion is proved similarly noting that these subgroups (resp. subspaces) are amphoric.
\ep
\subsection{Collation of cohomology classes}
\newcommand{\sS}{\mathcal{S}}
The following corollary of Theorem~\ref{pr:galois-cohomology-and-teich} will play a crucial role. It asserts that Galois cohomology classes provided by distinct arithmetic holomorphic structure may be collated in the standard Galois cohomology:
\bcor\label{cor:collation-classes} 
Let $i\geq 0$ be an integer. Let $\sS$ be  a collection of objects of $\fjxe$. For each $s\in \sS$, let $$\Psi^i_s\subset H^i(G_{E,s},\Zh(1)_s)\ \ (\Psi^i_s\subset H^i_f(G_{E,s},\Zh(1)_s) \text{ resp.} )$$ be a collection of Galois cohomology classes provided by each object $s\in\sS$. Then these classes may be assembled into a subset  $$\Psi\subset H^i(G_E,\Zh(1))\ \  (\Psi^i_f\subset H^i_f(G_{E},\Zh(1))\text{ resp.})$$
of the relevant standard Galois cohomology.
\ecor
\bp 
Let $\Psi\subset H^i(G_E,\Zh(1))$ ($\Psi^i_s\subset H^i_f(G_{E,s},\Zh(1)_s)$ resp.) be the union of the images of each element of $\Psi_s$ (for each $s\in\sS$) under all the isomorphisms provided by Theorem~\ref{pr:galois-cohomology-and-teich}{\bf(2)}.
\ep

\brem
The adelic version of this corollary (established in \cite{joshi-teich-def}, \cite{joshi-teich-IIIB-rosetta}) is the key principle which underlies Mochizuki's construction of his $\Theta$-values set in \iutf{IUT3}. For \moccor, one uses a specific choice of $\sS$ (for each prime of a fixed number field) whose construction requires the notions of $\Theta_{gau}$-links and $\mathfrak{log}$-links, both of which I have established in \cite{joshi-teich-estimates}. 
\erem

\section{Adelic Arithmetic Teichm\"uller Spaces} 
In this section, I want to indicate how the theory of arithmetic Teichm\"uller Spaces extends to an adelic theory in the presence of a number field. The constructions presented here are complemented and strengthened by the fact that a fixed number field itself has an Arithmetic Teichm\"uller Space.  The adelic theory of a number field is glued together by means of the product formula and this formula comes to play here as well. The  construction of an arithmetic Teichm\"uller Theory of a fixed number field is detailed in \cite{joshi-teich-def} and it clearly demonstrates the role the product formula of \cite{artin45} plays in the theory. This aspect also comes into play in \cite{joshi-teich-IIIB-rosetta}, \cite{joshi-teich-IIIC-moccor} which apply the said global+adelic theory to demonstrate the assertions of \iutf{IUT3, IUT4} using the theory of Arithmetic Teichm\"uller Spaces developed here and \cite{joshi-teich-def}.

\subsection{} Let $X/L$ be a geometrically connected, hyperbolic, smooth and quasi-projective variety over a  number field $L$ and let $\wp$ denote a non-archimedean prime of $L$ and let $\infty_1,\ldots,\infty_n$ be all the archimedean primes of $L$. Let $L_\wp$ (resp. $L_{\infty_i}$) be the 
completions of $L$ at $\wp$ (resp. $\infty_i$). By an \emph{\locvar} is the following collection: 
\benumlab
\item If $\wp$ is a non-archimedean place of $L$, then let
\begin{align*}
\fJ(X,L,\wp) &= \fJ(X,L_\wp), \text{ and},\\
\fJ(X,L,\wp)_F &=\fJ(X,L_\wp)_F.
\end{align*}

\item If $\wp=\infty_i$ for some $i$ (i.e. $\wp=\infty_i$ is an archimedean place of $L$) then let $$\fJ(X,L,\infty_i)=\fJ(X,L_{\infty_i}).$$
\eenum
Similar definition can be made of $\fJ(X,L,\wp)_F$. For $\wp={\infty_i}$ one takes  $\fJ(X,L,\wp)_F= \fJ(X,L,\infty_i)$ purely for notational symmetry.

\para
Properties of $\fJ(X,L,\wp)=\fJ(X,L_\wp)$ (resp. 
$\fJ(X,L,\wp)_F =\fJ(X,L_\wp)_F$ can be gleaned from the results established in \Cref{se:construct-att}, \Cref{se:connectedness}.

\begin{defn}
Let $X/L$ be a geometrically connected, smooth, quasi-projective variety over a number field $L$ with no real embeddings (I will make this restriction to avoid notational complexity). Then the \emph{adelic tempered fundamental group of $X/L$} is the group (equipped with product topology):
$$\widetilde{\pit{X/L}}=\prod_{0\neq\wp\in\spec(\O_L)}\pit{X/L_\wp}\times \prod_{i=1}^{n} \pi_1^{top}(X/L_{\infty_i}).$$ 
\end{defn}

The adelic Arithmetic Teichm\"uller Space will be an anabelian variation providing the adelic tempered fundamental group $\widetilde{\pit{X/L}}$.

\para In the global i.e. number situation, the local $p$-adic Teichm\"uller Spaces can be assembled into a global category.
\bthm\label{th:main4}
Let $X/L$ be a geometrically connected, smooth, quasi-projective, hyperbolic variety over a number field $L$. Assume $L$ has no real embeddings. Then there exists a category 
$\fTeich$, called  the \emph{Arithmetic Teichm\"uller Space associated to $X/L$} which  has the following properties:
\benumlab
\item $\fTeich$ is given as a product category:
$$\fTeich=\prod_{\wp}\fJ(X,L_\wp)$$
where $\wp$ runs over all the inequivalent, non-trivial valuations of $L$ and where $\fJ(X,L_\wp)$ is the $p$-adic Teichm\"uller Space associated to $X/L_\wp$ constructed in Theorem~\ref{th:main4.5-const-local-teich}.
\item $\fTeich$ is an anabelian variation providing the adelic tempered fundamental group 
$$\widetilde{\pit{X/L}}=\prod_{0\neq\wp\in\spec(\O_L)}\pit{X/L_\wp}\times \prod_{i=1}^{n} \pi_1^{top}(X/L_{\infty_i})$$
i.e. each object of $\fTeich$ provides an isomorph of $\widetilde{\pit{X/L}}$ and two such isomorphs may be distinguished from each other by means of the arithmetic-geometric labels provided by the objects of $\fTeich$ which give rise to the isomorphs.
\eenum
\ethm

\subsection{Adelic Theory and Diamonds}
\textit{This is a continuation of \Cref{re:diamonds} and should be read along with it.} Let $X/L$ be a geom. connected quasi-projective variety over a number field $L$ and for a prime $v\in\mathbb{V}_L$, let $L_v$ be the $v$-adic completion of $L$ and let $\xan_{L_v}$ be the analytic space associated to $X/L_v$. Then the Arithmetic Teichm\"uller Space (\Cref{def:arith-hol-space-local}) contains sets of the form
$$\left\{\left(X^{an}_{L_v}, (L_v\into K_v,\iota_v:K_v^\flat\isom F_v), *_{(K_v,K_v^\flat\isom F_v^\flat)}:\sM(K_v)\to \xan_{L_v}\right) \right\}/\text{suitable equivalence},$$
where $F_v$ is an algebraically  closed, perfectoid field of characteristic $p_v>0$ if  $v$ is non-archimdean
otherwise $F_v=\C$. 
The precise relationship between this and $X_{L_v}^\Diamond$ is given by \Cref{pr:local-diamond}. Thus working with the adelic Arithmetic Teichm\"uller Space $\fJ(X,L)$ is approximately akin to considering $$\prod_v  X_{L_v}^\Diamond.$$
\textit{But beware that \Cref{def:arith-hol-space-local} is more general than $X_{L_v}^\Diamond$.}

\section{Relationship to Mochizuki's Inter-Universal Teichm\"uller Theory}\label{se:relation-to-iut}
\para \emph{A more detailed treatment of the relationship between the theory of the present paper and that of \iut\ is found in the `Rosetta Stone' provided in \cite{joshi-teich-IIIB-rosetta}.} In this brief section, I show that the Arithmetic Teichm\"uller Theory of \present\ is fully compatible with \iut. Notably, the central deficiency \iutthr\ is that it provides no way of distinguishing between two distinct arithmetic holomorphic structures. As one knows from \Cref{thm:arith-hol-strs} this is not an issue in the present theory.

\subsection{Mochizuki's Key Principle of Inter-Universality}\label{ss:key-principle} Mochizuki's \iut\ asserts the existence of a Teichm\"uller Theory in the arithmetic context. The central principle which is key to \iut\ is \textit{Mochizuki's Key Principle of Inter-Universality \iutf{IUT1, \ssep I3, Pages 25--26}}. 

This principle asserts that the theory of \iut\ is a theory which works with the tempered fundamental groupoid (rather than the tempered fundamental group). Mochizuki's Principle is treated here in detail in \Cref{appendix-univ}. Notably,  a key assertion Mochizuki makes as a consequence of his principle is stated without proof in \iutf{IUT1, \ssep I3, Pages 25--26}. A proof that  consequence of his Key Principle is proved here as  \Cref{thm:zfc}. 

More precisely, {Mochizuki's Key Principle of Inter-Universality \iutf{IUT1, \ssep I3, Pages 25--26}} requires one to work with  tempered fundamental groups arising from arbitrary geometric base-points and also asserts that the domains and the codomains of the key operations of \iut\ are required to have distinct geometric base-points (in general). 

In the notation of this paper, Mochizuki's Key Principle of Inter-Universality requires one to work with 
a category whose objects are 
$$(\pit{X/E;K}, K, \text{``Frobenius Picture''}),$$
where the algebraically closed perfectoid field is provided by the $K$-geometric base-point used to calculate $\pit{X/E;K}$. The `Frobenius Picture' aspect is detailed in the `Rosetta Stone' \constrthrB{\ssep}{ss:mochizuki-dichotomy} and will not be recalled here (at this juncture it will suffice to say that this aspect is naturally incorporated in \Cref{def:arith-hol-strs}). On the other hand,  the perfectoid fields $K$ are not amenable to the anabelian reconstruction methods of \topics, \iutthr, and are completely ignored (while requiring them). This is the second important deficiency of \iutthr--namely \iutthr\ provides no clear way of disambiguating between distinct arithmetic holomorphic structures (while \iutf{IUT3, Theorem 3.11, Corollary 3.12} clearly rest upon disambiguating such structures). [This central difficulty of \iutthr\ disappears with a cleaner formulation of this notion provided by \Cref{def:arith-hol-strs}.] It should be clear to readers familiar with the theory of algebraically closed perfectoid fields that dealing with the missing information \textit{necessarily} requires one to take the path considered in \present.

In the present paper, I work with the category $\fjxe$ and  more precise objects $(Y/E',(E'\into K, K^\flat \isom F), *_K:\sM(K)\to \yan_{E'})\in\fjxe$ provided by \Cref{def:arith-hol-space-local}. By \Cref{cor:untilted-fun-grp} one obtains the natural class of topological groups
$$\pit{(Y/E'; (K\supset E',K^\flat\isom F),*_K:\sM(K)\to \yan_{E'})},$$
all of which are isomorphic to $\pit{X/E}$, and its relationship to Mochizuki's Theory and his Key Principle of Inter-Universality is established by \Cref{thm:arith-hol-strs}.

The present paper works with more precise geometric object of $\fjxe$ (\Cref{def:arith-hol-space-local}) and untilts of fundamental groups provided by \Cref{def:untilts-def} from which the functor to Mochizuki's category is obvious. 

These results should make applications of the theory of the \present\ to \iut\ unsurprising. \emph{However  there are  some important  differences between the two approaches--these are discussed in \ssep\ref{pa:diff-from-iut}.}

\para Let me remark that in \iut\ Mochizuki works with multiplicative groups as the anabelian approach considered in loc. cit. is inherently multiplicative. On the other hand, the approach taken in \Cref{th:main3} via \cite{fargues-fontaine} is necessarily additive and not  multiplicative. In the next few paragraphs I provide a translation between the two. This allows one to construct functors from $\fjxe$ to various categories used in \iut.

\para For the multiplicative description let me fix some notations. Let $E=\Q_p$, let $\sG/\Z_p$ be the Lubin-Tate formal group with formal logarithm given by $\sum_{n=0}^\infty \frac{T^n}{p^n}$, $F=\cpt$. Let $\bE=\bQ_p$ be the algebraic closure of $E=\Q_p$ in $\C_p$. 

\para\label{pa:artin-hasse-exp} The Artin-Hasse Exponential provides the following:
\blem\label{le:artin-hasse-exp}  Let $\sG$ be the Lubin-Tate formal group over $\Z_p$ with logarithm $\sum_{n=0}^\infty \frac{T^{p^n}}{p^n}$. Let $Exp_{AH}(T)$ be the Artin-Hasse exponential function. Then the homomorphism $a\mapsto Exp_{AH}(a)$ provides a natural isomorphism of topological $\Z_p$-modules $$Exp_{AH}:\sG(\O_{\cpt})\isom \hgm(\O_{\cpt})$$ and hence also of $$ \widetilde{\sG(\O_{\C_p})}\isom\sG(\O_{\cpt})\isom \hgm(\O_{\cpt})\isom \widetilde{\fgmcp}.$$
\elem
\bp 
See \cite[Example 4.4.7]{fargues-fontaine}.
\ep

\para Let $\hgm/\Z_p$ be the multiplicative formal group. Then one has for the multiplicative formal $\hGm$ one has 
\be\label{pa:multiplicative-group-facts1} 
\hGm(\O_{\C_p})=1+\frak{m}_{\O_{\C_p}} 
\ee
where $1+\frak{m}_{\O_{\C_p}} \subset \O_{\C_p}^*$ is the subgroup of units congruent to $1$ modulo the maximal ideal $\frak{m}_{\O_{\C_p}}\subset\O_{\C_p}$, and one also has from this that
\be\label{pa:multiplicative-group-facts2} 
\widetilde{\fgmcp}  = \left\{(x_n)_{n\in\Z}:x_n\in\hgm(\O_{\C_p})=1+\fm_{\C_p}, x_{n+1}^p=x_n\forall n\in\Z\right\} \ee

This fits into an exact sequence of $\Q_p$-Banach spaces (\cite[Proposition 4.5.14]{fargues-fontaine})
$$0\to T_p(\hgm)\tensor_{\Z_p}\Q_p\to \widetilde{\hgm(\O_{C_p})}\mapright{\log_{\hgm}} \C_p\to 0,$$
where $T_p(\hgm)$ is the $p$-adic Tate-module of $\hgm$ (note that $T_p(\hgm)$ is a rank one free $\Z_p$-module computed using the $p$-power roots of unity contained in $\C_p$), and where $\log_{\hgm}$ is the logarithm of the formal group $\hgm$. Explicitly $\log_{\hgm}$ is given in terms of the $p$-adic logarithm as follows (\cite[Proof of Proposition 4.5.9]{fargues-fontaine}). Let $x=(x_n)_{n\in\Z}\in\widetilde{\hgm(\O_{C_p})}$ then 
$$\log_{\hgm}(x)=\log(x_0),$$
where $$\log:\O_{\C_p}^*\to \C_p$$ is the $p$-adic logarithm.

\para\label{pa:times-mu-note} For a valued field $K\supset \Q_p$ let where $$\mu(K)\subset \O_K^*$$ be the (topological) subgroup of roots of unity in $K$ and write
\be\O_{K}^{\times\mu}=\O_{K}^*/\mu(K)\ee
for the quotient of $\O_K^*$ by $\mu(K)$.
Let \be\O_{K}^{\triangleright}=\O_K-\{0\}\ee be the multiplicative monoid of non-zero elements of $\O_K$. Both these notations are introduced and used extensively in \iut.

\para Let $\mu(\C_p)\subset \C_p^*$ be the subgroup of roots of unity contained in $\C_p$. Then one has the exact sequence of topological groups
$$0\to \mu(\C_p)\to \O_{\C_p}^*\mapright{\log} \C_p\to 0.$$

\brem  In \iut\ especially \iutf{IUT3}, Mochizuki works with $$\O_{ \bQ_p}^{\times \mu}=\O_{ \bQ_p}^*/\mu(\bQ_p)$$
and the exact sequence
$$1\to \mu(\bQ_p) \to \O_{ \bQ_p}^* \mapright{\log} \bQ_p\to 0.$$ 
Since one is dealing with $p$-adic logarithms, I prefer to work with the corresponding sequence obtained for the complete field $\C_p$ as opposed to $\bQ_p$. From the anabelian reconstruction point of view taken in \iut, complete, algebraically closed valued fields are not determinable from the absolute Galois group. This is one key difference between the theory of the present paper and Mochizuki's work. 
\erem

\blem
The inclusions $\O_{\bQ_p}^*\subset \O_{\C_p}^*$ and $\O_{\bQ_p}^{\times\mu}\subset \O_{\C_p}^{\times\mu}$ are dense inclusions.
\elem
\bp 
The first assertion is standard and it implies the second assertion as $\mu(\bQ_p)=\mu(\O_{\C_p})$.
\ep

\blem
One has an exact sequence  of topological $G_{\Q_p}$-modules 
$$1\to \mu_p(\bQ_p) \to (1+\fm_{\bQ_p})\to \O_{ \bQ_p}^{\times\mu}\to 1.$$
\elem
\bp 
Let $\mu'(\bQ_p)\subset \mu(\bQ_p)$ (resp. $\mu_p(\bQ_p)\subset \mu(\bQ_p)$) be the subgroup of roots of unity with orders coprime to $p$ (resp. the subgroup of roots of unity of order a power of $p$). Then one has
$$\mu'(\bQ_p)\times\mu_p(\bQ_p)\isom \mu(\bQ_p).$$ 
Note that for every $n\geq 1$, $$1-\zeta_{p^n}\cong0 \bmod{\fm_{\C_p}}.$$ So any $p$-power root of unity is contained in the group of $1$-units $1+\fm_{\bQ_p}$ and hence $$\mu_p(\bQ_p)\subset 1+\fm_{\C_p}.$$
Hence $$\O_{\bQ_p}^*=\mu'(\bQ_p)\times (1+\fm_{\bQ_p}),$$
and hence by definition $$\O_{ \bQ_p}^*/\mu(\bQ_p)\isom (1+\fm_{\bQ_p})/\mu_p(\bQ_p)\isom \O_{ \bQ_p}^{\times\mu},$$
and this provides the asserted exact sequence.
\ep

\para\label{pa:diff-from-iut} Let me explain the key difference between the theory described here and that of \iut. In \iut\ Mochizuki works with the pair $G_E$ and its action on $\O_{\bQ_p}^{\times\mu}$, that is  with $$G_E\act\O_{ \bQ_p}^{\times\mu}.$$ In \iut,  roughly speaking,   algorithms of Anabelian Reconstruction Theory,  automorphisms of $G_E\act\O_{ \bQ_p}^{\times\mu})$, the theory of log-link and theta-link are used to produced variation in the data of arithmetic line bundles. 

The present paper  can also be read in the multiplicative context using the isomorphism (\ssep\ref{pa:artin-hasse-exp})
$$G_E\act\widetilde{\hgm(\O_{\C_p})}\isom G_E\act\hgm(\O_{\cpt}),$$
and  variation of the data of arithmetic line bundles arises from existence Arithmetic Teichm\"uller Spaces (Theorem~\ref{th:main4.5-const-local-teich}) which arise from existence of deformations of analytic structure of $\xan_{\C_p}$ via deformations of $\C_p$.

\para The Virasoro type symmetry action is detailed in \Cref{th:main3} and \Cref{co:action-on-cpt-cat}; its relationship to Mochizuki's Indeterminacy `Ind2' is discussed in the `Rosetta Stone' \cite{joshi-teich-IIIB-rosetta}.

\para Let me now show that Theorem~\ref{thm:main}, Theorem~\ref{thm:main2}, \Cref{th:main3} provide functors to various anabelian categories considered in \iut. Notably these theorems show that there are geometrically distinguishable isomorphs of the tempered fundamental groups. The Arithmetic Teichm\"uller Space $\fjxe$ constructed in Theorem~\ref{th:main4.5-const-local-teich} includes all hyperbolic curves with tempered fundamental group topologically isomorphic to that of a given hyperbolic curve $X/E$.

\newcommand{\otmup}{\O_{\bE'}^{\times\mu}}

\para Let $Y/E'$ be a smooth, quasi-projective variety over a $p$-adic field $E'$ and let $\bE'$ be an algebraic closure of $E'$.

\para\label{pa:anab-functors} \textcolor{red}{This section is taken from `Rosetta Stone' \cite{joshi-teich-IIIB-rosetta} where readers will find a detailed discussion of \iut\ and its relationship to the \present.} Consider an arbitrary triple $(Y/E',E'\into K)\in \fjxe$. Then as $K$ is an algebraically closed, perfectoid field of characteristic zero,  consider the algebraic closure $\bE'\subset K$ (as a valued fields) of $E'\into K$.  Let $G_{E';K}$ be the absolute Galois group of $\bE'/E$. The one has the tautological action of $G_{E',K}$ on the various multiplicative submonoids of $\bE'$.  One has the following  functors from $\fJ(X,E)$ to various categories (called \textit{multiradial environments}) given by \iutf{IUT2, Example 1.8, Page 247} and used throughout \iut:
\benumlab 
\item 
$$
(Y/E',E'\into K)\mapsto \pit{X/E;K}\onto G_{E';K}
$$
\item and
$$
(Y/E',E'\into K)\mapsto \pit{Y/E';K}\act[K] \otmu,
$$ 
\item and
$$ (Y/E',E'\into K)\mapsto \pit{Y/E';K}\act[K] \O_{\bE'}^{\triangleright},$$
\item and
$$ (Y/E',E'\into K)\mapsto \pit{Y/E';K}\act[K] \bE^{'*}.$$
\eenum
[Here $\pit{Y/E';K}\act[K] \O_{\bE'}^*$ means that the field $K$ is used to compute the algebraic closure of $E'$, and the action of $\pit{Y/E';K}$ on $\O_{\bE'}^*$  through the quotient $\pit{Y/E'}\to G_{E'}$ is computed using the algebraic closure  $\bE'\subset K$. ]

Note that one has an isomorphism of topological groups $\pit{Y/E';K}\isom \pit{X/E}$. 

\brem
More details can be found in \cite{joshi-teich-estimates}, \cite{joshi-teich-IIIB-rosetta}.
\erem

\section{Appendix: On the Principle of Inter-Universality in \iut}\label{appendix-univ}
\nws
\newcommand{\sB}{\mathcal{B}}
\newcommand{\syqp}{\sY_{\cpt,\Q_p}}

The purpose of this appendix is to establish \emph{Mochizuki's Key Principle of Inter-Universality} which is stated without proof in \iutf{IUT1, \ssep I3, Pages 24-27}. As Mochizuki informs us, this is the key principle is the mathematical firmament upon which \iut\ is grounded.  

Specifically, Mochizuki's Key Principle asserts that working with a tempered fundamental group, and a geometric base-point used for its construction, yield a Grothendieck Universe.  A proof of this assertion is not provided in \iut, but is given here in \Cref{thm:zfc}. The key idea which Mochizuki enshrines in this principle is that these universes (arising from distinct geometric base-points) should serve as proxies for Arithmetic Holomorphic Structures, and this idea underlies \iut. Readers will find my discussion  of how geometric base-points for \'etale fundamental groups  lead to algebraic loop-spaces given in \cite[Answwers to Questions 16, 17 and 18]{joshi-faq} very useful.

However, as I have noted elsewhere \cite[Question 31]{joshi-faq}, Grothendieck universes are not actually used in the main portion of the proof of the $abc$-conjecture asserted in \iut--but without recognizing this Key Principle and its consequences, the key assertions of \iut, such as \moccor, and their proofs can appear to be problematic. The strategy employed in the present paper to avoid these difficulties  is to explicitly define Arithmetic Holormorphic Structures and Arithmetic Teichm\"uller Spaces.

But I want to demonstrate in \Cref{thm:zfc} (given below) that Mochizuki's Key Principle of Inter-Universality  does lead to models of set theories as asserted by Mochizuki. The techniques of the proof of this theorem are certainly well-known among topos theorists (so I am not claiming anything original here) but I give a proof here for completeness. \textit{However,} readers should beware that  models of ZFCs are not used in the main portion of \iut.  

Let $X/E$ be a geometrically connected, smooth hyperbolic curve (e.g. the affine curve $y^2=x(x-1)(x-p)$).
Consider an arithmetic holomorphic structure (\Cref{def:arith-hol-strs}). This consists of the following data:
\be \label{eq:geom-base-point}y=(X/E,\ \overbrace{\left(E\into K_y,K_y^\flat\isom \cpt\right)}^{\text{an untilt of }\cpt},\  \overbrace{*_{K_y}:\sM(K_y)\to \xan_{E}}^{\text{geometric base-point}}\ ).\ee
This gives  a point of the local Arithmetic Teichm\"uller Spaces (\Cref{def:arith-hol-space-local}) i.e. an object of the category $\fjxe$ given by (\Cref{def:arith-hol-space-local}).

More generally, applying \Cref{cor:untilted-fun-grp} and \Cref{def:arith-hol-space-local} to each object  $(Y/E',(E'\into K, K^\flat \isom F), *_K:\sM(K)\to \yan_{E'})\in\fjxe$,  one obtains the natural class of topological groups
$$\pit{(Y/E'; (K\supset E',K^\flat\isom F),*_K:\sM(K)\to \yan_{E'})}$$
all of which are isomorphic to $\pit{X/E}$. The relationship to Mochizuki's Theory and his Key Principle of Inter-Universality is established by \Cref{thm:arith-hol-strs}.

In \eqref{eq:geom-base-point}, one thinks of the untilt datum $(E\into K_y,K_y^\flat\isom \cpt)$ given by $y$ as arising from a closed classical point of the Fargues-Fontaine curve $y\in \syfqpe{\cpt}$, with residue field $K_y$. [As is pointed out in \cite[Remark 2.3]{joshi-untilts}, \cite[Proposition 4.3.1]{joshi-teich}, $y$ also arises (via Huber's functor to Berkovich spaces) from some 
\be\label{eq:geom-base-point2} x=((E\into K_y,K_y^\flat\isom \cpt),\  *_{K_y}:{\rm Spa}(K_y)\to X^{ad}_{E})\in X^\Diamond(\cpt),
\ee where $X^\Diamond$ is the diamond of $X/E$ constructed in \cite{scholze-diamonds}.]
Let $$\Pi_y=\pi_1^{temp}(\xan_{E},*_{K_y}:\sM(K_y)\to \xan_{E})$$ be the tempered fundamental group of $\xan_{E}$ computed using the datum of the geometric base-point given by $y$. From \cite{andre03}, one knows that $\Pi_y$ is a pro-discrete group.
\iut\ does not prove the precise relationship between Mochizuki's Key Principle of Inter-Universality and Universes or Set Theories. I claim the following:
\bthm\label{thm:zfc} Let the notation and assumptions be as in the above paragraph.
Each $y$ \eqref{eq:geom-base-point} gives rise to a model  $\mathrm{ZF}_y$ of Zermelo-Fraenkel set theory with Zermelo-Fraenkel axioms and an inner model $\mathrm{ZFC}_y$ of Zermelo-Fraenkel set theory plus axiom of choice.
\ethm

\bp Let $\sB_y=B\Pi_y$ be the category of sets equipped with discrete topology and continuous action of $\Pi_y$. Then I claim that for each $y$, $\sB_y$ is an atomic, boolean Grothendieck topos. Here boolean topos means its internal logic is Classical Logic and not Intuitionistic Logic. Now it has been well-understood from the work of William Lawvere that a topos may be considered as a universe in which to formulate mathematics. [So this is where universes arise in Mochizuki's context.]

That $\sB_y$ has these properties is well-known among topos theorists. From \cite[Chap. V, Section 6, Theorem 1]{macLane_moerdijk1992} one knows that $\sB_y$ is a topos;
from  \cite[III, Prop. 1.2.8, Page 110]{Giraud1971} (or other works) one knows that $\sB_y$ is a Grothendieck topos. From \cite[Section 7 (2)]{barr1980} (or other works) one knows that $\sB_y$ is an atomic topos; from  \cite[Theorem A]{barr1980} one knows that any atomic Grothendieck topos is a boolean topos.

As was remarked by \cite{freyd1980}, to prove that $\sB_y$ provides a model $\mathrm{ZF}_y$ of Zermelo-Fraenkel set theory, one uses \cite{fourman1980}. By \cite{godel1939} one knows that any model, $\mathrm{ZF}_y$, of Zermelo-Fraenkel set theory has an inner model, $\mathrm{ZFC}_y$, of Zermelo-Fraenkel set theory with the Axiom of Choice.

This completes the proof of my claim. Hence, one sees that Mochizuki's Key Principle of Inter-Universality \iutf{\S I3, Pages 24-26}, does lead to models of set theory indexed by geometric base-points (in the sense of \eqref{eq:geom-base-point}), as claimed by Mochizuki.
\ep

\brem\label{re:zfc}\ 
\benumlab
\item Universes or different models of set theories are not used in the main body \iutf{Pages 40--701} of Mochizuki's paper. [To assert that these models are distinct, Mochizuki's Key Principle \iutf{IUT1, \ssep I3, Pages 24-27} suggests using the geometric base-point information (codified here in the datum $y$).] \item Mochizuki asserts in \iutf{Remark 3.11.1} that in \iutf{Theorem 3.11, Page 573} one compares $\theta$-values arising from one arithmetic holomorphic structure with respect to those arising from another arithmetic holomorphic structure. \item A careful reading of Mochizuki's discussion of his Key Principle  \iutf{\S I3, Pages 24--26} and how the geometric base-point information is used in \iut\ (for example \cite[Page 149]{mochizuki-gaussian} and other places therein) makes it clear that  \iut\ requires one to treat the fields (e.g. $\Q_p,\bQ_p$) provided by each geometric base-points separately.
\item It is precisely because of {\bf(1)-(3)} that my approach eschews all these nice meta-mathematical considerations and works directly by precisely quantifying the notion of local (resp. global)  arithmetic holomorphic structures--by \cite[Definition 4.1.1]{joshi-teich}, this is the arithmetic-geometric datum given by each $y$ for each prime $p$, (resp. \cite[\ssep 6]{joshi-teich-def}, \cite[\ssep 9]{joshi-teich}), and shows that this approach gives rise to a canonical Arithmetic Teichm\"uller Theory satisfying (the precise versions of) all the global (and local) properties needed for  the proof of  the $abc$-conjecture asserted by Mochizuki.
\item As pointed out above \eqref{eq:geom-base-point2},   my datum $y$ \eqref{eq:geom-base-point} also arises from $X^\Diamond(\cpt)$ (\cite[Section 15]{scholze12-perfectoid-ihes}). So anyone suggesting, as \cite{boyd2025} does, that Mochizuki's approach is artificial, or my approach to construction of Arithmetic Teichm\"uller Theory is artificial, stems from a poor understanding of underlying mathematics.
\eenum
\erem 

\section{Appendix: Existence of algebraically closed perfectoid fields}\label{pa:existence-perfectoid}
In the early days of my investigations of arithmetic Teichm\"uller Theory (\cite{joshi-untilts-2020}, \cite{joshi-teich}), my ideas were shaped by my reading of \cite{schmidt33} and \cite[Theorem 8]{kaplansky42}.
Eventually I replaced \cite{schmidt33} by \cite{kedlaya18} which is rather fundamental to my theory as it asserts the existence of topologically non-isomorphic perfectoid fields with tilts isometric with $\cpt$.  Nevertheless, an understanding of \cite{schmidt33} will still be useful to many readers of my theory. 
A more natural way of thinking about the theorem below is to think in terms of the absolute Riemann-Zariski Space $RZ(K,\Z)$  of an algebraically closed, complete valued field $K$, but I will take the more prosaic approach of \cite{schmidt33}.

\subsection{}
 The following is an essentially refined (but less general) version of \cite[Satz 2]{schmidt33}. The proof given here is similar to the one given there:
\bthm\label{th:schmidt}
Fix a prime number $p$. Let $(K,\abs{-}_K)$ be an algebraically closed, complete, rank-one valued field equipped  and suppose that $(\Q,\abs{-}_K)$ is equivalent to the standard $p$-adic valuation on $\Q$ (for example $K=\C_p$). Then there exist uncountably many rank-one valuations, $(K,\abs{-}_v)$ on $K$, with respect to each of which $K$ is complete, algebraically closed but each of these valuations $\abs{-}_v$ is not equivalent to the given valuation $(K,\abs{-}_K)$, while its restriction $(\Q,\abs{-}_v)$ is equivalent to the standard $p$-adic valuation on $\Q$. [In particular each $(K,\abs{-}_v)$ is an algebraically closed, perfectoid field with an isometric embedding of $\Q_p$.]
\ethm
\bp
Since $K$ is complete with respect to a non-trivial valuation, so $K$ is an uncountable field. The hypothesis implies that $K\supset\Q$. Choose a transcendence basis, $T=\{ t_\alpha\}_\alpha$, for the extension $K/\Q$. By replacing each $t_\alpha\in T$, by a multiple of the form $p^{c_\alpha}\cdot t_\alpha$  for a suitable $c_\alpha\in\Z$, one can assume additionally that $\abs{t_\alpha}_K\neq1$ (such a scaling evidently still provides a transcendence basis). Then $K\supset \Q(T)\supset \Q$ with $K\supset \Q(T)$ an algebraic extension while $\Q(T)\supset \Q$ is a purely transcendental extension. 

Since any valuation on a field extends to any of its algebraic extensions, it suffices to construct valuation $\abs{-}_v$ of $\Q(T)\supset \Q$ and extend it to a valuation $\abs{-}_v$ of $K$.

Let $S\subset T$ be any non-empty subset. Clearly, for any choice of $\emptyset\neq S\subset T$, there exist  a choice of integers $\{n_s\in\Z: s\in S\}$ such that for at least one $s\in S$, one has $n_s<0$.

Define an automorphism $\sigma:\Q(T)\to\Q(T)$ by $$\sigma(s)= s^{n_s}\qquad\forall s\in S,$$ and 
$$\sigma(t)=t \qquad \text{whenever } t\in T-S.$$  

Let $\abs{-}_v=\abs{-}_{v_{\sigma,S}}:\Q(T)\to \R$ be the pull-back, via $\sigma$, of the valuation $\abs{-}_K$ of $\Q(T)$. By construction, for any $x\in \Q(T)$: $$\abs{x}_v=\abs{\sigma(x)}_{K}.$$ 
From this one has
$$
\abs{x\cdot y}_v=\abs{\sigma(x\cdot y)}_{K}=\abs{\sigma(x)}_K\cdot\abs{\sigma(y)}_{K}=\abs{x}_v\cdot\abs{y}_v,
$$
and
$$
\abs{x+ y}_v=\abs{\sigma(x+y)}_{K}=\abs{\sigma(x)+\sigma(y)}_{K}\leq \max(\abs{\sigma(x)}_K,\abs{\sigma(y)}_K)=\max(\abs{x}_v,\abs{y}_v).
$$
So $\abs{-}_v$ is indeed a valuation on $K$. Moreover, since $\sigma$ is a $\Q$-automorphism of $\Q(T)$, $\sigma$ is identity on $\Q$, so for every $x\in \Q$, one has $$\abs{x}_v=\abs{\sigma(x)}_{K}=\abs{x}_{K}.$$

So the restriction of $\abs{-}_v$ to $\Q$ is equivalent to the standard $p$-adic valuation. Also one has for each $s\in S$: $$\abs{s}_v=\abs{\sigma(s^{n_s})}_{\Q(T)}=\abs{s}_{\Q(T)}^{n_s}.$$
Now $\abs{s}_K\neq1$ for all $s\in S$ and $n_s<0$  for at least one $s\in S$. Thus, for some $s\in S$,  one has either $\abs{s}_v<1$ and $\abs{s}_K>1$ or vice versa. Hence by \cite[Chap XII, Propostion 1.1]{lang-algebra}, the two valuations $\abs{-}_v$ and $\abs{-}_K$ are not equivalent on $\Q(T)$ and hence on $K$.   The completion of $(K,\abs{-}_v)$ has the same cardinality and the same characteristic as $K$ and both are algebraically closed, and hence the completion can be identified with the field $K$ and with this identification $K$ is algebraically closed, complete  with respect to both $\abs{-}_v$ and $\abs{-}_K$ and these two valuations on $K$ are inequivalent. At any rate, each $v$ provides a complete, algebraically closed field equipped with an isometric embedding of $\Q_p$. That such  algebraically closed, fields are perfectoid is immediate from Lemma~\ref{lem:perfectoid}. This proves the assertion. 
\ep

\newcommand{\cptmax}{{\C_p^{\flat\,\rm max}}}
\section{Appendix: Self similarity of $\cpt$ and its consequences}\label{se:self-similarity}
An important property of classical Teichm\"uller (and Riemann's moduli) Theory is the existence of fundamental domains. To take the genus one example, the upper half-plane $\fH\subset \C$ is tiled with copies of a fundamental domain under the  action of the (Teichm\"uller) modular group $\PSL_2(\Z)$. One may think of this as a self-similar tiling of the upper half-plane. One may ask if a similar self-similarity phenomenon exists in the Arithmetic Teichm\"uller Theory of this paper, and the answer is yes, and this self-similarity phenomenon arises from self-similarity of $\cpt$. This point is elaborated in this appendix.

\para Let me begin with the following  reformulation of an important result of \cite[Th\'eor\`eme 2 and \S3 Remarque 2]{matignon84} (and also \cite{kedlaya18}).
\bthm\label{th:self-similarity-of-cpt} Let $p$ be any prime number, $\C_p=\widehat{\overline{\Q}}_p$ be the completion of an algebraic closure of $\Q_p$. Let $\cpt$ be the tilt of $\C_p$. Then
\benumlab
\item There exists an isomorphism $\cpt \isom \fpc{x}$.
\item There exists $y\in \fpc{x}$ such that $\fpc{y}\subsetneqq \fpc{x}$ (more precisely $x\not\in\fpc{y}$). 
\item Hence $\fpc{x}$ is a self-similar valued field.
\item Hence $\cpt$ is a self-similar valued field i.e. it contains infinitely many proper subfields which are topologically isomorphic to $\fpc{t}\isom\cpt$ for some variable $t$.
\eenum
\ethm
\bp 
The first assertion is standard (see \cite{scholze12-perfectoid-ihes}). The second assertion is a consequence of the main theorem of \cite[Th\'eor\`eme 2 and \S3 Remarque 2]{matignon84} or the proof of \cite[Theorem 1.2]{kedlaya18}. The remaining assertions are immediate from the first two.
\ep

\para By a labeled copy of $\cpt$ will mean an identification $\cpt\isom \fpc{t}$ for some variable $t$. I will write $\cptl{t}$ for a copy of $\cpt$ labeled by the variable $t$.

\para A fundamental consequence of Theorem~\ref{th:self-similarity-of-cpt} is the following:
\bthm\label{th:self-similarity-of-ff-curve} 
The Fargues-Fontaine curve $\syfqp{}$ is a self-similar curve. More precisely, for every pair of elements $x,y\in\cpt$ as in Theorem~\ref{th:self-similarity-of-cpt}, there exists infinitely many strict inclusions 
$$\abs{\syfqp{y}}\into \abs{\syfqp{x}}$$
arising from the strict isometric inclusions $\cptl{y}\into\cptl{x}$.
\ethm
\bp 
Since $\abs{\syfqp{}}$ is identified by \cite[Th\'eor\`eme 2.4.1 and Corollaire 2.4.2]{fargues-fontaine} with the set of primitive degree one elements of $W(\O_{\cpt})$. By \cite[Corollaire 2.2.9]{fargues-fontaine} any primitive element of degree one in $W(\O_{\cptl{y}})$ can be written, up to multiplication by a unit in $W(\O_{\cptl{y}})$, as
$$[\alpha]-p,$$
for some element $\alpha\in\O_{\cptl{y}}$ with $v(\alpha)>0$   
and $\syfqp{y}\into \syfqp{x}$ is given  by sending the primitive element $[\alpha]-p\in W(\O_{\cptl{y}})$ to the primitive degree one element $[\alpha]-p\in W(\O_{\cptl{x}})$ and at the level of ideals $$([\alpha]-p)W(\O_{\cptl{y}}) \mapsto ([\alpha]-p)W(\O_{\cptl{x}}).$$
One has
$$W(\O_{\cptl{y}})\subsetneq W(\O_{\cptl{x}}),$$ as $x\notin \cptl{y}$ so $[x]-p\not\in W(\O_{\cptl{y}})$ and hence there is a primitive element of degree one of $W(\O_{\cptl{x}})$ which is not contained in the set of primitive elements of degree one of $W(\O_{\cptl{y}})$. So the inclusion of $\abs{\syfqp{y}}\into \abs{\syfqp{x}}$ is strict.  This proves the assertion.
\ep

\para Before proceeding it will be useful to understand this self-similarity in terms of  Classical Teichm\"uller Theory. In the classical Teichm\"uller Theory (i.e Teichm\"uller Theory for Riemann Surfaces), the Teichm\"uller space is tiled by isomorphs of a fundamental domain for the mapping class group or modular group actions. To put it differently the Teichm\"uller space is equipped with a self-similar tiling (not unique in general).

\para So the question arises if there is a theory of fundamental domains in the Arithmetic Teichuller Theory constructed here. The answer to this question is yes and is given by \Cref{th:main6-fundamental-domain}. The notion of fundamental domains in Arithmetic Teichm\"uller Theory and this arises precisely from the fact that  that $\cpt$ is a self-similar valued field i.e. $\cpt$ (see Theorem~\ref{th:self-similarity-of-cpt} and Theorem~\ref{th:self-similarity-of-ff-curve}).

\para\label{pa:fja-defined} Let $E$ be a $p$-adic field and let $X/E$ be a smooth, geometrically connected, smooth, quasi-projective variety. Let $\fja{X,E}$ be the category of $\fJ(X,E)$ defined by \Cref{def:subcat-F} for $F=\cpt$ i.e.
\be
\fja{X,E}=\fjxe_{\cpt}
\ee

\para A \emph{fundamental domain for $\fja{X,E}$} is the full subcategory $\fja{X,E}_{\cptl{x}}$ for some variable $x$. The following is now a tautology:

\bthm\label{th:main6-fundamental-domain}
Let $X$ be a geometrically connected, smooth, quasi-projective variety over a $p$-adic field $E$.
\benumlab
\item
Any object of $\fja{X/E}$ belongs to some  fundamental domain $\fJ(X,E)_{\cptl{t}}$.
\item The identification $y=x$ (of variables) provides a tautological equivalence of categories
$$\fja{X,E}_{\cptl{y}}\isom \fja{X,E}_{\cptl{x}}.$$
\item In particular $\fja{X,E}$ is tiled with copies of the fundamental domains $\fJ(X,E)_{\cptl{x}}$.
\eenum
\ethm

\para The self-similarity of $\cpt$ discussed in Section~\ref{se:self-similarity} at once implies that the group $\sG(\O_{\cpt})$ is self-similar and this in turn has the following important consequence whose proof is clear from Theorem~\ref{th:self-similarity-of-cpt} and the preceding definitions and discussion:
\bthm\label{th:self-similarity-of-theta-values}
Assume $F=\cpt$. Then
\benumlab
\item $\ttxl\subset \bpip\isom \sG(\O_F)$.
\item The set $\ttxl$ is a self-similar subset of $\bpip$: more precisely let 
$$\ttxlt{t}=\left\{z\in\sG(\O_F): \eta_{K}(z)\in f_{\theta}(C[\ell]-\{O\}) \subset K \text{ for some perfectoid field with } K^\flat=\cptl{t}\right\}$$
Then $$\ttxl=\bigcup_{\cptl{t}} \ttxlt{t}$$
where the union runs over all the (isometric) identifications $\cpt=\cptl{t}$. 
\item To put it colloquially, the $\theta$-torsion value sets $\ttxl$ form a fractal in $\bpip\isom\sG(\O_F)$.
\eenum
\ethm

\bp 
After Theorem~\ref{th:self-similarity-of-cpt} it is enough to note that any $z\in\ttxl$ lives in some $\ttxlt{t}$ and so the assertion is immediate.
\ep

\newcommand{\tK}{\widetilde{K}}
\newcommand{\tO}{\widetilde{\O}}

\section{Appendix: Anabelian Hyperbolic Varieties}\label{se:number-field-case}
\emph{\textcolor{red}{This section may be skipped on the initial reading. Reader may simply work with hyperbolic curves instead.} The discussion of hyperbolic varieties of dimension bigger than one is included only to illustrate my point that the construction provided in this paper works in higher dimensions.} 

\para Let $E\supset \Q$ be any complete valued field, which is either an archimedean  or a non-archimedean with a rank one valuation inducing a $p$-adic valuation on $\Q$, and $\ebh$ be its completed algebraic closure. Let $X/E$ be a geometrically connected, smooth projective variety over $E$. I will say that $X/E$ is a \emph{hyperbolic variety} if the analytic space $X^{an}/\ebh$, is a (Brody) hyperbolic variety (see  \cite{lang86} for the archimedean case, \cite{javanpeykar18} for the non-archimedean case). 

\para If $\dim(X)=1$ then $X/E$ is hyperbolic (in the above sense) if and only if $X\times_E\C$ is a hyperbolic Riemann surface. 

\para In \cite{lang86}, Serge Lang has conjectured that any hyperbolic variety $X$ defined over a number field $E$ has a finite number of $E$-rational points.

\brem\ 
\benumlab
\item \emph{For concrete applications, beyond $dim(X)=1$, one may also need to assume that $X/E$ is a $K(\pi,1)$-space in some suitable sense.}
\item As an aside let me say that I have used the case of \emph{hyperbolic varieties} here because of their relevance in Diophantine Geometry via Lang's Conjecture \cite{lang86} and \cite{faltings91}. One can also use other related hypothesis as a substitute for the hyperbolic hypothesis, for example, one can work with \emph{groupless varieties} instead of hyperbolic varieties (see \cite{javanpeykar20} for other related hypothesis which may be used here instead).
\eenum
\erem
\para An important property in many anabelian considerations is the following: a profinite group $\Pi$ is said to be a \emph{slim profinite group} (or simply  $\Pi$ is \emph{slim}) if every open subgroup of $\Pi$ has trivial center. By \cite[Def. 0.1 and Remark 0.1.3]{mochizuki04}, $\Pi$ is slim if and only if the centralizer of any open subgroup of $\Pi$ is trivial.

\para  A hyperbolic variety over a finitely generated field has the following anabelian property.
\bpro
Let $X/E$ be a geometrically connected, smooth, hyperbolic variety over a finitely generated field $E$ of characteristic zero or a $p$-adic field. Then there exists a basis of Zariski open subsets $\{U\}$ of $X$ such that
\benumlab 
\item every $U\not=\emptyset$ in this basis is hyperbolic, and
\item for every pair $U,V$ nonempty opens in this basis one has
$${\rm Isom}_E(U,V)\isom {\rm Isom}_{G_E}^{out}(\pi_1(U/E),\pi_1(V/E)),$$
\item and for every $U\neq\emptyset$, $\pi_1(U/E)$ is slim.
\eenum
\epro
\bp 
From \cite[Corollary 1.7]{stix16}, if  $E$ is finitely generated field and by \cite[Theorem C]{hoshi-poly}, if $E$ is a $p$-adic field, it follows that a basis of Zariski open sets satisfying property (2) exist. From \cite{lang86} it is immediate that every non-empty open subvariety of $X$ is hyperbolic. The last assertion is proved by induction: for $\dim(U)=1$, $U$ is a hyperbolic curve and so the slim-ness property is immediate from the aforementioned results. For $\dim(U)>1$ one uses induction on dimension and the proofs of \cite{stix16} or \cite{hoshi-poly} from the fibration structure $U$ is equipped with by the construction of this basis.
\ep

\begin{defn}\label{def:anab-hyp}
Let $X/E$ be an hyperbolic variety over a $p$-adic field or a finitely generated field $E$. Let $\bE$ be an algebraic closure of $E$. Then $X$ is said to be an \textit{anabelian hyperbolic variety} if the \'etale fundamental group $\pi_1(X/E)$ and	the geometric tempered \'etale fundamental group $\pi_1(X/\bE)$ are both slim.
\end{defn}

\brem
\emph{There exist   hyperbolic varieties which are not anabelian hyperbolic non-slim (see  \cite{ihara97}).} 
As an addendum to \cite{ihara97}, which will not be used in the rest of this paper, let me remark that if $A/E$ is a simple abelian variety over a number field $E$, and suppose $D\subset A$ over $E$ is an irreducible, smooth, ample divisor, then by \cite[Theorem 4.1]{debarre95}, one has an isomorphism $\pi_1(D\times \bE)\isom \pi_1(X\times\bE)$ of the geometric fundamental groups. In particular $\pi_1(D\times\bE)$  is not slim. But as $A$ is simple, $D$ is not a translate of an a sub-abelian variety of $A$, and hence by \cite{lang86}, $D$ is hyperbolic. \textit{Hence $D$ is not an anabelian hyperbolic variety.} So this provides another class of examples of non-slim hyperbolic varieties complementing the non-slim examples of \cite{ihara97}. Note that by \cite[Theorem 1]{faltings91}, $D$ has finitely many rational points. 
\erem

\bpro\label{pa:hyp-curve-slim-ness} If $X/E$ is a hyperbolic curve over a $p$-adic field or a finitely generated field $E$ then $X/E$ is an anabelian hyperbolic variety.
\epro
\bp 
This is immediate from \Cref{def:anab-hyp} and \cite[Corollary 1.3.3 and Lemma 1.3.1]{mochizuki04}.
\ep

\para The absolute Grothendieck Conjecture for anabelian hyperbolic varieties is the following:
\bcon\label{con:anab-hyper}
Let $X/L$ be any smooth, quasi-projective and anabelian hyperbolic variety over a finitely generated field  or  a $p$-adic field $L$.  If $Y/L'$ is a smooth, quasi-projective anabelian hyperbolic variety which is anabelomorphic to $X/L$ the one has an isomorphism $Y\isom X$ of $\Z$-schemes.
\econ

\brem\  
\benumlab
\item This is a formulation, for anabelian hyperbolic varieties, of  the celebrated anabelian conjecture of Grothendieck. 
\item For $\dim(X)=1$, \Cref{con:anab-hyper} is a well-known theorem of  \cite{mochizuki96-gconj}, \cite{tamagawa97-gconj}. 
\eenum
\erem

\section{Appendix: Anabelian variations providing $\Pi$}\label{se:comments}
\numberwithin{equation}{subsection}

\setlist[enumerate]{wide, labelwidth=\parindent}
\setlist[enumerate,1]{labelindent=0pt, labelwidth=\parindent, label=\arabic*.}
\setlist[enumerate,2]{wide=.825cm, label={{\bf(\alph{*})}}}
\setlist[enumerate,3]{wide=1.25cm, label=(\arabic*)}

\para I want to present some elementary considerations which will prove useful in understanding the problem of constructing Arithmetic Teichm\"uller Spaces using the tempered fundamental group. 
\newcommand{\sC}{\mathscr{C}}
\newcommand{\schz}{\mathcal{Sch}_\Z}
\newcommand{\prodis}{\mathcal{ProD}_{\Pi}}
\newcommand{\prof}{\mathcal{ProF}_{\widehat{\Pi}}}
\newcommand{\qpan}{\mathcal{An}_{\Q_p}}
\newcommand{\anc}{{\mathcal{An}_{\C}}}
\para Fix a pro-discrete group $\Pi$ (for example $\Pi=\pit{X/E}$). Let $\prodis$ be the isomorphism class of the pro-discrete group $\Pi$, 
and for each pro-discrete group $H\in\prodis$, let  $\widehat{H}$ be its profinite completion. One can make $\prodis$ into a category with the class $\prodis$ as objects and morphisms between any pair of objects are isomorphisms. This makes the category $\prodis$ into a groupoid whose every object is isomorphic to $\Pi$. Let $\schz$ be the category of schemes.  Let $\qpan$ (resp. $\anc$) be the category of $\Q_p$-analytic spaces  in the sense of Theorem~\ref{thm:main} (resp. $\C$-analytic spaces).

\para  Let  $S$ be one of categories $\schz,\anc,\qpan$. Suppose $\sC$ is a category and $\Pi$ is a fixed pro-discrete group. I will say that  \emph{$\sC$ is an ananabelian variation providing $\Pi$ with base $S$} if the following conditions are satisfied:
\benumlab
\item  For every $V$ in $\sC$ there exists an isomorphism of pro-discrete groups $\alpha_V:\Pi_V\mapright{\isom} \Pi$ i.e.  one is given a functor $${\rm ob}(\sC) \to \prodis$$  written $V\mapsto \Pi_V$ from the class of objects of $\sC$ to the category $\prodis$ of $\Pi$.
\item There is a  functor $\sC\to S$ denoted $V\mapsto [V]\in S$. In this case  $[V]$ called the scheme (resp. $\C$-analytic space, $\Q_p$-analytic space) underlying $V\in\sC$.
\eenum

\para I will often simply say ``anabelian variation providing $\Pi$'' instead of ``anabelian variation providing $\Pi$ with base $S$.'' Hopefully there will be no confusion. The categories constructed here will come with functors to all the three values of $S$.

\para I will say that an anabelian variation providing $\Pi$ is a \emph{trivial anabelian variation providing $\Pi$} if any pair of  objects $V,V'\in\sC$ are isomorphic. Obviously one is interested in constructing non-trivial anabelian variations providing $\Pi$. 

\para Suppose $\sC$ is a non-trivial anabelian variation providing $\Pi$. Then the function $V\mapsto \Pi_V$ can be thought of as  providing labeled isomorphs of $\Pi$ and the function ${\rm ob}(\sC)\to \prodis$ ($V\mapsto \Pi_V$) will be called the \emph{labeling function of $\sC$}.

\para\label{pa:classical-teichmuller-space} \emph{The classical Teichm\"uller Space of any connected, hyperbolic Riemann surface  $\Sigma$ is a non-trivial anabelian variation providing $\Pi$, the topological fundamental group of $X$.} Indeed suppose $\Sigma$ is a connected, hyperbolic Riemann surface with $\Pi=\pi_1(\Sigma)$ be its \'etale fundamental group and consider the  Teichm\"uller space $T_\Sigma$ of $\Sigma$ \cite[Chapter V]{lehto-book}. Indeed one may think of $T_\Sigma$ as a category consisting of pairs $(\Sigma, f:\Sigma\to\Sigma')$ where $f$ is a quasi-conformal mapping of $\Sigma$ onto a Riemann surface $\Sigma'$. Then one has $(\Sigma, f)\mapsto \pi_1(\Sigma')\isom \pi_1(\Sigma)=\Pi$; and as any connected hyperbolic Riemann surface is obviously a $\C$-analytic space so one has an obvious functor to $T_\Sigma\to\anc$ given by $(\Sigma,f:\Sigma\to \Sigma')\mapsto \Sigma'$.  Hence $T_{\Sigma}$ is an anabelian variation providing $\Pi$ with base $\anc$.

\para\label{pa:geometric-anabelian-var} I will say that an anabelian variation providing $\Pi$ with base $S$ is a \emph{geometric anabelian variation providing $\Pi$ with base $S$} if there exists $V,V'\in\sC$ such that $[V]\neq[V']$ (in $S$).

\para Obviously any geometric anabelian variation providing $\Pi$ is non-trivial. Hence if one constructs geometric anabelian variations providing $\Pi$ then one automatically gets a non-trivial anabelian variation providing $\Pi$.

\para The moduli $\sM_{g}/\C$ stack of smooth, proper curves of genus $g\geq 2$ over $\C$ is an example of a geometric anabelian variation providing $\Pi$. Indeed it is clear that this is an anabelian variation providing $\Pi=\pi_1^{top}(X(\C))$ is the topological fundamental group of the Riemann surface $X(\C)$ where $X/\C$ is any smooth curve, proper curve  of genus $g$. This is also a geometric anabelian variation providing $\Pi$ because one can obviously find two non-isomorphic smooth, proper curves of genus $g$ over $\C$. Such a pair of curves cannot be isomorphic as $\Z$-schemes as well. In particular $\sM_g$ is a non-trivial anabelian variation providing $\Pi$. 

\para These examples should convince the reader that, when a non-trivial (or even a geometric) anabelian variation providing $\Pi$ exists, then it can serve as an \emph{anabelian stand-in for a variation of (mixed) Hodge structures.} 

\para In some sense the presence of this sort of a structure (i.e. an anabelian variation providing $\Pi$ with base $S$) should be understood as a manifestation of Kodaira-Spencer classes! Constructing similar structures in  the $p$-adic setting leads to  an \emph{$p$-adic Teichm\"uller Landscape} or an \emph{$p$-adic Teichm\"uller Theory} presented here. Assembling such data for each valuation of a number field, leads to a global  \emph{Arithmetic Teichm\"uller Landscape} or an \emph{Arithmetic Teichm\"uller Theory} (also presented here)  in which one can hope to contemplate applications to global Diophantine problems as is done in \iut. 

\para Now the results of \Cref{se:construct-att} and the preceding lead to the following:
\bthm\label{th:anab-var} 
Let $X/E$ be a geometrically connected, smooth, quasi-projective variety over a $p$-adic field $E$ or let $E=\C$. If $E$ is a $p$-adic field, let $S=\ufjxe$ and if $E=\C$, let $S$ be the category of complex Banach spaces with isometries as morphisms. If $E$ is a $p$-adic field, let $\Pi=\pit{X/E}$ be the tempered fundamental group and if $E=\C$ let $\Pi$ be the topological fundamental group.

\benumlab
\item If $E$ is a $p$-adic field,  then one has a natural functor (given by forgetting the geometric base-point)
$$\fjxe\to S=\ufjxe, $$
and  the tempered fundamental group functor:
$$\fjxe\to \prodis.$$
\item If $E=\C$ one has the Royden algebra functor $\fjxe\to S$ and the topological fundamental group functor
$$\fjxe\to \prodis.$$
\eenum
In both the cases, these two functors make $\fjxe$ into a non-trivial anabelian variation providing $\Pi$.
\ethm
\bp 
The proof is clear from the results of \Cref{se:construct-att}.
\ep

	\bibliographystyle{dinat}

\Address
\end{document}